\documentclass[preprint,12pt]{elsarticle}

\usepackage[letterpaper, margin = 2.5cm]{geometry} 

\usepackage{amssymb}
\usepackage{amsmath}
\usepackage{lipsum}
\usepackage{amsfonts}
\usepackage{graphicx,subfigure}
\usepackage{epstopdf}
\usepackage{algorithmic}
\usepackage{booktabs}
\usepackage{multirow} 
\usepackage{enumitem}
\usepackage{stmaryrd}

\usepackage{hyperref}
\usepackage{xcolor}
\hypersetup{
    colorlinks=true,
    citecolor=green, 
    linkcolor=red,   
    urlcolor=cyan    
}
\newcommand{\ie}{\emph{i.e.} }

\newcommand{\C}{\mathbb{C}}

\newcommand{\cu}{\boldsymbol}

\newtheorem{assumption}{Assumption}[section]
\newtheorem{theorem}{Theorem}
\newproof{lemma}{Lemma}
\newtheorem{remark}{Remark}

\makeatletter
\def\ps@pprintTitle{%
 \let\@oddhead\@empty
 \let\@evenhead\@empty
 \let\@oddfoot\@empty
 \let\@evenfoot\@empty
}
\makeatother

\begin{document}

\begin{frontmatter}

\title{A Hybrid Iterative Neural Solver Based on Spectral  Analysis for Parametric PDEs}

\author[a]{Chen Cui}
\ead{cuichensx@gmail.com}
\author[a]{Kai Jiang}
\author[a]{Yun Liu}
\author[a]{Shi Shu}

\affiliation[a]{organization={Hunan Key Laboratory for Computation and Simulation in Science and Engineering,
Key Laboratory of Intelligent Computing and Information Processing of Ministry
of Education, School of Mathematics and Computational Science, Xiangtan University},
            city={Xiangtan},
            postcode={411105}, 
            state={Hunan},
            country={China}}

\begin{abstract}
     Deep learning-based hybrid iterative methods (DL-HIM) have emerged as a promising approach for designing fast neural solvers to tackle large-scale sparse linear systems. DL-HIM combine the smoothing effect of simple iterative methods with the spectral bias of neural networks, which allows them to effectively eliminate both high-frequency and low-frequency error components. However, their efficiency may decrease if simple iterative methods can not provide effective smoothing, making it difficult for the neural network to learn mid-frequency and high-frequency components. This paper first conducts a convergence analysis for general DL-HIM from a spectral viewpoint, concluding that under reasonable assumptions, DL-HIM exhibit a convergence rate independent of grid size $h$ and physical parameters $\cu{\mu}$. To meet these assumptions, we design a neural network from an eigen perspective, focusing on learning the eigenvalues and eigenvectors corresponding to error components that simple iterative methods struggle to eliminate. Specifically, the eigenvalues are learned by a meta subnet, while the eigenvectors are approximated using Fourier modes with a transition matrix provided by another meta subnet. The resulting DL-HIM, termed the Fourier Neural Solver (FNS), can be trained to achieve a convergence rate independent of PDE parameters and grid size within a local neighborhood of the training scale by designing a loss function that ensures the neural network complements the smoothing effect of the damped Jacobi iterative methods. We verify the performance of FNS on five types of linear parametric PDEs.
 \end{abstract}

\begin{keyword}
     Hybrid iterative method\sep
     Neural solver\sep
     Preconditioning\sep
     Convergence analysis\sep 
     Spectral bias



\end{keyword}

\end{frontmatter}



\section{Introduction}\label{sec:01}
Large-scale sparse linear system of equations
\begin{equation}
    \cu{A} \cu{u}=\cu{f},
    \label{eq:linear_system}
\end{equation}
are ubiquitous in scientific and engineering applications, particularly those arising from the discretization of partial differential equations (PDEs). Developing efficient, robust, and scalable numerical methods to solve these equations remains a significant challenge for researchers in applied mathematics.
Iterative methods \cite{saad2003iterative} are effective methods for solving such systems. Starting with an initial guess $\cu{u}^{(0)}$, the simplest form is relaxation
\begin{equation}
    \cu{u}^{(k+1)} = \cu{u}^{(k)} + \omega \left(\cu{f} - \cu{A u}^{(k)}\right),
    \label{eq:relaxation}
\end{equation}
where $\omega$ is the relaxation parameter. This method is known as the Richardson iterative method. Although the computational cost per step is cheap, its convergence rate is quite slow in practical applications. Consequently, various acceleration techniques have been developed.
For example, one can select the optimal $\omega$ using Chebyshev polynomials \cite{golub1961chebyshev}; search for solutions along the direction of conjugate gradient instead of the negative gradient \cite{hestenes1952methods}; introduce the momentum term, such as Nesterov acceleration \cite{nesterov1983method,luo2022differential}; split the unknown $\cu{u}$ into blocks for block coordinate descent \cite{beck2013convergence}; split $\cu{A}$ into the sum of different operators for alternating-direction implicit iteration \cite{birkhoff1962alternating}; use information from previous steps, such as Anderson acceleration \cite{walker2011anderson}; and project into the Krylov subspace to find the solution, such as GMRES \cite{saad1986gmres}, etc.

However, the convergence speed of these acceleration techniques is still constrained by the condition number of $\cu{A}$. To mitigate this limitation, preconditioning techniques are introduced
\begin{equation}
    \cu{u}^{(k+1)} = \cu{u}^{(k)} + \cu{B}\left(\cu{f} - \cu{A u}^{(k)}\right),
    \label{eq:prec}
\end{equation}
where $\cu{B}$ is called the preconditioner, designed to approximate $\cu{A^{-1}}$ and should be computationally easy to obtain.
Simple preconditioners, such as damped Jacobi, Gauss-Seidel, and successive over-relaxation methods \cite{saad2003iterative}, also exhibit slow convergence rates. Commonly used preconditioners include incomplete LU (ILU) factorization \cite{chow2015fine}, multigrid (MG) \cite{trottenberg2000multigrid}, and domain decomposition methods (DDM) \cite{toselli2004domain}. 
For example, MG methods are optimal for elliptic equations. In more complex cases, preconditioners often need to be combined with acceleration methods, such as flexible conjugate gradient (FCG) \cite{notay2000flexible} and flexible GMRES (FGMRES) \cite{saad1993flexible}. However, when tackling challenging problems, effective methods frequently require problem-specific parameters that must be determined by experts.

In recent years, deep learning techniques have emerged as innovative approaches to solving PDEs, offering new perspectives in scientific computing. There are three main applications:  
First, as a universal approximator for solving complex (\emph{e.g.}, high-dimensional) PDEs, known as \textit{neural pde} \cite{cuomo2022scientific}.  
Second, as a discretization-invariant surrogate model for parametric PDEs (PPDE) that maps infinite-dimensional parameter spaces to solution spaces, referred to as \textit{neural operator} \cite{kovachki2023neural}.  
Third, in designing fast iterative methods for discretized systems, which we refer to as \textit{neural solver}.
Neural solvers primarily evolve from two aspects. The first involves automatically learning problem-specific parameters for existing iterative methods. For example, in acceleration techniques, neural networks have been utilized to learn parameters in Chebyshev acceleration for improved smoothing effects \cite{cui2024neural}, or to learn iterative directions in place of conjugate gradient \cite{kaneda2022deep}.
For parameters in preconditioners, neural networks can be used to correct ILU preconditioners \cite{trifonov2024learning}; learn smoothers \cite{katrutsa2020black, huang2022learning, chen2022meta}, transfer operators \cite{greenfeld2019learning, luz2020learning, kopanivcakova2024deeponet}, coarsening \cite{taghibakhshi2021optimization,caldana2023deep, zou2023autoamg} in MG; learn adaptive coarse basis functions \cite{klawonn2024learning},
interface conditions and interpolation operators \cite{taghibakhshi2022learning, knoke2023domain,taghibakhshi2023mg, kopanivcakova2024deeponet} in DDM, etc.

The other approach involves combining the traditional preconditioner $\cu{B}$ with the neural network $\mathcal{H}$ to form the following deep learning-based hybrid iterative method (DL-HIM)
\begin{subequations}
    \begin{align}
    \cu{u}^{(k+\frac{1}{2})} &= \cu{u}^{(k)} + \cu{B} (\cu{f} - \cu{A u}^{(k)}) \quad \text{(smoothing iteration, repeat } M \text{ times)}, \label{eq:hybrid1} \\
    \cu{u}^{(k+1)} &= \cu{u}^{(k+\frac{1}{2})} + \mathcal{H} (\cu{f} - \cu{A u}^{(k+\frac{1}{2})}) \text{ (neural iteration)}. \label{eq:hybrid2}
    \end{align}
    \label{eq:hybrid}
\end{subequations}
The primary motivation is that $\mathcal{H}$ tends to fit low-frequency functions due to spectral bias \cite{rahaman2019spectral, hong2022activation, xu2019frequency}, while simple preconditioners $\cu{B}$, such as damped Jacobi and Gauss-Seidel, are effective for eliminating high-frequency error components. 
We refer to \eqref{eq:hybrid1} as the smoothing iteration because it is primarily intended to eliminate oscillatory error components. However,  it can not fully achieve this goal in certain problems, so the term ``smoothing iteration'' is used in a broader sense, as long as it can reduce some parts of the error.
The neural network $\mathcal{H}:\mathbb{C}^N \rightarrow \mathbb{C}^N$ need be applicable for any $N \in \mathbb{Z}^+$. Therefore, it is typically designed as a discretization-invariant neural operator.
The iterative method \eqref{eq:hybrid} was first proposed in \cite{zhang2022hybrid}, where DeepONet \cite{lu2019deeponet} was used as $\mathcal{H}$ and combined with different preconditioners $\cu{B}$ to test the compressibility  for error components of various frequencies. This approach, referred to as the hybrid iterative numerical transferable solver (HINTS), was later extended to handle different geometries \cite{kahana2023geometry} and to solve the Helmholtz equation \cite{zou2024large}. Around the same time, the authors of \cite{cui2022fourier} employed local Fourier analysis (LFA) \cite{brandt1977multi} to estimate which frequency components the selected $\cu{B}$ could effectively eliminate.
Inspired by the fast Poisson solver \cite{fortunato2020fast,strang2007computational}, they designed $\mathcal{H}$ based on the fast Fourier transform (FFT) to eliminate error components difficult for $\cu{B}$ to handle.
In 2023, the authors of \cite{xie2023mgcnn} designed a multilevel structure that mimics MG as $\mathcal{H}$, with shared parameters across different levels. This approach demonstrated good computational efficiency for convection-diffusion equations. In 2024, the authors of \cite{rudikov2024neural} used the SNO \cite{fanaskov2023spectral} as $\mathcal{H}$ and further accelerated it with FCG. In the same year, the authors of \cite{hu2024hybrid} used MIONet \cite{jin2022mionet}  as $\mathcal{H}$ and analyzed the convergence rate of the corresponding DL-HIM for the Poisson equation for the first time.
Also in 2024, the author of \cite{chen2024graph} used a graph neural network to directly approximate $\cu{A}^{-1}$ as a nonlinear preconditioner for FGMRES.

In this paper, we consider the linear systems arising from discretizing the steady-state linear PPDE
\begin{equation}
    \mathcal{L}[u(\boldsymbol{x}, \boldsymbol{\mu}); \boldsymbol{\mu}] = f(\boldsymbol{x}, \boldsymbol{\mu}), \quad (\boldsymbol{x}, \boldsymbol{\mu}) \in \Omega \times \mathcal{P},
    \label{ppde}
\end{equation}
where $\mathcal{L}$ is the differential operator, $u$ is the solution function, $f$ is the source term, and $\boldsymbol{\mu}$ is the parameter in a compact space.
It is assumed that the equation \eqref{ppde} is well-posed under appropriate boundary conditions.
By dividing $\Omega$ into a discrete grid $\mathcal{T}_h$ and applying numerical discretization methods, such as the finite element method (FEM), we obtain the following linear systems of equations
\begin{equation}\nonumber
    \cu{A_\mu} \cu{u_\mu} = \cu{f_\mu},
\end{equation}
where $\cu{A_\mu} \in \C^{N \times N}$ and $\cu{f_\mu} \in \C^{N}$.
In many practical scenarios, such as inverse problems, design, optimization, and uncertainty quantification, it is crucial to examine the behavior of the physics-based model across various parameters $\cu{\mu}$. In these cases, high-fidelity simulations, which involve solving the linear system multiple times, can become prohibitively expensive.
Therefore, we aim to utilize DL-HIM to reduce the computational cost. For simplicity in writing, we omit the superscripts related to specific parameters $\cu{\mu}$ and consider linear systems of the form \eqref{eq:linear_system}.

The main contents of this paper are as follows. 
First, we do a convergence analysis for the general DL-HIM \eqref{eq:hybrid} from a spectral viewpoint. For a diagonalizable matrix $\cu{A}$, under reasonable assumptions on $\cu{B}$ and $\mathcal{H}$, we can obtain a DL-HIM with a convergence rate independent of the mesh size $h$ and physical parameters $\cu{\mu}$.
To ensure that $\mathcal{H}$ meets its assumptions, we then design $\mathcal{H}$ from an eigen perspective. Building on our previous work \cite{cui2022fourier}, we improve $\mathcal{H}$ by replacing the use of a Fourier matrix as the eigenvector matrix with a transition matrix that maps from the Fourier modes to the eigenvector basis. By introducing two meta subnets to learn the inverse of the eigenvalues and the transition matrix, we transform the requirements on $\mathcal{H}$ into requirements on the meta subnets. This approach avoids the constraints of spectral bias, enabling $\mathcal{H}$ to approximate not only low-frequency error components but also error components that $\cu{B}$ struggles to eliminate.
The improved solver, still called Fourier Neural Solver (FNS), is implemented as an end-to-end, differentiable neural solver  for linear systems derived from structured grids.
Finally, we design a reasonable loss function and construct easily accessible and targeted training data to train FNS and test its performance on a variety of classical PDE discrete systems.
In the numerical experiments, we first verify the validity of $\mathcal{H}$ for the Poisson equation using analytic eigenvalues and eigenvectors. This results in an FNS where the convergence rate remains independent of the problem size.
For random diffusion equations, the selected $\cu{B}$ effectively eliminates high-frequency errors, while a well-trained $\mathcal{H}$ can learn low-frequency and mid-frequency error components. This enables FNS to achieve a convergence rate independent of both physical parameters and grid size within a local neighborhood of the training scale.
For multi-scale problems, such as discrete systems arising from anisotropic diffusion, convection-diffusion, and jumping diffusion equations, FNS maintains a convergence rate independent of physical parameters and grid size on medium scales.
In the case of the Helmholtz equation, even though the chosen $\cu{B}$ amplifies the lowest-frequency errors, the trained $\mathcal{H}$ effectively corrects these errors, ensuring a convergent FNS.
In summary, the advantages of FNS over other DL-HIM are as follows:
\begin{enumerate}[label=(\arabic*)]
\item FNS can handle error components with different frequencies. Several DL-HIMs \cite{zhang2022hybrid, kahana2023geometry, hu2024hybrid, zou2024large, rudikov2024neural} train $\mathcal{H}$ independently, following an approach similar to that used in neural operators. However, the trained $\mathcal{H}$ tends to capture only low-frequency components due to spectral bias. In contrast, FNS mitigates this issue by learning in the frequency domain and employing an end-to-end training strategy, which enables better performance.
\item The training data for FNS is easy to obtain as it does not require solving PDEs. The loss function is designed to ensure that $\mathcal{H}$ focuses on error components that $\cu{B}$ struggles to eliminate.
\item The nonlinearity of FNS only exists in the meta subnets; $\mathcal{H}$ is linear and therefore scale-equivariant. During the iteration process, the input (residual) of $\mathcal{H}$  can have vastly different scales, and the output (error) should also change proportionally. DL-HIM with nonlinear $\mathcal{H}$ \cite{hsieh2019learning, zhang2022hybrid, kahana2023geometry, zou2024large, rudikov2024neural, chen2024graph} cannot naturally guarantee this property.
\item The asymptotic convergence rate of FNS is independent of the right-hand sides (RHS).
\end{enumerate}
 
The rest of this paper is organized as follows: Section \ref{sec:02} establishes a convergence analysis framework for the general DL-HIM; Section \ref{sec:03} proposes FNS and designs reasonable training data and loss function under the guidance of convergence analysis; Section \ref{sec:04} conducts numerical experiments on several types of second-order linear PDEs to verify the theoretical analysis; and Section \ref{sec:final} gives a summary and outlook.

\section{Convergence analysis}\label{sec:02}
The error propagation matrix of \eqref{eq:hybrid} is given by
$$
\cu{E} = (\cu{I} - \mathcal{H} \cu{A})(\cu{I - BA})^M,
$$
then the iterative method \eqref{eq:hybrid} converges if and only if
$$
\rho(\cu{E}) < 1,
$$
where $\rho$ denotes the spectral radius. If $\rho(\cu{I - BA}) < 1$, \eqref{eq:hybrid} will certainly converge as $M \rightarrow \infty$ \cite{hu2024hybrid}. 
Next, we analyze the convergence of the smoother $\cu{B}$ and the neural operator $\mathcal{H}$ from a spectral perspective for a fixed $M$.

\subsection{Convergence of the smoother $\cu{B}$}

We use LFA \cite{brandt1977multi} to perform convergence analysis on the smoothing iteration \eqref{eq:hybrid1}. 
Define the infinite grid 
\begin{equation}
    \cu{G}_h=\{\cu{x=kh}:=(k_1h_1,k_2h_2),\mathrm{~}\cu{k}\in\mathbb{Z}^2\},
\end{equation}
and the Fourier modes on it
\begin{equation}
    \varphi(\cu{\theta,x})=e^{i\cu{\theta}\cdot \cu{x}/\cu{h}}:=e^{i\theta_1x_1/h_1}e^{i\theta_2x_2/h_2}\quad\mathrm{~for~}\cu{x}\in\cu{G}_h,
\end{equation}
where $\cu{\theta}=(\theta_1,\theta_2)\in\mathbb{R}^2$ denotes the Fourier frequencies, which can be restricted to $\Theta = [-\pi, \pi)^2 \subset \mathbb{R}^2$, due to the periodic property
$$\varphi(\cu{\theta}+2\pi, \cu{x}) = \varphi(\cu{\theta,x}).$$
Consider a general discrete operator defined on $\cu{G}_h$ 
\begin{equation}
L_hu_h(\cu{x})=\sum_{\cu{k} \in V}c_{\cu{k}}(\cu{x})u_h(\cu{x+k h}),\quad \cu{x}\in\cu{G}_h,
\label{eq:L_h}
\end{equation}
where the coefficients $c_{\cu{k}}(\cu{x})\in\mathbb{C}$, and $V$ is a finite index set. 
When $ c_{\cu{k}}(\cu{x}) $ is independent of $\cu{x}$, the operator $ L_h $ corresponds to a Toeplitz matrix, which can be diagonalized using Fourier modes. This leads to the following lemma.
\begin{lemma}[Lemma 4.2.1 \cite{trottenberg2000multigrid}]
For any $\cu{\theta} \in \Theta,$ the Fourier modes $\varphi(\cu{\theta,x})$ are formal eigenfunctions of the discrete operator $L_h$ with constant stencil
\begin{equation}
    L_h\varphi(\cu{\theta,x})=\tilde{L}_h(\cu{\theta})\varphi(\cu{\theta,x}),\quad \cu{x}\in\cu{G}_h,
    \label{eq:formal_eig}
\end{equation}
where
\begin{equation}
    \tilde{L}_h(\cu{\theta})=\sum_{k}c_{\cu{k}} e^{i\cu{\theta}\cdot \cu{k}},
    \label{eq:symbol}
\end{equation}
is called the formal eigenvalue or \textit{symbol} of $L_h$.
\end{lemma}
When $ c_{\cu{k}} $ depends on $\cu{x}$, the symbol $\tilde{L}_h(\cu{\theta})$ becomes a Fourier matrix function rather than a scalar function. The smoothing effect of $ L_h $ on different frequencies can still be analyzed through appropriate transformations. For further details, see \cite{bolten2018fourier, kumar2019local, brown2019local}.

Assume that the $k$-th iteration error of \eqref{eq:hybrid} is $\cu{e}^{(k)} = \cu{u - u}^{(k)}$, and $e^{(k)}(\cu{x})$ is its periodic extension function on the infinite grid $\cu{G}_h$.  Since $\{\varphi(\cu{\theta, x}), \cu{\theta} \in \Theta\}$ form an orthonormal basis for the Fourier space, $e^{(k)}(\cu{x})$ can be expanded in terms of these bases as
\begin{equation}
    e^{(k)}(\cu{x}) = \sum_{\cu{\theta} \in \Theta} \mu^{(k)}_{\cu{\theta}} \varphi(\cu{\theta, x}), \quad \cu{x} \in \cu{G}_h,
    \label{eq:err_k}
\end{equation}
and satisfies Parseval's identity
\begin{equation}
    \left\|e^{(k)}(\cu{x})\right\|^2_2 = \sum_{\cu{\theta} \in \Theta} \left|\mu^{(k)}_{\cu{\theta}}\right|^2.
    \label{eq:parseval}
\end{equation}

Since smoothing iteration \eqref{eq:hybrid1} is always implemented as a simple stationary iterative method, the corresponding error propagation matrix $\cu{E_B = I-BA}$ can be wirtten down as a discrete operator in the form of \eqref{eq:L_h}. Thus, 
\begin{equation}
    \cu{E_B}\varphi(\cu{\theta, x}) = \tilde{E}_B(\cu{\theta})\varphi(\cu{\theta, x}),
\end{equation}
then applying \eqref{eq:hybrid1} yields
\begin{equation}
    e^{(k+\frac{1}{2})}(\cu{x}) = \sum_{\cu{\theta} \in \Theta} \mu^{(k)}_{\cu{\theta}} \cu{E}^M_{\cu{B}} \varphi(\cu{\theta, x}) = \sum_{\cu{\theta} \in \Theta} \mu^{(k)}_{\cu{\theta}} \tilde{E}_B^M(\cu{\theta}) \varphi(\cu{\theta, x}), \quad \cu{x} \in \cu{G}_h.
    \label{eq:B_err}
\end{equation}
It can be seen that when $|\tilde{E}_B(\cu{\theta})| \ll 1$, the smoothing iteration effectively reduces the error component with frequency $\cu{\theta}$. However, when $|\tilde{E}_B(\cu{\theta})| \approx 1$ or greater than 1, the corresponding error components decay slowly or may even amplify.
Based on this fact, we make the following assumptions.
\begin{assumption}[Smoothing effect of $\cu{B}$]\label{ass:B}
    Assume that $\Theta$ can be split into $\Theta = \Theta^{\cu{B}} \cup \Theta^{\mathcal{H}}$, which satisfy
    \begin{equation}
    \begin{cases}
        |\tilde{E}_B(\cu{\theta})| \leq \mu_{\cu{B}}, \quad \text{if} \,\, \cu{\theta} \in \Theta^{\cu{B}}, \\
        \mu_{\cu{B}} < |\tilde{E}_B(\cu{\theta})| \leq 1 + \varepsilon_{\cu{B}}, \quad \text{if} \,\, \cu{\theta} \in \Theta^{\mathcal{H}},
    \end{cases}
    \label{eq:24}
    \end{equation}
    where the smoothing factor $\mu_{\cu{B}} \in (0,1)$ and the small positive constant $\varepsilon_{\cu{B}}$ are independent of the mesh size $h$ and physical parameters $\cu{\mu}$.
\end{assumption}
Under this assumption, smoothing error \eqref{eq:B_err} can be decomposed into
\begin{equation}
    e^{(k+\frac{1}{2})}(\cu{x}) = \sum_{\cu{\theta} \in \Theta^{\cu{B}}} \mu^{(k)}_{\cu{\theta}} \tilde{E}_B^M(\cu{\theta}) \varphi(\cu{\theta, x}) + \sum_{\cu{\theta} \in \Theta^{\mathcal{H}}} \mu^{(k)}_{\cu{\theta}} \tilde{E}_B^M(\cu{\theta}) \varphi(\cu{\theta, x}), \quad \cu{x} \in \cu{G}_h.
    \label{eq:err_B2}
\end{equation}

\subsection{Convergence of the neural operator $\mathcal{H}$}
Denote the restriction of $ e^{(k+\frac{1}{2})}(\cu{x}) $ and $ \varphi(\cu{\theta,x}) $ on $\mathcal{T}_h$ as
$$
\cu{e}^{(k+\frac{1}{2})} = e^{(k+\frac{1}{2})}(\cu{x})|_{\mathcal{T}_h} \quad \text{and} \quad \cu{\varphi(\theta)} = \varphi(\cu{\theta,x})|_{\mathcal{T}_h},
$$
then Eq. \eqref{eq:err_B2} becomes
\begin{equation}
    \cu{e}^{(k+\frac{1}{2})} = \sum_{\cu{\theta} \in \Theta^{\cu{B}}} \mu^{(k)}_{\cu{\theta}} \tilde{E}_B^M(\cu{\theta}) \cu{\varphi(\theta)} + \sum_{\cu{\theta} \in \Theta^{\mathcal{H}}} \mu^{(k)}_{\cu{\theta}} \tilde{E}_B^M(\cu{\theta}) \cu{\varphi(\theta)}.
    \label{eq:err_part_restricted}
\end{equation}
Applying the neural iteration \eqref{eq:hybrid2}, we have
\begin{equation}
    \begin{aligned}
        \cu{e}^{(k+1)} &= (\cu{I} - \mathcal{H} \cu{A}) \cu{e}^{(k+\frac{1}{2})} \\
        &= (\cu{I} - \mathcal{H} \cu{A}) \sum_{\cu{\theta} \in \Theta^{\cu{B}}} \mu^{(k)}_{\cu{\theta}} \tilde{E}_B^M(\cu{\theta}) \cu{\varphi}(\cu{\theta}) + (\cu{I} - \mathcal{H} \cu{A}) \sum_{\cu{\theta} \in \Theta^{\mathcal{H}}} \mu^{(k)}_{\cu{\theta}} \tilde{E}_B^M(\cu{\theta}) \cu{\varphi}(\cu{\theta}).
    \end{aligned}
\end{equation}

We make the following assumption on $\mathcal{H}$.
\begin{assumption}\label{ass:H}
Assume that $\mathcal{H}$ satisfies the following properties: $\exists \,\,  \epsilon >0$ and bounded constant $\mu_{\mathcal{H}}$, such that
\begin{equation}
    \begin{cases}(\cu{I}-\mathcal{H}\cu{A})\cu{\varphi}(\cu{\theta}) = \cu{\varphi}(\cu{\theta}),\quad\text{if}\,\,\cu{\theta}\in \Theta^{\cu{B}},\\
        \|(\cu{I}-\mathcal{H}\cu{A})\cu{\varphi}(\cu{\theta})\|_2\leq \mu_{\mathcal{H}} / N^{(1/2+\epsilon)},\quad\text{if}\,\,\cu{\theta}\in \Theta^{\mathcal{H}}.
      \end{cases}
      \label{eq:assH}
\end{equation}
\end{assumption}
Using the above assumption and computing the $\ell^2$ norm of $\cu{e}^{(k+1)}$, we obtain
\begin{equation}
    \begin{aligned}
    \|\cu{e}^{(k+1)}\|_2^2 &=\left\|\sum_{\cu{\theta} \in \Theta^{\cu{B}}}\mu^{(k)}_{\cu{\theta}}\tilde{E}_B^M(\cu{\theta})\cu{\varphi}(\cu{\theta})+\sum_{\cu{\theta} \in \Theta^{\mathcal{H}}}\mu^{(k)}_{\cu{\theta}}\tilde{E}_B^M(\cu{\theta})(\cu{I}-\mathcal{H}\cu{A})\cu{\varphi}(\cu{\theta})\right\|_2^2\\
    &\leq 2 \left\|\sum_{\cu{\theta} \in \Theta^{\cu{B}}}\mu^{(k)}_{\cu{\theta}}\tilde{E}_B^M(\cu{\theta})\cu{\varphi}(\cu{\theta})\right\|_2^2  +2\left\|\sum_{\cu{\theta} \in \Theta^{\mathcal{H}}}\mu^{(k)}_{\cu{\theta}}\tilde{E}_B^M(\cu{\theta})(\cu{I}-\mathcal{H}\cu{A})\cu{\varphi}(\cu{\theta})\right\|_2^2\\
    &\leq 2 \sum_{\cu{\theta} \in \Theta^{\cu{B}}}|\mu^{(k)}_{\cu{\theta}}\tilde{E}_B^M(\cu{\theta})|^2+2\left\|\sum_{\cu{\theta} \in \Theta^{\mathcal{H}}}\mu^{(k)}_{\cu{\theta}}\tilde{E}_B^M(\cu{\theta})(\cu{I}-\mathcal{H}\cu{A})\cu{\varphi}(\cu{\theta})\right\|_2^2\\
    &\leq 2\mu_{B}^{2M}\sum_{\cu{\theta} \in \Theta^{\cu{B}}}|\mu^{(k)}_{\cu{\theta}}|^2+2(1+\varepsilon_{\cu{B}})^{2M}\left(\sum_{\cu{\theta} \in \Theta^{\mathcal{H}}}|\mu^{(k)}_{\cu{\theta}}|\left\|(\cu{I}-\mathcal{H}\cu{A})\cu{\varphi}(\cu{\theta})\right\|_2\right)^2\\
    &\leq 2\mu_{B}^{2M}\sum_{\cu{\theta} \in \Theta^{\cu{B}}}|\mu^{(k)}_{\cu{\theta}}|^2+2(1+\varepsilon_{\cu{B}})^{2M}\left(\sum_{\cu{\theta} \in \Theta^{\mathcal{H}}}|\mu^{(k)}_{\cu{\theta}}|^2\right)\left(\sum_{\cu{\theta} \in \Theta^{\mathcal{H}}}\left\|(\cu{I}-\mathcal{H}\cu{A})\cu{\varphi}(\cu{\theta})\right\|_2^2\right)\\
    &\leq 2\mu_{B}^{2M}\sum_{\cu{\theta} \in \Theta^{\cu{B}}}|\mu^{(k)}_{\cu{\theta}}|^2+2(1+\varepsilon_{\cu{B}})^{2M}(\frac{\mu_{\mathcal{H}}}{N^{(1/2+\epsilon)}})^2|\Theta^{\mathcal{H}}|\left(\sum_{\cu{\theta} \in \Theta^{\mathcal{H}}}|\mu^{(k)}_{\cu{\theta}}|^2\right)\\
    &\leq \max{\left\{2\mu_{B}^{2M},C/N^{2\epsilon}\right\}}\|\cu{e}^{(k)}\|_2^2.
    \end{aligned}
\end{equation}
where $C= 2(1+\varepsilon_{\cu{B}})^{2M}\mu_{\mathcal{H}}^2\frac{|\Theta^{\mathcal{H}}|}{N}$ is bounded.

In summary, we obtain the following theorem.
\begin{theorem}[Main Result]\label{th:main}
    Assume that Assumption \ref{ass:B} and Assumption \ref{ass:H} are hold,  then the iterative error of the DL-HIM \eqref{eq:hybrid} satisfies
    \begin{equation}
        \|\cu{e}^{(k+1)}\|_2 \leq \eta \|\cu{e}^{(k)}\|_2,
    \end{equation}
    where the convergence rate given by
    $$
    \eta \leq \sqrt{\max{\left\{2\mu_{B}^{2M},C/N^{2\epsilon}\right\}}}.
    $$
     is independent of the grid size $h$.
\end{theorem}

Since most commonly used smoothers $\cu{B}$ always satisfy Assumption \ref{ass:B}, our focus is on designing a neural operator $\mathcal{H}$ that meets Assumption \ref{ass:H}. Existing DL-HIM frequently use standard neural operators as $\mathcal{H}$ \cite{hsieh2019learning, zhang2022hybrid, kahana2023geometry, hu2024hybrid, zou2024large, rudikov2024neural}, which exploit their spectral bias to effectively learn low-frequency error components. However, these methods often struggle to eliminate error components with other frequencies, resulting in inefficiencies when $\Theta^{\mathcal{H}}$ includes non-low frequencies. For instance, even in the case of Poisson equation, increasing discretization scales can lead to $\Theta^{\mathcal{H}}$ encompassing many intermediate frequencies, which typically requires additional operators or levels \cite{kopanivcakova2024deeponet} and can be costly. 
To address this challenge, we next propose an enhanced $\mathcal{H}$ based on our original FNS.

\section{Fourier neural solver}\label{sec:03}
The neural operator $\mathcal{H}$ used in the original FNS \cite{cui2022fourier} is given by
\begin{equation}
    \mathcal{H} = \mathcal{F}\cu{\tilde{\Lambda}}\mathcal{F}^{-1},
    \label{eq:matrix_H}
\end{equation}
where $\mathcal{F}$ is the Fourier matrix and $\cu{\tilde{\Lambda}}$ is a learned diagonal matrix. This design is inspired by the fast Poisson solver \cite{fortunato2020fast,strang2007computational}, where $\mathcal{F}$ is intended to approximate the eigenvector matrix of $\cu{A}$, and $\cu{\tilde{\Lambda}}$ approximates the inverse of the eigenvalue matrix. Although FNS demonstrated better convergence in some examples compared to existing neural solvers at that time, it has relatively obvious defects. Specifically, using $\mathcal{F}$ as a replacement for the eigenvector matrix of $\cu{A}$ is not always appropriate.
To illustrate this, consider the following two-dimensional (2D) Poisson equation
\begin{equation}
    \left\{
    \begin{aligned}
        -\Delta u &= f, \quad \cu{x} \in \Omega= (0,1)^2, \\
    u(\cu{x}) &= 0, \quad \cu{x} \in \partial \Omega.
    \end{aligned}\right.
    \label{eq:poisson}
\end{equation}
Using the five-point difference method on a uniform mesh with grid spacing $h=1/(n+1)$ in both the $x$- and $y$-directions, the eigenvalues of the resulting coefficient matrix are given by
$$
\lambda_{jx,jy} = \frac{4}{h^2} \sin^2\left(\frac{\pi j_x}{2(n+1)}\right) + \frac{4}{h^2} \sin^2\left(\frac{\pi j_y}{2(n+1)}\right),
$$
where $j_x=1,\ldots,n$ and $j_y=1,\ldots,n$.
The corresponding eigenvectors $\cu{\xi}_{jx,jy}$ are expressed as
$$
v_{ix,iy,jx,jy} = \frac{2}{n+1} \sin\left(\frac{i_{x} j_{x} \pi}{n+1}\right)\sin\left(\frac{i_{y} j_{y} \pi}{n+1}\right),
$$
where the multi-index $j_x,j_y$ pairs the eigenvalues and the eigenvectors, while the multi-index $i_x,i_y$ determines the location of the value of each eigenvector on the regular grid.

Let $N=n^2$ and $j=(j_x,j_y)$. Solving \eqref{eq:poisson} can be divided into the following three steps:
\begin{enumerate}[label=(\arabic*)]
    \item Expand $\cu{f}$ as a combination of the eigenvectors
    $$
    \cu{f}=a_{1} \cu{\xi}_{1}+\cdots+a_{N} \cu{\xi}_{N}.
    $$
    \item Divide each $a_{j}$ by $\lambda_{j}$.
    \item Recombine the eigenvectors to obtain $\cu{u}$
    $$
    \cu{u}=\left(a_{1} / \lambda_{1}\right) \cu{\xi}_{1}+\cdots+\left(a_{N} / \lambda_{N}\right) \cu{\xi}_{N}.
    $$
\end{enumerate}
Since the eigenvectors $\cu{\xi}_j$ are discrete sine functions, the first and third steps can be accelerated using the discrete sine transform, while the second step corresponds to a diagonal matrix multiplication (Hadamard product).

However, the eigenvectors of most PDE discrete operators are not discrete trigonometric functions, so directly using $\mathcal{F}$ as a replacement for eigenvector matrices is not suitable. On the other hand, both the Fourier basis and the eigenvector basis constitute orthonormal bases of $\mathbb{C}^{N}$. To better approximate the eigenvectors, we introduce a transition matrix $\cu{T} \in \mathbb{C}^{N \times N}$ such that $\cu{Q} = \mathcal{F}\cu{T}$. The target $\tilde{\mathcal{H}}$ then becomes
\begin{equation}
    \tilde{\mathcal{H}}=\cu{Q\Lambda}^{-1}\cu{Q}^{-1}.
\end{equation}
Since $\cu{Q}$ and $\mathcal{F}$ are both unitary operators, $\cu{T}$ is also a unitary operator. Thus,
$$
\cu{Q^{-1}} = \cu{T}^{-1}\mathcal{F}^{-1} = \cu{T}^{*}\mathcal{F}^{*},
$$
where $\cu{T}^{*}$ and $\mathcal{F}^{*}$ are the conjugate transposes of $\cu{T}$ and $\mathcal{F}$, respectively.

Moreover, we found that $\cu{T}$ has a sparse circulant structure for the Poisson equation, which inspired us to use a convolutional neural network (denoted as $\mathcal{C}$) to approximate $\cu{T}$.
In summary, the new improved $\mathcal{H}$ to approximate $\tilde{\mathcal{H}}$ is
\begin{equation}
    \mathcal{H} = \mathcal{F}\mathcal{C}\cu{\tilde{\Lambda}}\mathcal{C}^{*}\mathcal{F}^{-1}, \quad \mathcal{C} \approx \cu{T}.
    \label{eq:newH}
\end{equation}
\begin{remark}
    Note that the frequency of the $j$-th Fourier basis
    \begin{equation}
        F^j_k = \frac{1}{\sqrt{N}}e^{-\frac{2\pi i}{N} jk},\quad 0 \leq j,k \leq N-1,
    \end{equation}
    is $\theta = 2\pi j$, while the frequency of $\cu{\xi}_{j}$ is $\theta = \pi j$. In practical calculations, to approximate the eigenvector with the Fourier basis function and transition matrix easily, we extend $\cu{f} \in \mathbb{C}^N$ to $\tilde{\cu{f}} \in \mathbb{C}^{N^{\prime}}$, with $N^{\prime} = 2(N+1)$, using odd reflection symmetry about each wall. 
    Then the entry of Fourier matrix becomes $\exp(-i2\pi jk/N^{\prime}) = \exp(-i\pi jkh)$, whose frequency is also $\theta = \pi j$. The specific process will be demonstrated in the next section.
\end{remark}

Next, we estimate $(\cu{I}-\mathcal{H}\cu{A})\cu{\varphi}(\cu{\theta})$ to discuss whether the improved $\mathcal{H}$ \eqref{eq:newH} can satisfy Assumption \ref{ass:H}. Note that $\cu{\varphi}(\cu{\theta}) \in \C^N$ can be expressed as a linear combination of the eigenvector basis. Specifically, there exists a set of coefficients ${t_{\cu{\theta}}^i \in \C}$ such that
\begin{equation}
    \cu{\varphi}(\cu{\theta}) = \sum _{i=1}^N\cu{\xi}_{i}t_{\cu{\theta}}^i:=\cu{Qt_{\theta}},\,.
    \label{eq:phi_rep}
\end{equation}
where the vector $\cu{t_{\theta}} = [t_{\cu{\theta}}^1,\cdots,t_{\cu{\theta}}^N]^T$ is one column of the sparse matrix $\cu{T}$, i.e. there exists an index set $V_{\cu{\theta}}$ with $|V_{\cu{\theta}}| = O(1)$, such that $t_{\cu{\theta}}^i \neq 0$ if and only if $i \in V_{\cu{\theta}}$.
Using \eqref{eq:phi_rep} and assume that $\cu{A}= \cu{Q} \cu{\Lambda} \cu{Q}^{-1}$,
we have
\begin{equation}
    (\cu{I}-\mathcal{H}\cu{A})\cu{\varphi}(\cu{\theta}) =(\cu{I}-\mathcal{H}\cu{A})\cu{Qt_{\theta}} 
    =(\cu{Q}-\mathcal{H}\cu{Q\Lambda Q^{-1}Q})\cu{t_{\theta}}
    =(\cu{Q}-\mathcal{H}\cu{Q\Lambda})\cu{t_{\theta}},
\label{eq:30}
\end{equation}

According to Assumption \ref{ass:H}, we divide the index set $I = [1, \cdots, N]$ into the following two parts
\begin{equation}
    I^{\mathcal{H}}=\bigcup_{\cu{\theta} \in \Theta^{\mathcal{H}}}V_{\cu{\theta}}, \quad
    I^{\cu{B}}=I \setminus I^{\mathcal{H}}.
\end{equation}
A sufficient condition for $\mathcal{H}$ to satisfy Assumption \ref{ass:H} is
\begin{equation}
    \begin{cases}\cu{l}_i=\cu{\xi}_i,\,\,i \in I^{\cu{B}},\\
        \|\cu{l}_i\|_2\leq \frac{\mu_{\mathcal{H}}}{|V_{\cu{\theta}}|N},\,\,i \in I^{\mathcal{H}}.
      \end{cases}
      \label{eq:assl}
\end{equation}

In Section 4.1, we take the Poisson equation as an example to numerically verify the feasibility of the above assumptions.
Two subnetworks, Meta$-\lambda$ and Meta$-T$, are introduced to leverage the information from the PDE parameters $\cu{\mu}$ of the discrete system to provide $\cu{\tilde{\Lambda}}$ and the convolution kernel of $\mathcal{C}$, respectively. Specifically, Meta$-\lambda$ employs the widely used Fourier Neural Operator (FNO) \cite{li2020fourier}, while Meta$-T$ uses a convolutional neural network. By leveraging the powerful  approximation capabilities of neural networks \cite{hornik1989multilayer,kovachki2023neural}, these assumptions are ensured to hold.

\begin{figure}[!htbp]
    \centering
    \includegraphics[width=0.8\textwidth]{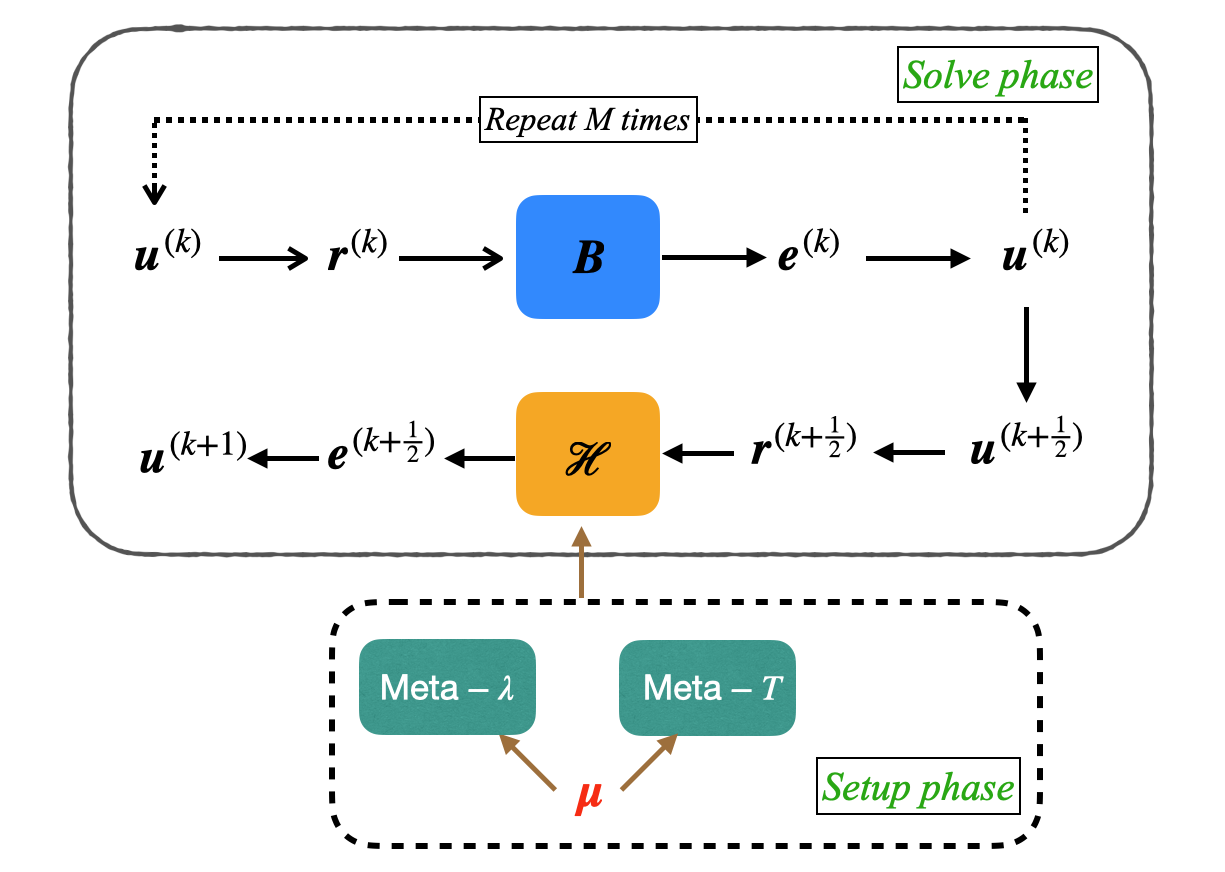}
    \caption{The schematic diagram for the calculation flow of FNS.}
    \label{fig:FNS}
\end{figure}

The schematic diagram of FNS is illustrated in Figure \ref{fig:FNS}, which consists of two stages: setup phase and solve phase. During the setup phase, two meta subnets are used to give the parameters required for the $\mathcal{H}$. The specific operation of $\mathcal{H}$ will be detailed in the subsequent section. The meta subnets incorporate the nonlinear activation function to improve their expressive ability, which does not break the linearity of $\mathcal{H}$ with respect to its input. The linearity of $\mathcal{H}$ with scale-equivalent property can give the correction error associated with residuals with different orders of magnitudes, thereby mitigating the issue of poor  generalization in neural networks for data exceeding two orders of magnitude \cite{CSIAM-AM-4-13}. The computational complexity of a single-step iteration during the solve phase is $ O(N \log N) $.

Returning to the question raised at the end of the previous section, the improved $\mathcal{H}$ in \eqref{eq:newH} overcomes spectral bias by learning the corresponding eigenvalues and eigenvectors in the frequency domain with the help of additional meta subnets. This provides a structural advantage over the $\mathcal{H}$ used in existing DL-HIM. Next, we construct appropriate training data and loss functions for learning FNS to meet the requirements of Theorem \ref{th:main}.

\subsection{Training data and loss function}
Some DL-HIM \cite{zhang2022hybrid,kahana2023geometry,hu2024hybrid,zou2024large} train $\mathcal{H}$ separately using supervised learning with a loss function of the form
$$
\mathcal{L} = \frac{1}{N_{\text{train}}}\sum_{i=1}^{N_{\text{train}}}\frac{\|\cu{u}_i - \mathcal{H}\cu{f}_i\|}{\|\cu{u}_i\|}.
$$
Here, the data pairs $(\cu{f}_i, \cu{u}_i)$ are generated by solving PDEs. This approach is not only computationally expensive but also contain all frequency error components without accounting for the role of $\cu{B}$.
In training neural solvers, one common approach is to sample $\cu{u_i}$ from a certain distribution, such as $\mathcal{N}(0, \cu{I})$, and obtain $\cu{f_i}$ by computing $\cu{f_i} = \cu{A u_i}$. However, this method results in $\cu{f_i} \sim \mathcal{N}(0, \cu{AA^T})$, which is biased toward the dominant eigen-subspace of $\cu{AA^T}$. Consequently, this can lead to poor performance of the neural solver when dealing with inputs that are close to the bottom eigen-subspace \cite{chen2024graph}.

The loss function we employed is the relative residual
\begin{equation}
    \mathcal{L} = \frac{1}{N_{\text{train}}}\sum_{i=1}^{N_{\text{train}}}\frac{\|\cu{f}_i - \cu{A}_i \cu{u}_i^K\|}{\|\cu{f}_i\|},
    \label{eq:loss_func}
\end{equation}
where $\cu{u}_i^K$ is the iterative solution obtained after applying \eqref{eq:hybrid} $K$ times, starting from a zero initial guess, and $\cu{f}_i$ is sampled from $\mathcal{N}(0, \cu{I})$.
This approach ensures that the input to $\mathcal{H}$ follows the distribution $\mathcal{N}(0, \cu{(I - AB)} \cu{(I - AB)^T})$. When $\cu{B}$ is selected as a simple preconditioner, such as the diagonal Jacobi preconditioner, this distribution tends to be skewed toward the dominant eigen-subspace of $\cu{I - BA}$, which presents challenges for $\cu{B}$ to effectively handle. Moreover, instead of constructing data for each fixed $\cu{A_i}$, we generate various $\cu{A_i}$ matrices from the discretized PPDE \eqref{ppde}. Since the discrete systems considered in this paper are based on structured grids, and all matrix-vector multiplications can be performed using convolution, there is no need to store the sparse matrix $\cu{A_i}$. Our training data consists of tuples $(\cu{\mu}_i, \cu{f}_i)$.

\section{Numerical experiments}\label{sec:04}
In this section, we evaluate the performance of FNS on discrete systems derived from several types of PPDE. These systems vary in their algebraic properties and the challenges they present, including:
\begin{enumerate}[label=(\arabic*)]
    \item Poisson equation: SPD.
    \item Random diffusion equation: Poisson-like elliptic problem, SPD.
    \item Anisotropic diffusion equation: Multiscale, SPD.
    \item Convection-diffusion equation: Non-symmetric, positive definite.
    \item Jumping diffusion equation: Multiscale with different scales at different grid points, SPD.
    \item Helmholtz equation: Complex, indefinite, non-Hermitian.
\end{enumerate}

\begin{figure}[!htb]
    \centering
    \includegraphics[width=\textwidth]{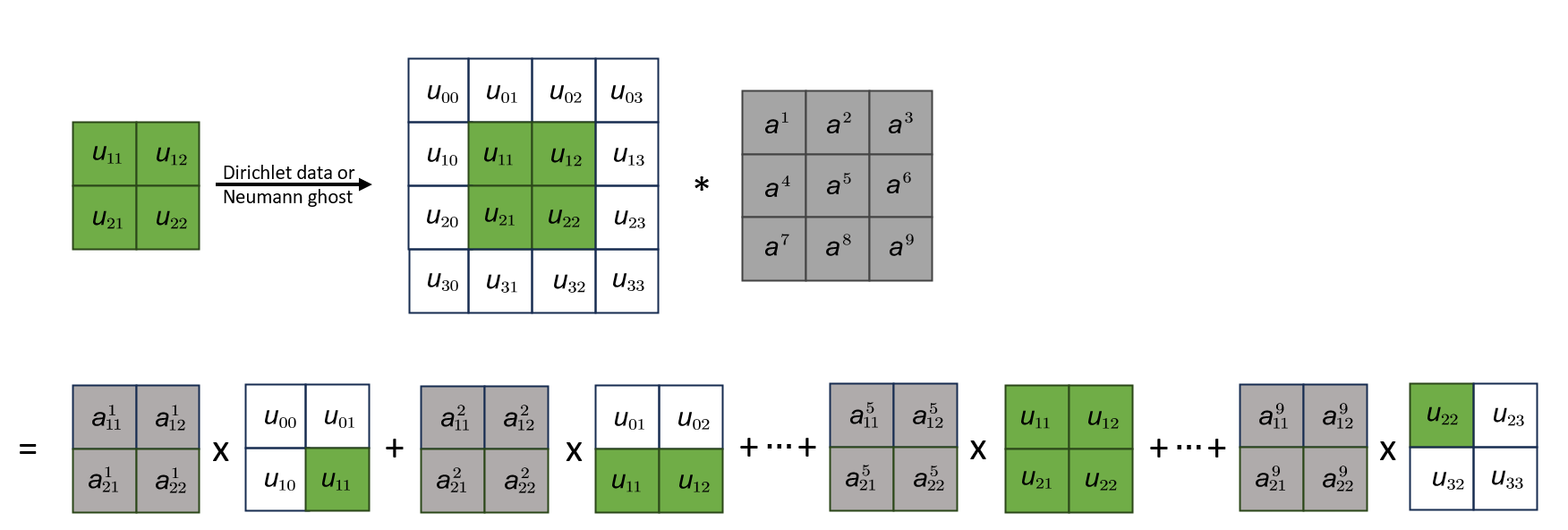}
    \caption{Matrix-vector product with a variable stencil.}
    \label{fig:conv}
\end{figure}

We implement FNS using the PyTorch deep learning framework \cite{Paszke2019PyTorchAI} and conduct numerical experiments on an Nvidia A100-SXM4-80GB GPU. For general matrix-vector multiplication, we adopt a matrix-free implementation approach.
For a general 9-point scheme, the computation of $f_{ij}$ is given by:
$$
\begin{aligned}
    f_{i j} & =a_{i-1, j+1} u_{i-1,j+1}+a_{i, j+1} u_{i, j+1}+a_{i+1, j+1} u_{i+1, j+1} \\
    & +a_{i-1,j} u_{i-1, j}+a_{i, j} u_{i, j}+a_{i+1, j} u_{i+1, j} \\
    & +a_{i-1, j-1} u_{i-1,j-1}+a_{i, j-1} u_{i,j-1}+a_{i+1, j-1} u_{i+1,j-1}
    \end{aligned}.
$$
As illustrated in Figure \ref{fig:conv}, we first compute each term of the product for all nodes $i, j$, resulting in nine channels of tensors. These nine channels are then summed element-wise to obtain $\cu{f}$.

\subsection{Poisson equation}
\begin{figure}[!htbp]
    \centering
    \subfigure[$\omega=1/2$]{\includegraphics[width=0.3\textwidth]{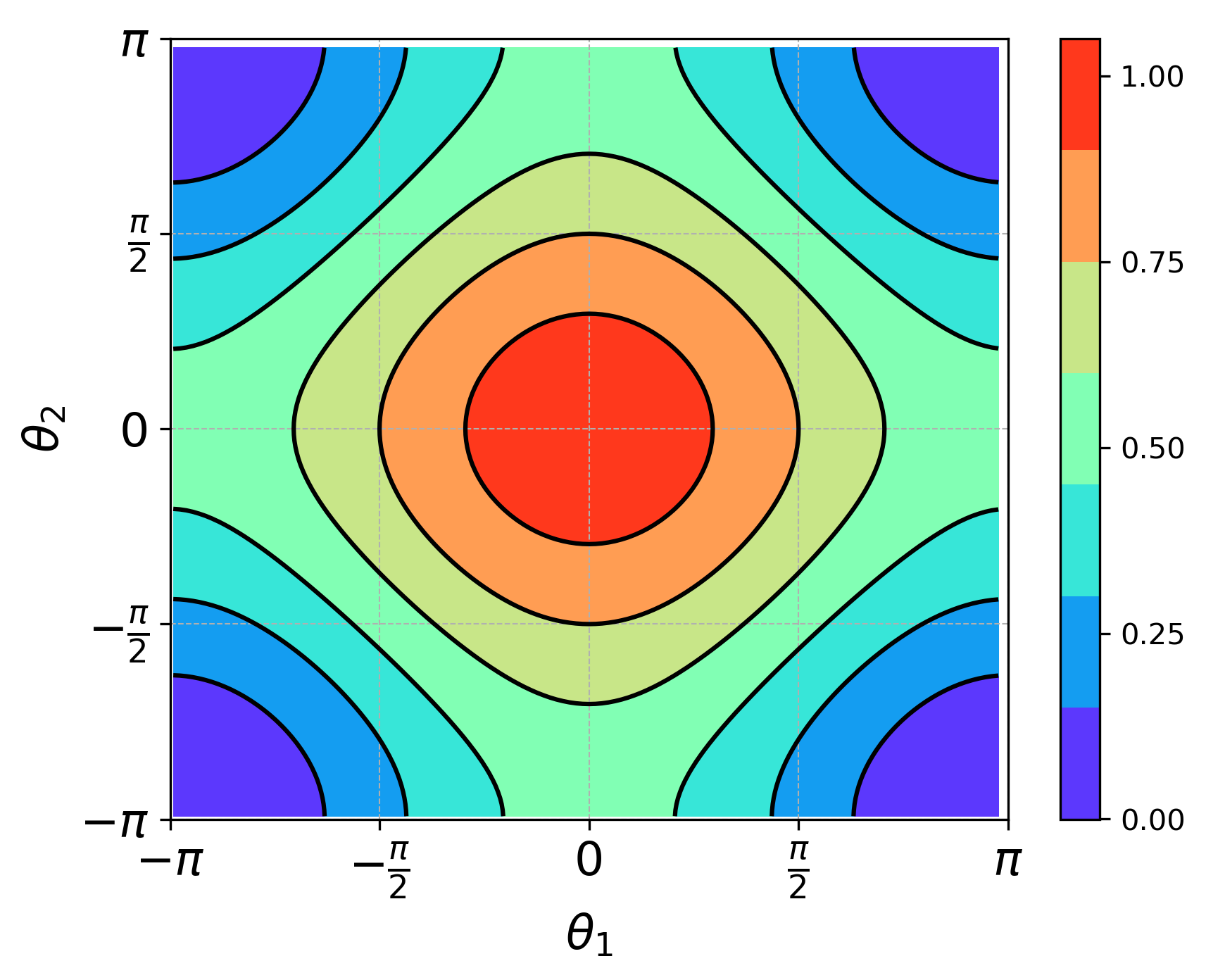}}\quad
    \subfigure[$\omega=3/4$]{\includegraphics[width=0.3\textwidth]{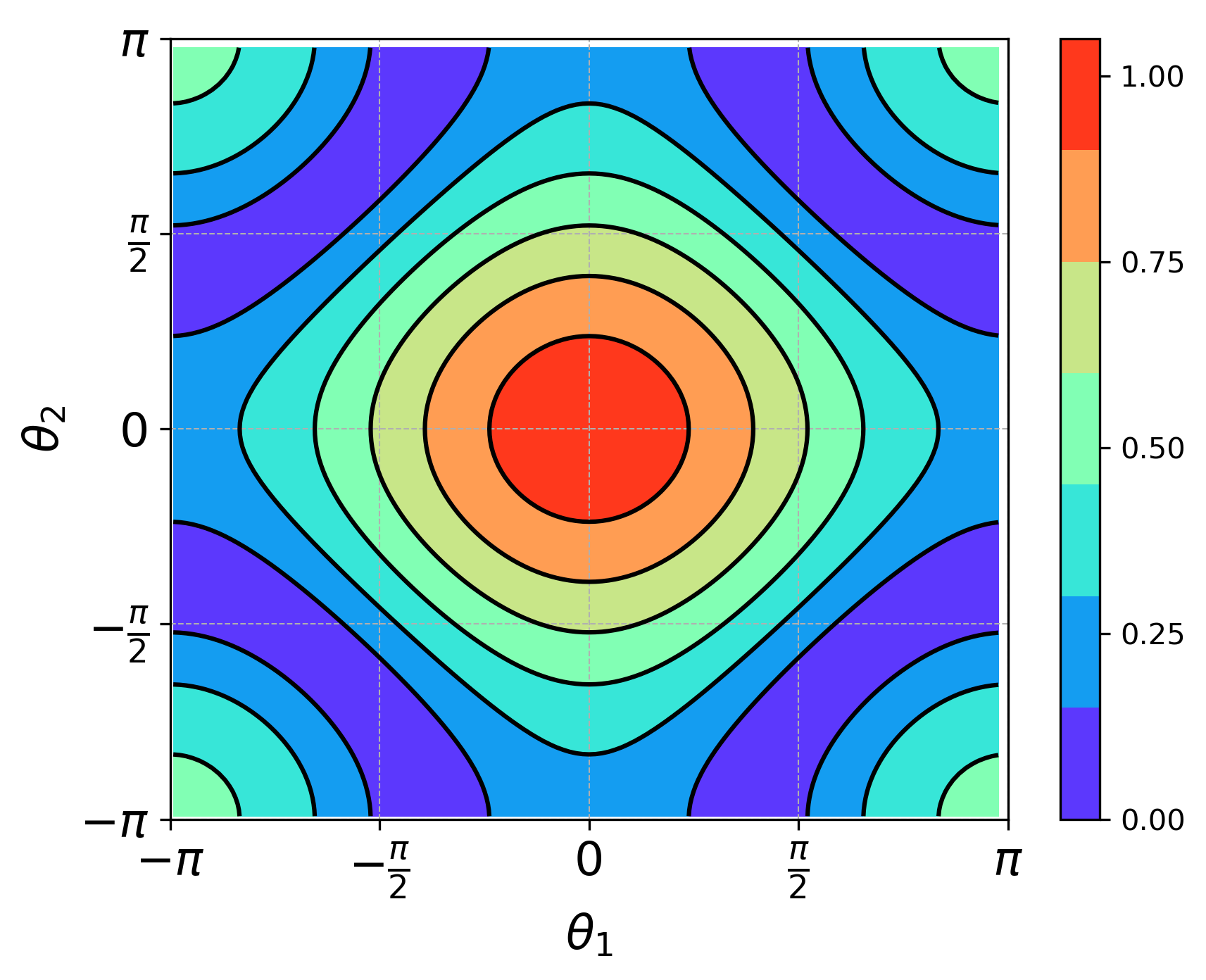}}\quad
    \subfigure[$\omega=1$]{\includegraphics[width=0.3\textwidth]{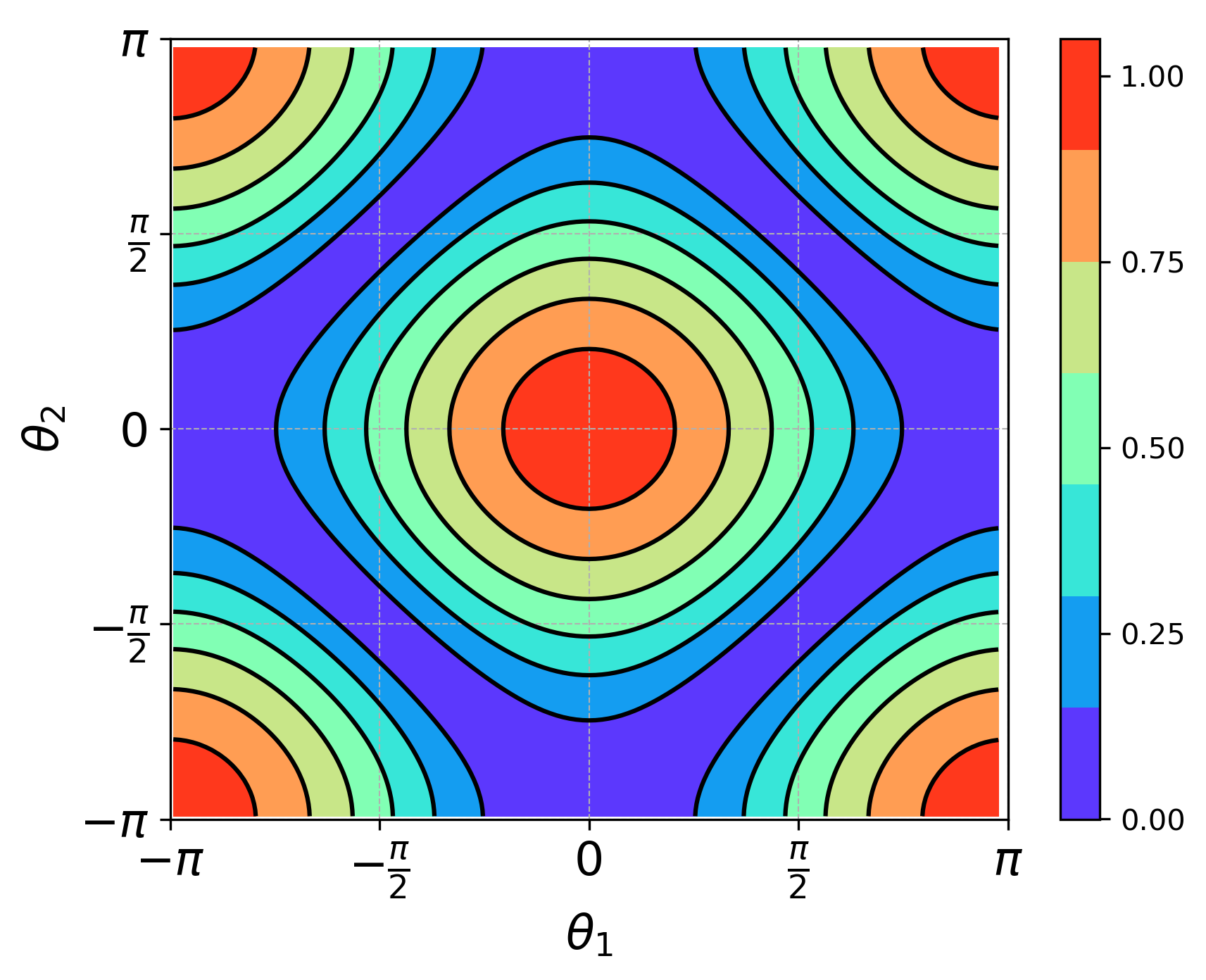}}
    \caption{Jacobi symbol with different weights when solving the Poisson equation.}
    \label{fig:LFA_poisson}
\end{figure}

We first verify the performance of FNS on the Poisson equation \eqref{eq:poisson}.
Taking $\cu{B}$ as the damped Jacobi method, we use LFA to evaluate its compressibility for error components with different frequencies. Figure \ref{fig:LFA_poisson} shows the modulus of the formal eigenvalue of the damped Jacobi method (abbreviated as the Jacobi symbol) with different weights when solving the discrete system. It can be observed that the damped Jacobi method exhibits poor compressibility for error components with frequencies in $[- \pi/2, \pi/2)^2$, regardless of the weights $\omega$. Let $\Theta^{\mathcal{H}} = [-\pi/2, \pi/2)^2$ and $\Theta^{\cu{B}} = [-\pi, \pi)^2 \setminus \Theta^{\mathcal{H}}$, it can be found that when $\omega = 3/4$, the damped Jacobi method achieves the best compressibility for the error components with frequencies in $\Theta^{\cu{B}}$. Therefore, we choose $\cu{B}$ as the damped Jacobi method with $\omega = 3/4$, which can meet Assumption \ref{ass:B}. To achieve a better smoothing effect, we set $M=10$.

\begin{figure}[!htb]
    \centering
    \subfigure[Low frequency]{
    \includegraphics[width=0.4\textwidth]{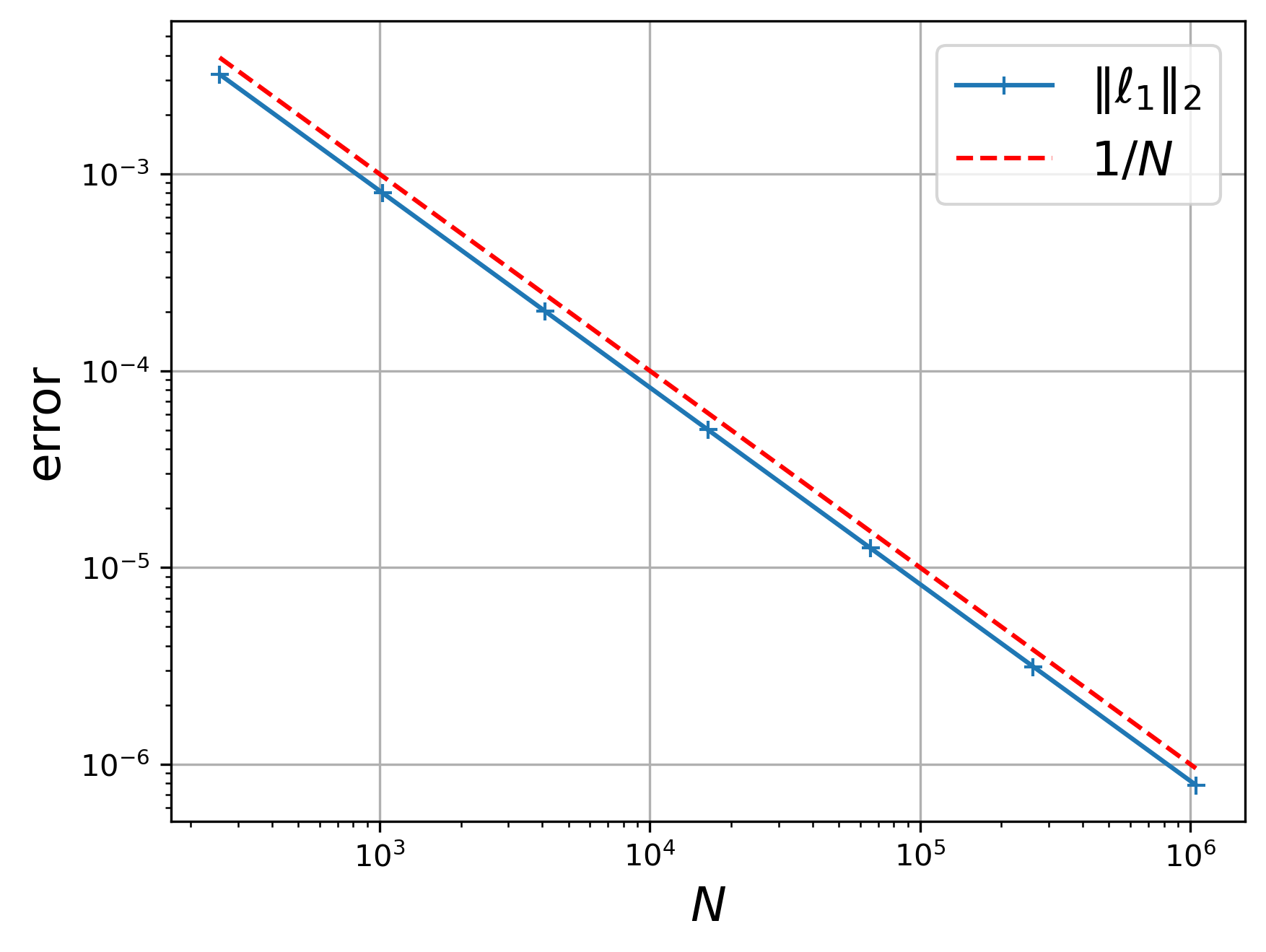}}\quad
    \subfigure[High frequency]{
    \includegraphics[width=0.4\textwidth]{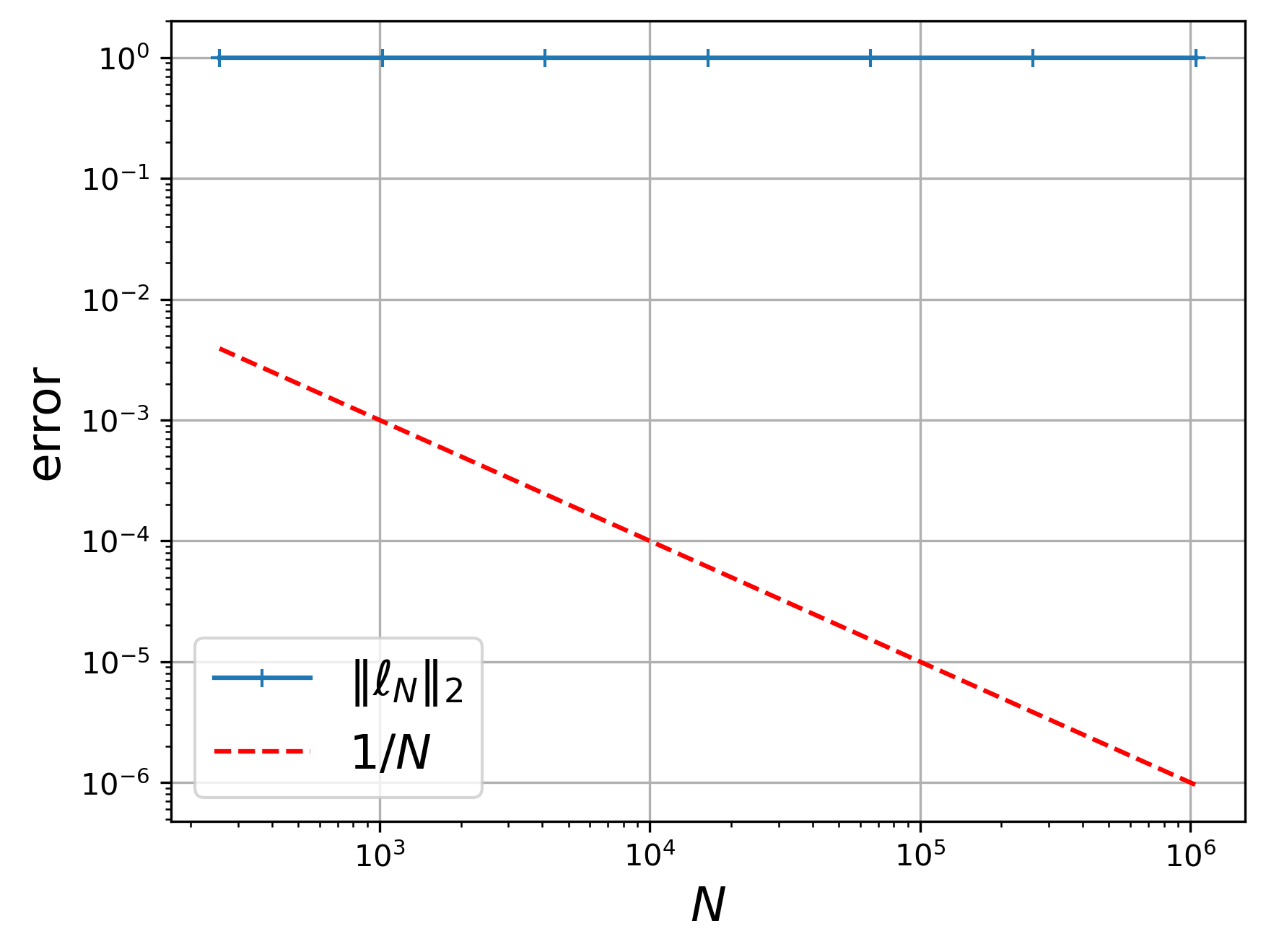}}
    \caption{Approximation error of $\mathcal{H}$ when applied to the low- and high-frequency eigenvectors of the Poisson equation.}
    \label{fig:poisson_H}
    \end{figure}

To construct $\mathcal{H}$ to eliminate the error components with frequencies belonging to $\cu{\Theta}^{\mathcal{H}}$, following the fast Poisson solver \cite{fortunato2020fast}, we take $\mathcal{C}$ as the identity operator, and
$$
\tilde{\cu{\Lambda}}\left(\theta_1, \theta_2\right)=\left\{\begin{array}{cl}
 \frac{h^2}{\theta_1^2+\theta_2^2} & \left(\theta_1, \theta_2\right) \neq 0, \left(\theta_1, \theta_2\right) \in \cu{\Theta}^{\mathcal{H}} \\
1 & \theta_1=\theta_2=0 \\
0 & \text { otherwise }
\end{array}.\right.
$$
We then verify whether $\mathcal{H}$ satisfies assumption \eqref{eq:assl} numerically. Figure \ref{fig:poisson_H} illustrates the approximation error of $\mathcal{H}$ for representative low-frequency and high-frequency eigenvectors. It can be observed that for the low-frequency eigenvector ($j=1$, on the left), $\mathcal{H}$ achieves an $\mathcal{O}(1/N)$ approximation rate, whereas it does not alter the high-frequency eigenvector ($j=N$, on the right). This behavior is consistent with assumption \eqref{eq:assl} with $\epsilon=1/2$.

\begin{figure}[!htbp]
    \centering
    \includegraphics[width=0.6\textwidth]{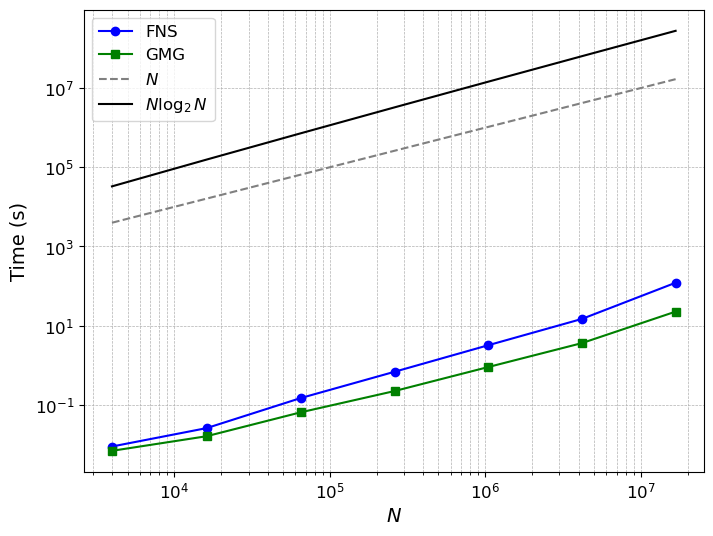}
    \caption{Solving time of FNS and GMG for the Poisson equation at different scales.}
    \label{fig:FNSvsGMG}
\end{figure}

Finally, we evaluate the performance of FNS in solving discretized systems at various scales and compare it with geometric multigrid (GMG). For each scale, 10 random RHS are generated, and the average iteration count and computation time required to achieve a relative residual below $10^{-6}$ are recorded.
Experiment results show that the number of iterations for FNS consistently remains at 9, regardless of the problem scale or the specific RHS. Figure \ref{fig:FNSvsGMG} illustrates the computation time required by both solvers as the problem scale increases. The results indicate that both solvers demonstrating approximately linear solution time scaling with the length of the solution vector.

\subsection{Random diffusion equations}
\begin{figure}[!htb]
     \centering
     \subfigure[$a(\cu{x})$]{\label{fig:a}
     \includegraphics[width=0.3\textwidth]{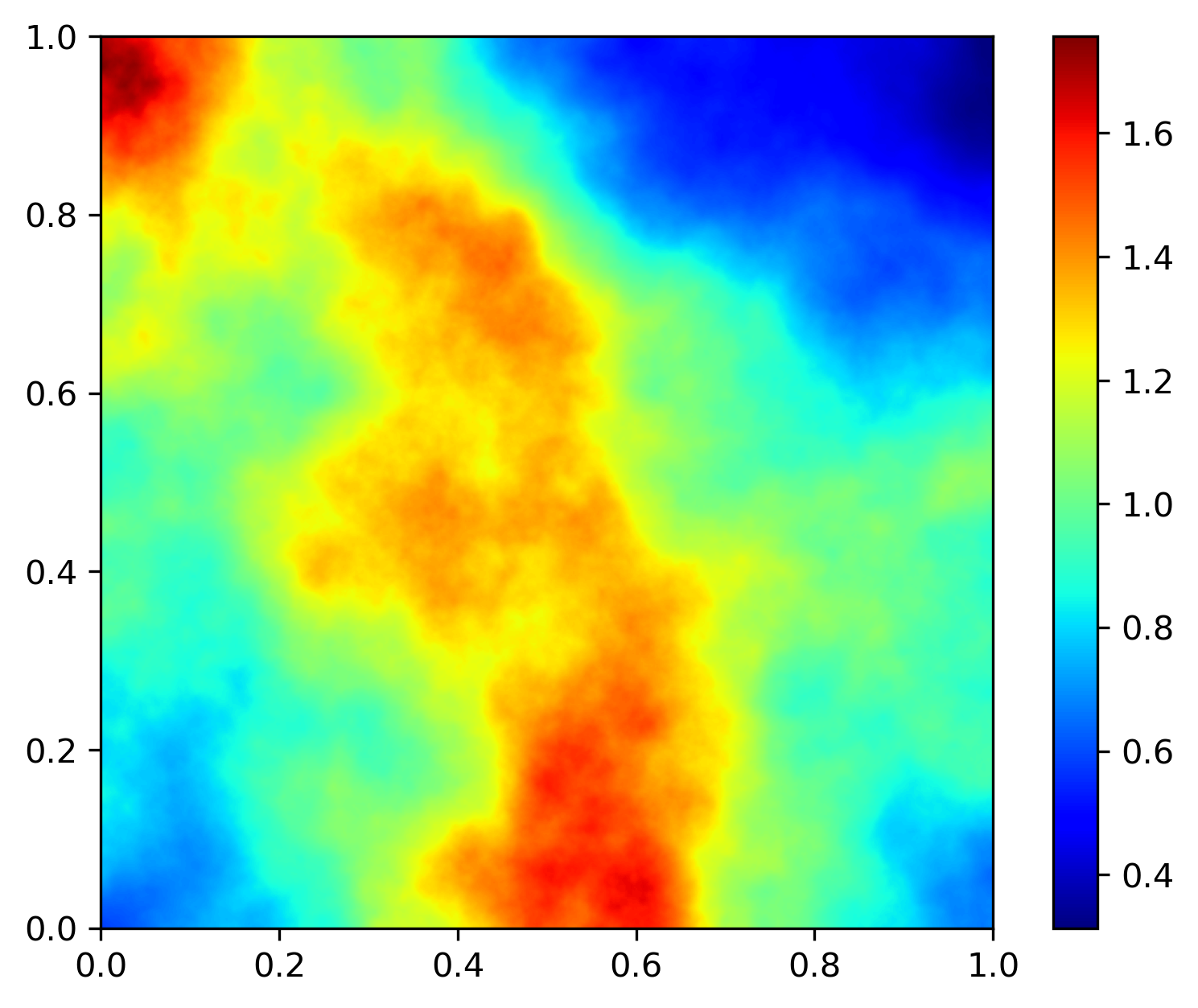}}
     \subfigure[$u(\cu{x})$]{\label{fig:u}
     \includegraphics[width=0.3\textwidth]{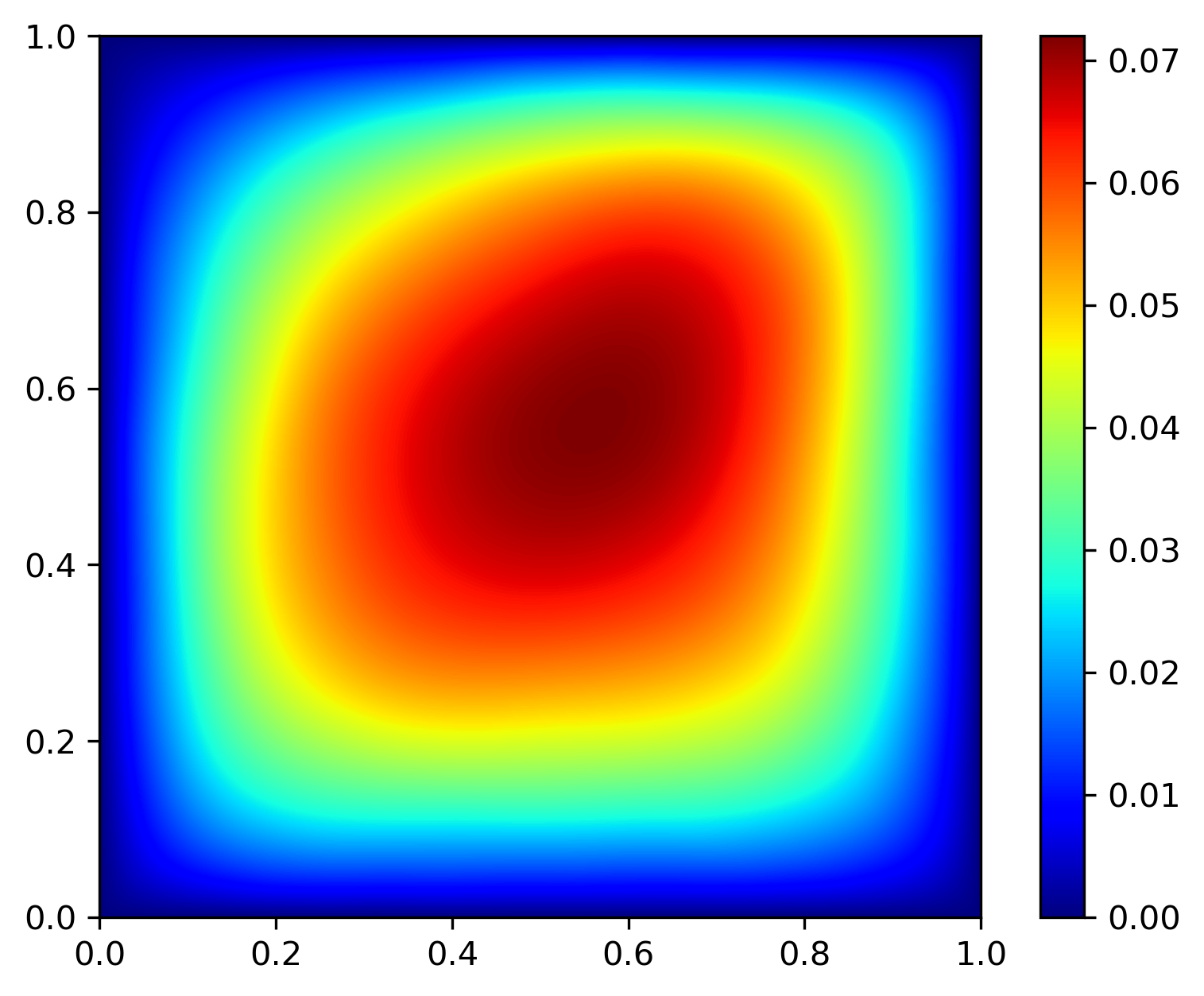}}
     \subfigure[Grid points]{\label{fig:grid_point}
     \includegraphics[width=0.3\textwidth]{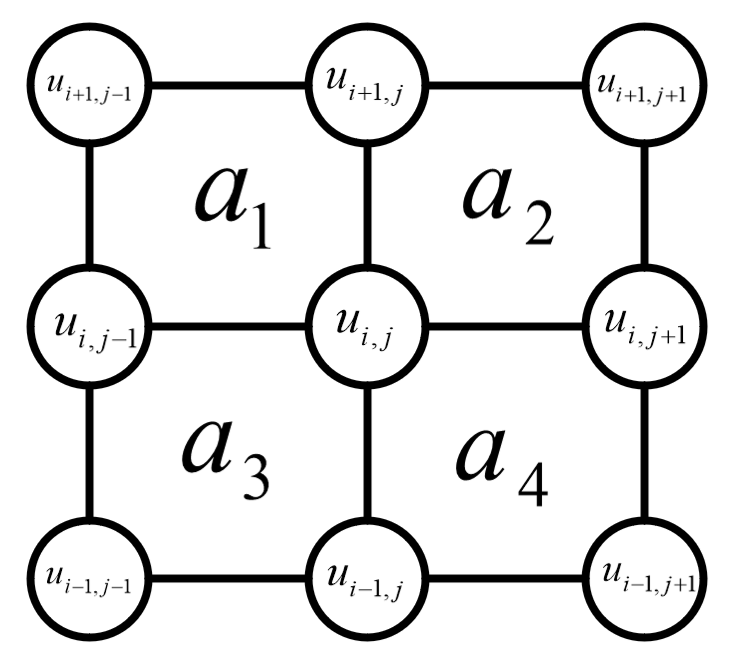}}
     \caption{An example of the random diffusion equation and discretization grid points.}
\end{figure}

Consider the following 2D random diffusion equation
\begin{equation}
    \begin{aligned}
    -\nabla \cdot (a(\cu{x})\nabla u(\cu{x})) &= f, \quad \cu{x} \in \Omega, \\
    u(\cu{x}) &= 0, \quad \cu{x} \in \partial \Omega,
    \end{aligned}
    \label{eq:random}
\end{equation}
where the diffusion coefficient $a \sim \psi_{\#}\mathcal{N}(0,(-\Delta + 9I)^{-2})$, and $\psi$ is the exponential function \cite{li2020fourier}. 
Figure \ref{fig:a} shows an example of $a$, and Figure \ref{fig:u} shows the corresponding numerical solution when $f=1$. 

We use bilinear FEM on squares with sides of length $ h = 1/(n+1) $. The discrete equation at an interior mesh point with indices $(i, j)$ is
\begin{equation}
    \begin{aligned}
    & \frac{2}{3}(a_{1} + a_{2} + a_{3} + a_{4})u_{i,j} - \frac{1}{3}(a_{1}u_{i+1,j-1} + a_{2}u_{i+1,j+1} + a_{3}u_{i-1,j-1} + a_{4}u_{i-1,j+1}) \\
    & - \frac{1}{6}((a_{3} + a_{4})u_{i-1,j} + (a_{1} + a_{3})u_{i,j-1} + (a_{2} + a_{4})u_{i,j+1} + (a_{1} + a_{2})u_{i+1,j}) = h^2f_{i,j},
    \end{aligned}
\end{equation}
where $a_1, a_2, a_3,$ and $a_4$ are constant random coefficients corresponding to the four neighboring elements of index $(i, j)$, as shown in Figure \ref{fig:grid_point}.

\begin{figure}[!htbp]
    \centering
    \subfigure[$\omega=1/2$]{\includegraphics[width=0.3\textwidth]{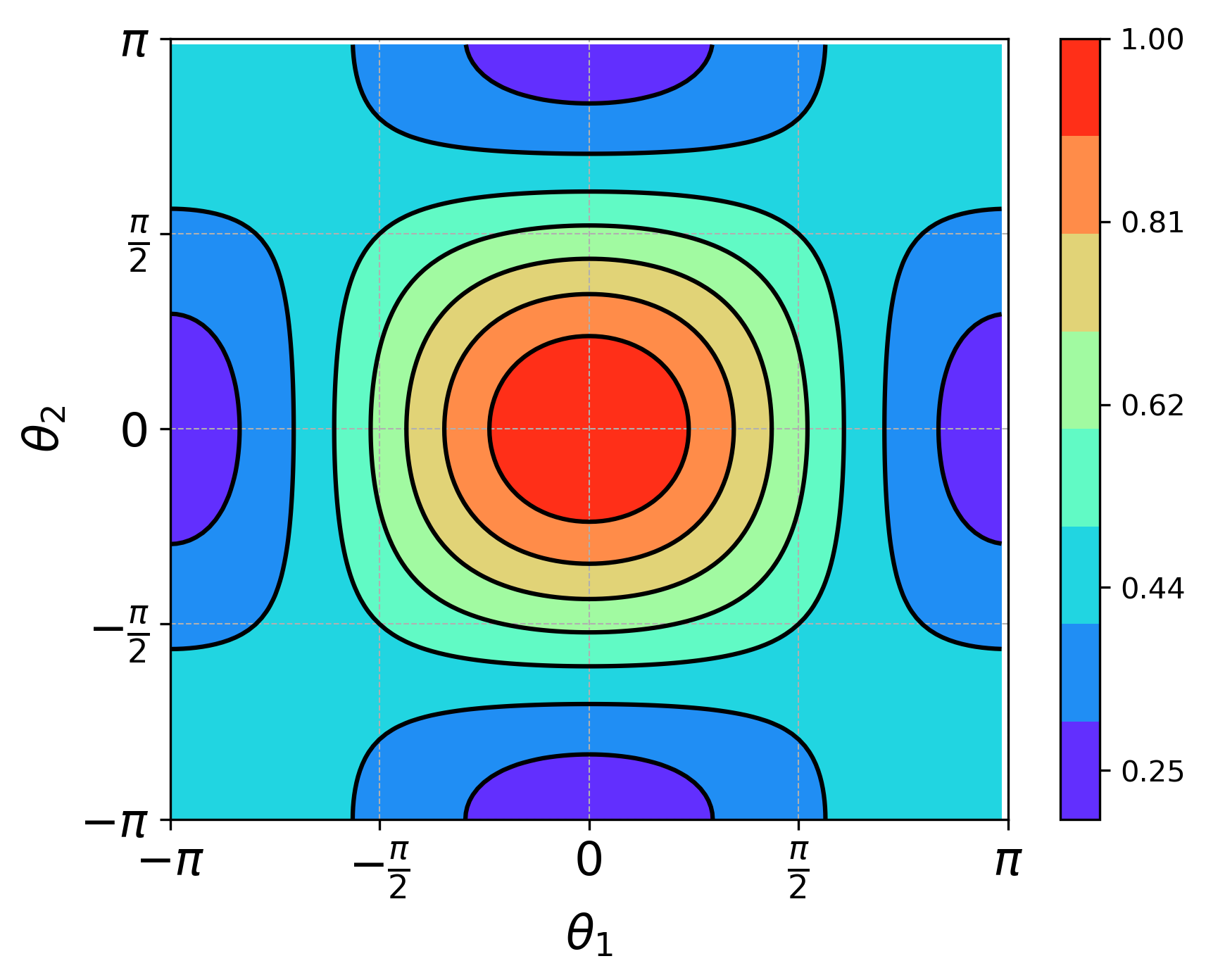}}\quad
    \subfigure[$\omega=3/4$]{\includegraphics[width=0.3\textwidth]{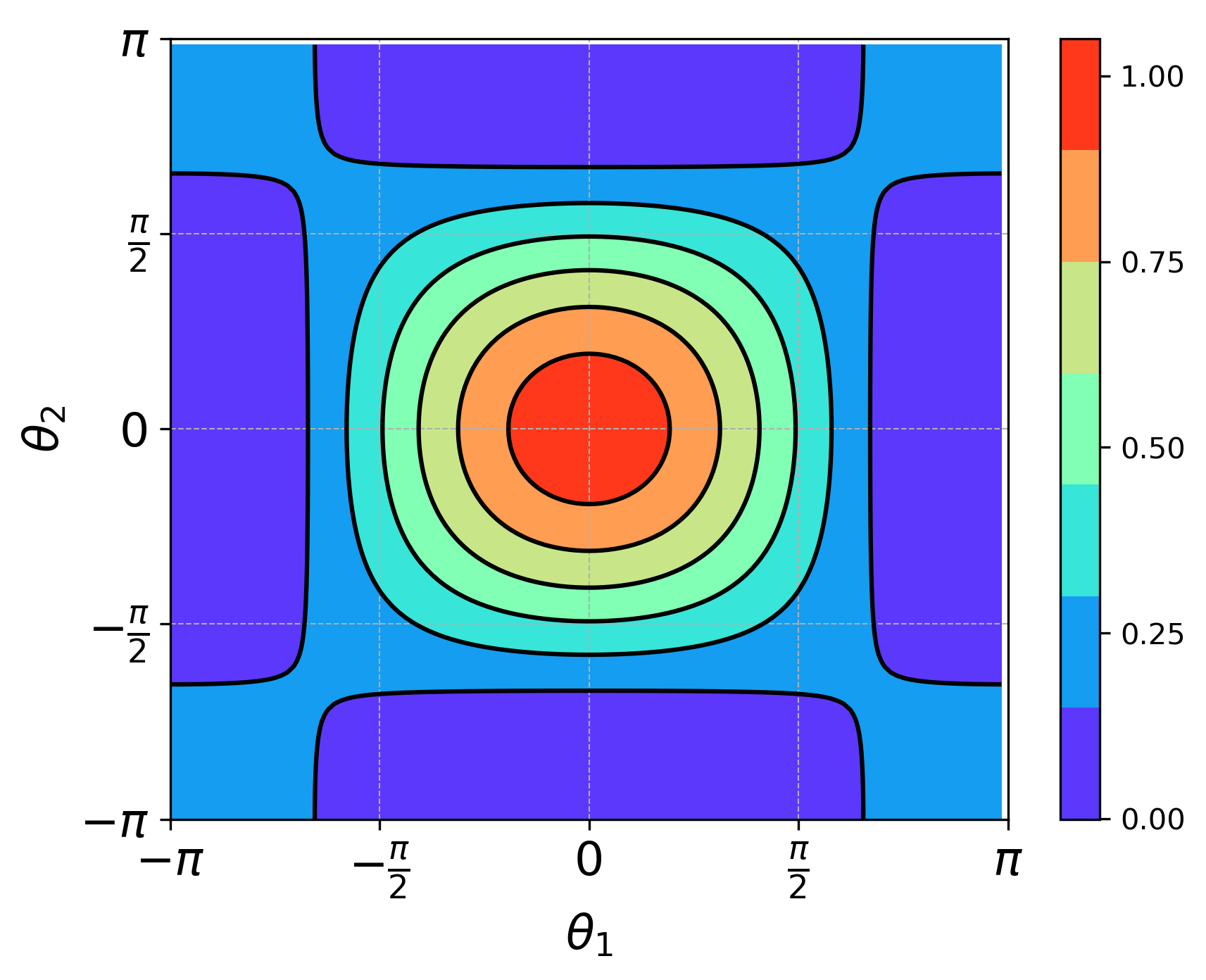}}\quad
    \subfigure[$\omega=1$]{\includegraphics[width=0.3\textwidth]{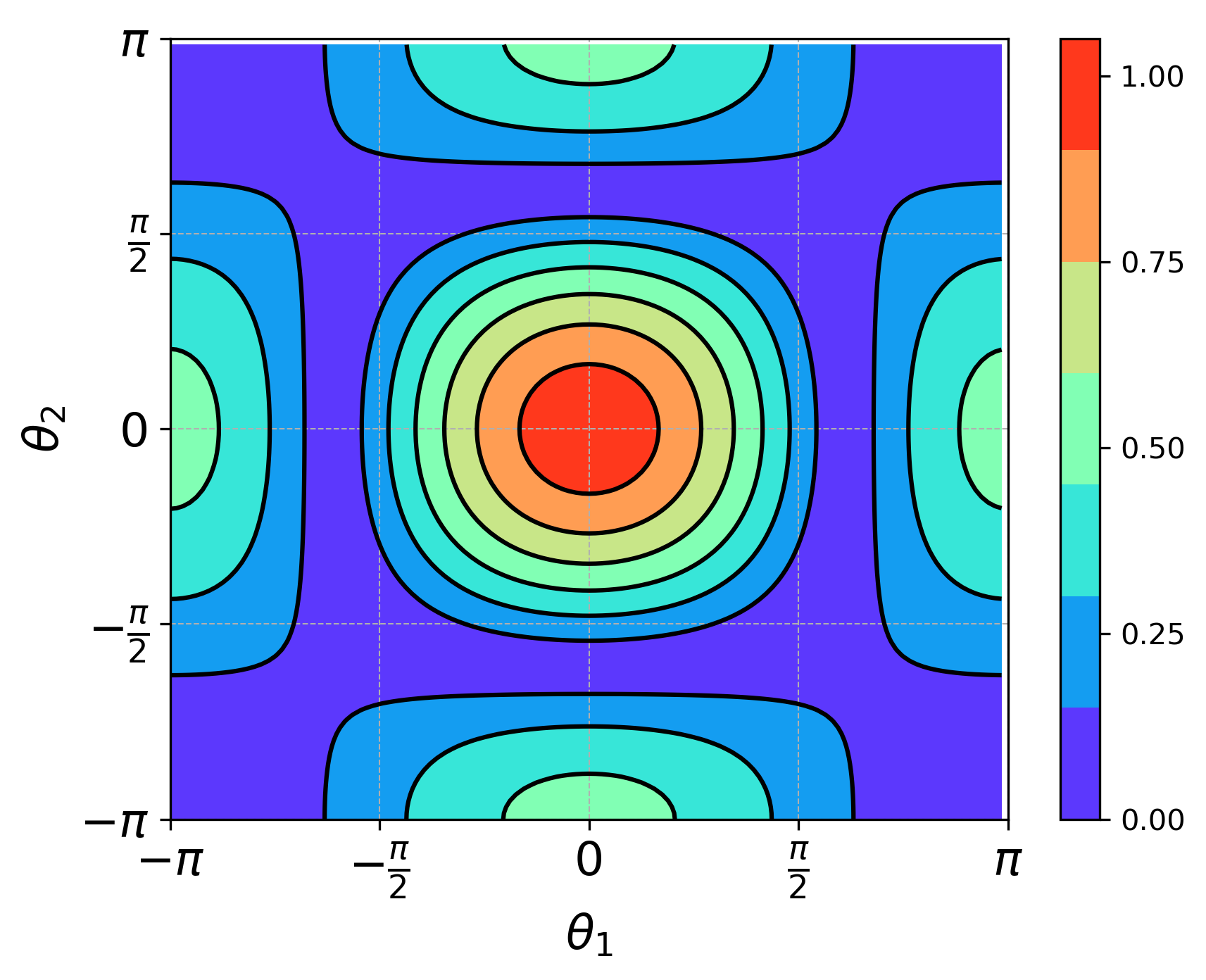}}
    \caption{Jacobi symbol with different weights when solving the random diffusion equation.}
    \label{fig:LFA1}
\end{figure}

Selecting $\cu{B}$ as the damped Jacobi method, we use non-standard LFA \cite{bolten2018fourier} to evaluate its compressibility for error components with different frequencies. Figure \ref{fig:LFA1} shows the modulus of the formal eigenvalue of the damped Jacobi method with different weights when solving the discrete system corresponding to $a$ depicted in Figure \ref{fig:a}. It can be observed that the damped Jacobi method exhibits behavior similar to that observed in the Poisson equation. Let $\Theta^{\mathcal{H}} = [-\pi/2, \pi/2)^2$ and $\Theta^{\cu{B}} = [-\pi, \pi)^2 \setminus \Theta^{\mathcal{H}}$. We choose $\cu{B}$ as the damped Jacobi method with $\omega = 3/4$ and set $M=10$, which satisfies Assumption \ref{ass:B}.

\begin{figure}[!htbp]
    \centering
    \includegraphics[width=0.5\textwidth]{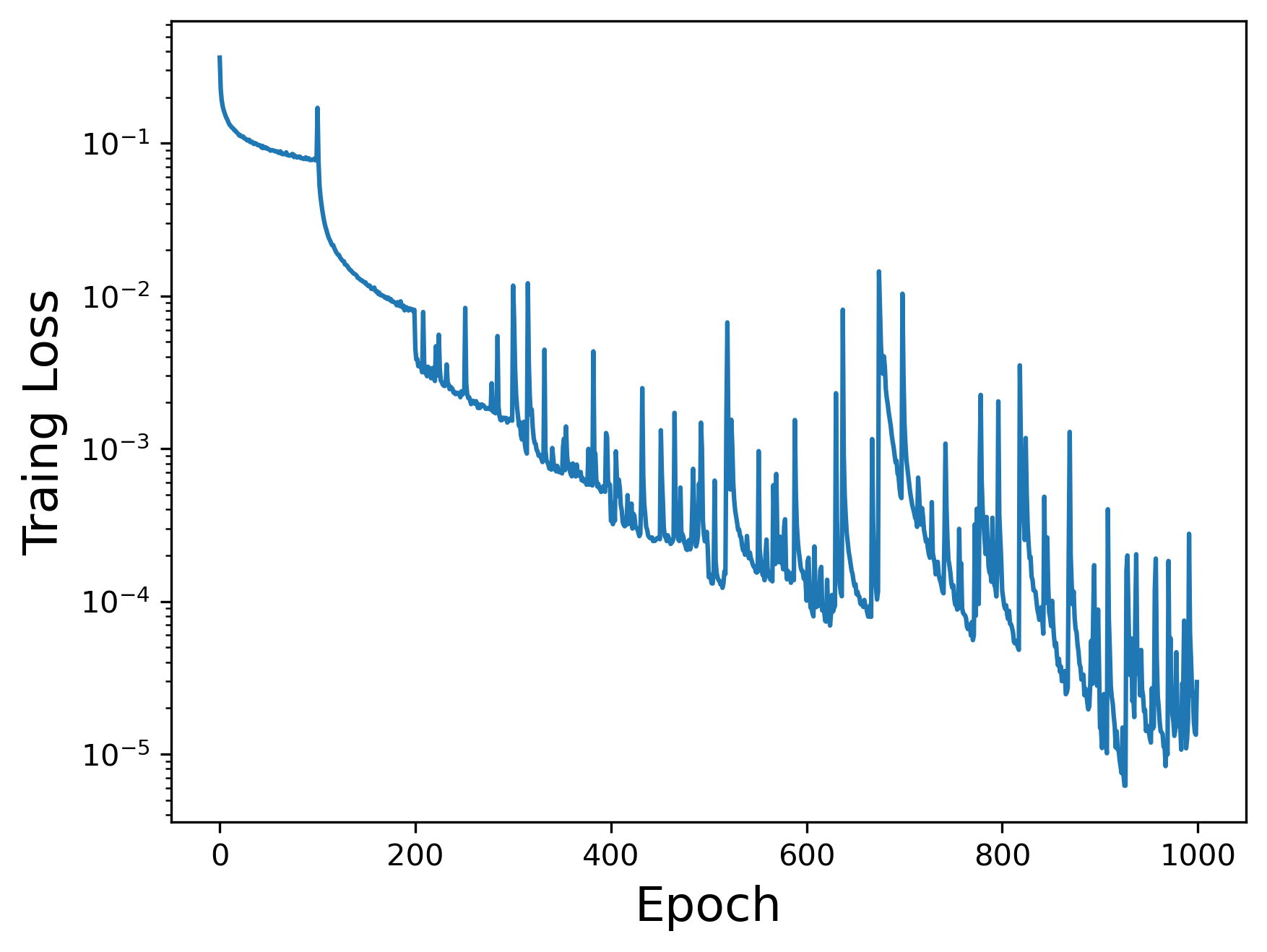}
    \caption{Training loss of FNS when solving random diffusion equations.}
    \label{fig:train_loss}
\end{figure}

Next, we train FNS to ensure $\mathcal{H}$ satisfies assumption \eqref{eq:assl} as much as possible. For this problem, we aim for the trained FNS to generalize well for both $a$ and $f$. According to previous analysis, the asymptotic convergence rate of FNS is independent of RHS (as we will verify below), so $a$ is the only parameter of interest, corresponding to $\boldsymbol{\mu}$ in PPDE \eqref{ppde}.
We generate 10,000 random diffusion coefficients $a_i$ and sample a random RHS $\cu{f_i} \sim \mathcal{N}(0, \cu{I})$ for each parameter $a_i$. We employ the loss function \eqref{eq:loss_func} for unsupervised training,  utilize the Adam optimizer with an initial learning rate of $10^{-4}$ and a learning rate schedule that halves the learning rate every 100 epochs. Additionally, we increase $K$ by one every 100 epochs. The training loss is depicted in Figure \ref{fig:train_loss}.

\begin{figure}[!htbp]
    \centering
    \includegraphics[width=\textwidth]{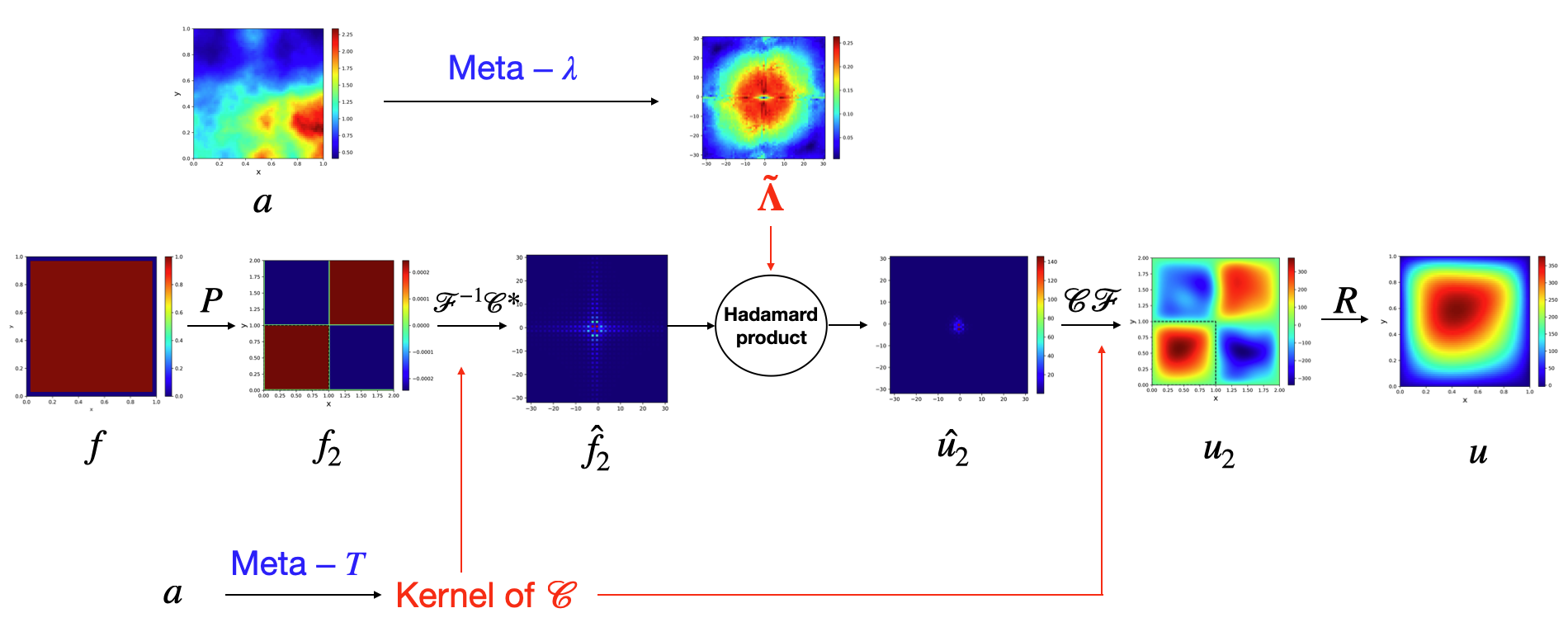}
    \caption{Calculation flow of $\mathcal{H}$ when solving the random diffusion equation.}
    \label{fig:Darcy_flow}
\end{figure}

Now we test the trained FNS. Figure \ref{fig:Darcy_flow} illustrates the calculation flow of $\mathcal{H}$ when it receives an input pair $(a, f)$. It can be observed that the $\cu{\tilde{\Lambda}}$ given by Meta$-\lambda$ is large in $\Theta^{\mathcal{H}}$ and small in $\Theta^{\cu{B}}$, which aligns with our expectations. By inputting the RHS into $\mathcal{H}$, we obtain an initial value close to the reference solution, indicating that $\mathcal{H}$ effectively captures the low-frequency components of the solution.

We then evaluated the performance of FNS across different discretization scales. For testing, we sampled 10 distinct coefficients $a_i$ and fixed the RHS $f = 1$. Table~\ref{tab:random_iters} presents the mean and standard deviation of the iteration counts required by FNS (trained at specific scales) to reduce the relative residual below $10^{-6}$. The results reveal the following:
\begin{enumerate}
\item FNS achieves optimal performance at the training scale, with iteration counts increasing at both smaller and larger scales.
\item FNS exhibits parameter-independent convergence that is only weakly dependent on discretization scale within a certain range near the training scale.
\item FNS generalizes better to scales smaller than the training scale than to significantly larger scales.
\end{enumerate}
These observations confirm that while FNS achieves scale-independent convergence near the training regime, it does not yet demonstrate universal scale invariance. At present, this issue can be mitigated by training on larger-scale problems. Developing more cost-effective approaches to enhance scale generalization will be a focus of our future work.

\begin{table}[!htb]
    \centering
    \caption{FNS iteration counts for FEM discretizations of random diffusion equations at different grid size. Each entry represents the mean $\pm$ std over 10 random coefficient samples, with a fixed RHS $f = 1$.}
    \label{tab:random_iters}
    \footnotesize{
    \begin{tabular}{lcccccc}
    \toprule
    Grid size $n$  & 31 & 63 & 127 & 255 & 511  \\ \midrule
    Model trained at $n=63$  & $15.6 \pm 1.74$ & $\mathbf{12.3 \pm 0.78}$ & $13.7\pm 0.78$ & $25.0 \pm 1.34$ & $45.2 \pm 3.67$ \\
    Model trained at $n=255$ & $15.0 \pm 2.88$ & $14.3 \pm 2.05$ & $12.6 \pm 1.43$ & $\mathbf{12.5 \pm 1.43}$ & $15.9 \pm 1.14$ \\
    \bottomrule
    \end{tabular}}
\end{table}

\begin{remark}[Generalization of FNS on different RHS]
    When using iterative methods, we often start with a zero initial guess. In this case, the frequencies present in the initial error are determined by $\cu{u}$, and those in the initial residual are determined by $\cu{f}$.
    Below, we compare the behavior of FNS when solving four different RHS:
    \begin{enumerate}[label=(\arabic*)]
        \item Single frequency function: 
        $$
        f_1(x,y)=\sin(\pi x)\sin(3\pi y),
        $$
        \item Gaussian function: 
        $$
        f_2(x,y)=\exp(-200((x-0.6)^2+(y-0.55)^2)),
        $$
        \item Constant function: 
        $$
        f_3(x,y)=1,
        $$
        \item Random function: 
        $$
        f_4=Au, \quad u\sim \mathcal{N}(0,1).
        $$
    \end{enumerate}
    Figure \ref{fig:res} shows the comparison results. It can be seen that although the error reduction of the first step differs due to the different RHS, the asymptotic convergence rate remains almost the same.
    \begin{figure}[!htbp]
        \centering
        \subfigure[$n=31$]{\includegraphics[width=0.3\textwidth]{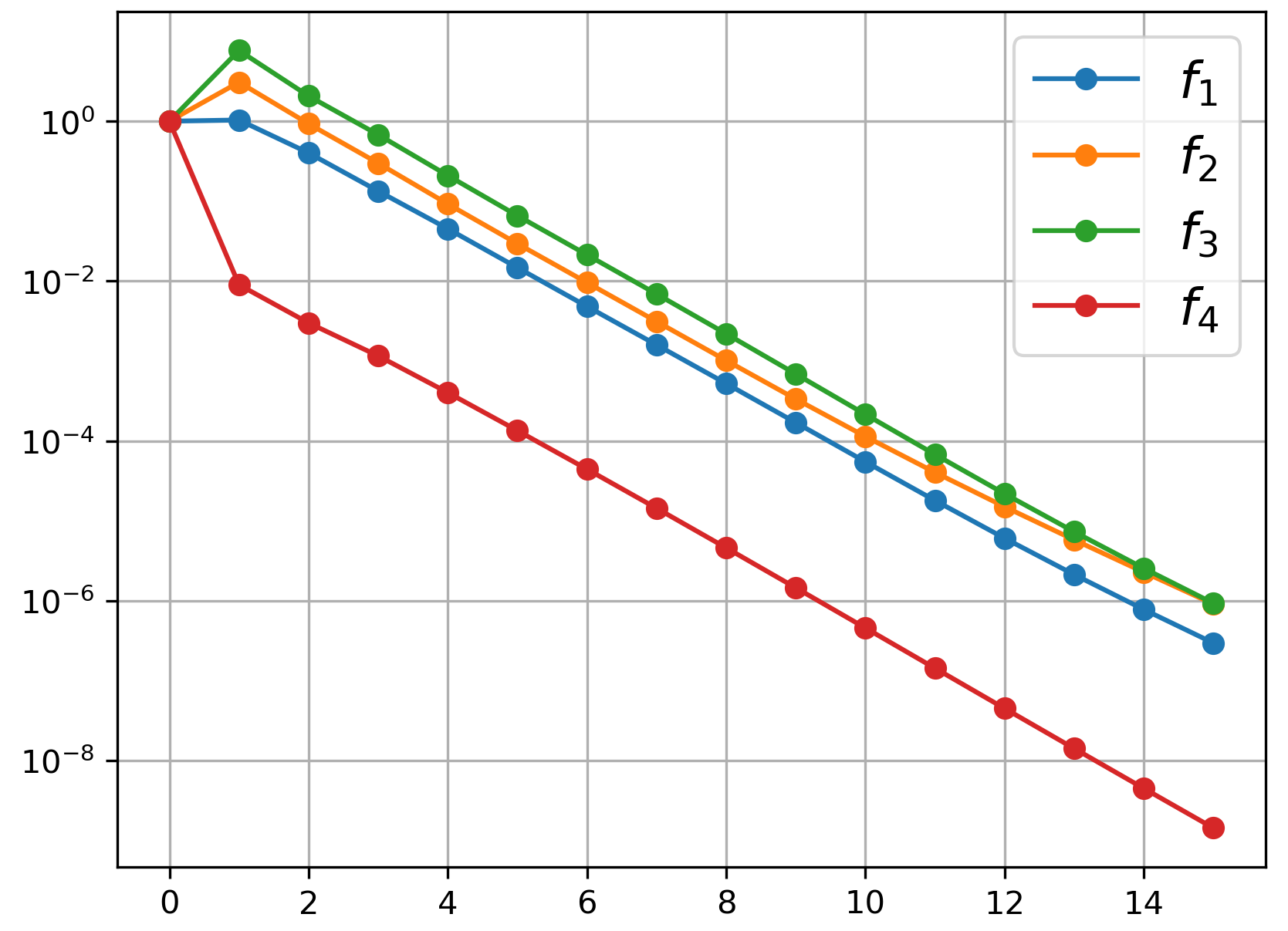}}\quad
        \subfigure[$n=63$]{\includegraphics[width=0.3\textwidth]{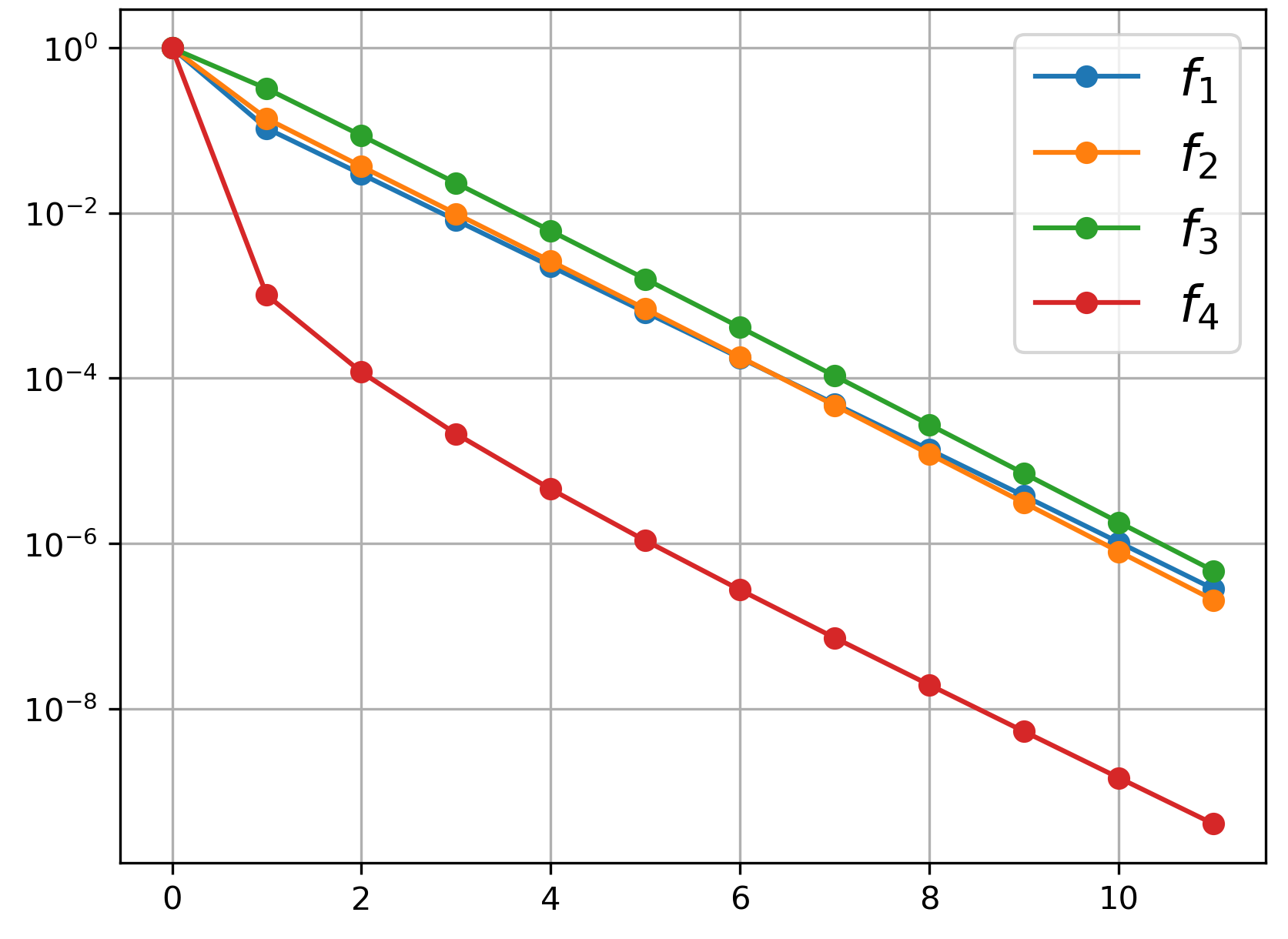}}\quad
        \subfigure[$n=127$]{\includegraphics[width=0.3\textwidth]{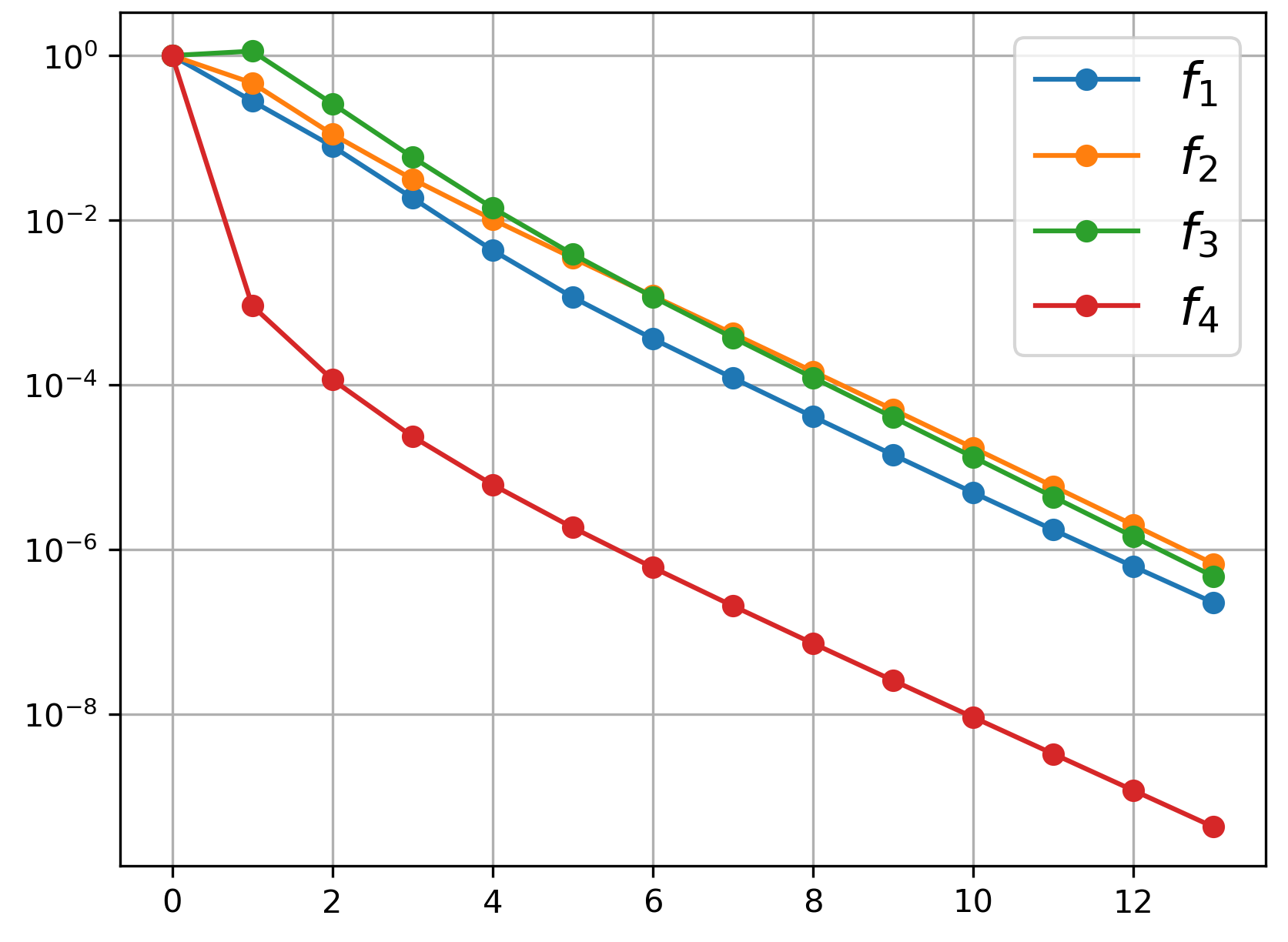}}
        \caption{Convergence history of FNS for solving random diffusion equations with different RHS.}
        \label{fig:res}
    \end{figure}
\end{remark}

\subsection{Anisotropic diffusion equations}
Consider the following 2D anisotropic diffusion equation
\begin{equation}
    \left\{\begin{aligned}
    -\nabla \cdot(C \nabla u) &= f, \quad \cu{x} \in \Omega,\\
    u &= 0, \quad \cu{x} \in \partial \Omega,
    \end{aligned}\right.
    \label{eq:anistropy}
\end{equation}
where the diffusion coefficient is a constant matrix
\begin{equation}
    C = \left(\begin{array}{ll}
    c_1 & c_2 \\
    c_3 & c_4
    \end{array}\right)
    = \left(\begin{array}{cc}
    \cos \theta & -\sin \theta \\
    \sin \theta & \cos \theta
    \end{array}\right)
    \left(\begin{array}{ll}
    1 & 0 \\
    0 & \xi
    \end{array}\right)
    \left(\begin{array}{cc}
    \cos \theta & \sin \theta \\
    -\sin \theta & \cos \theta
    \end{array}\right),
\end{equation}
$0 < \xi < 1$ is the anisotropic strength, and $\theta \in [0, \pi]$ is the anisotropic direction, $\Omega=(0,1)^2$. Using bilinear FEM on a uniform rectangular mesh with size $h = 1/(n+1)$, the resulting stencil is
\begin{equation}
    c_1\left[\begin{array}{ccc}
    -\frac{1}{6} & \frac{1}{3} & -\frac{1}{6} \\
    -\frac{2}{3} & \frac{4}{3} & -\frac{2}{3} \\
    -\frac{1}{6} & \frac{1}{3} & -\frac{1}{6}
    \end{array}\right] + (c_2+c_3)\left[\begin{array}{ccc}
    -\frac{1}{4} & 0 & \frac{1}{4} \\
    0 & 0 & 0 \\
    \frac{1}{4} & 0 & -\frac{1}{4}
    \end{array}\right] + c_4\left[\begin{array}{ccc}
    -\frac{1}{6} & -\frac{2}{3} & -\frac{1}{6} \\
    \frac{1}{3} & \frac{4}{3} & \frac{1}{3} \\
    -\frac{1}{6} & -\frac{2}{3} & -\frac{1}{6}
    \end{array}\right].
    \label{eq:anistencil}
\end{equation}

\begin{figure}[!htbp]
    \centering
    \subfigure[$\omega=1/2, M=1$]{\label{fig:ani1}
    \includegraphics[width=0.3\textwidth]{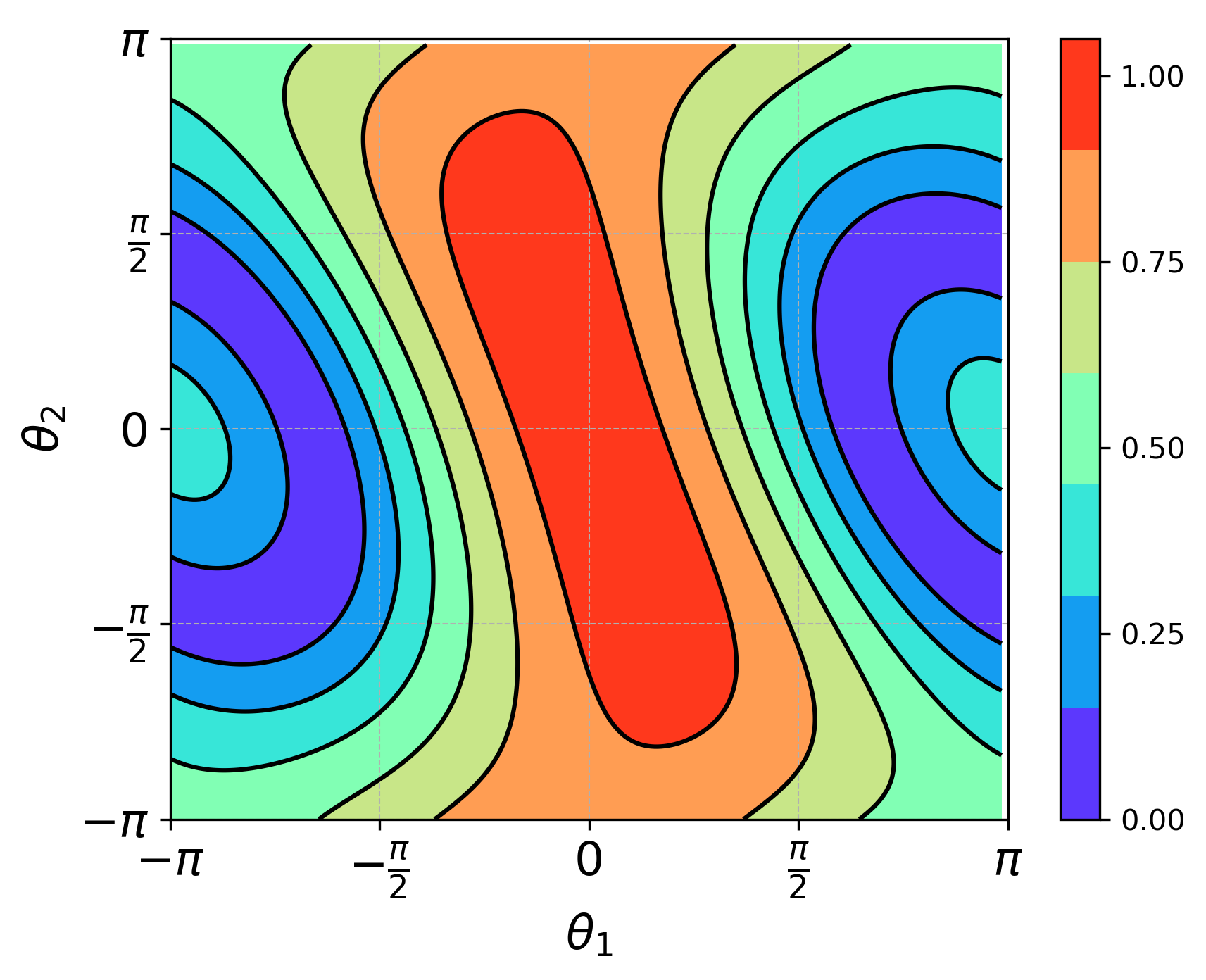}}\quad
    \subfigure[$\omega=1/2, M=5$]{\label{fig:ani2}
        \includegraphics[width=0.3\textwidth]{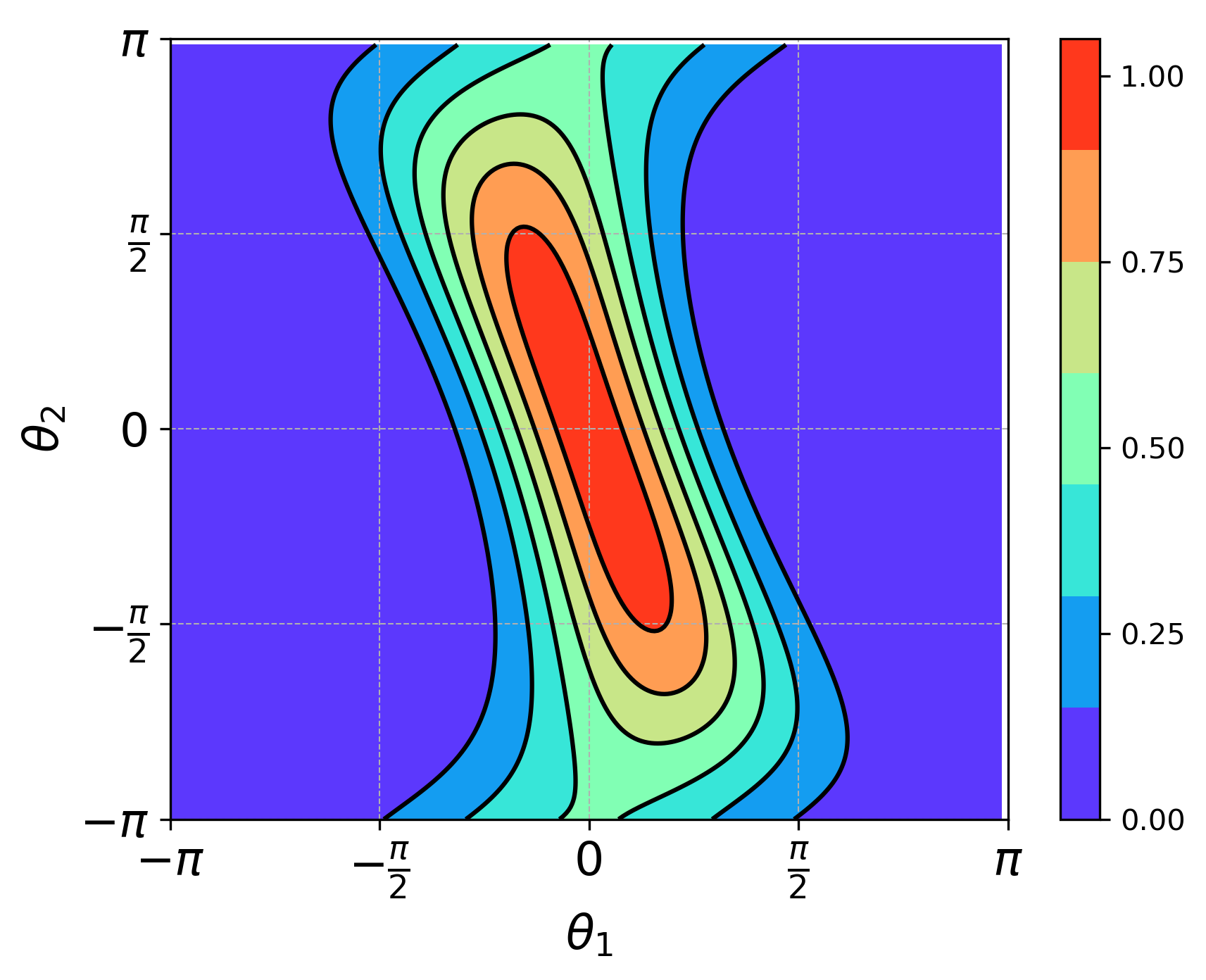}}
    \subfigure[$\omega=2/3, M=5$]{\label{fig:ani3}
        \includegraphics[width=0.3\textwidth]{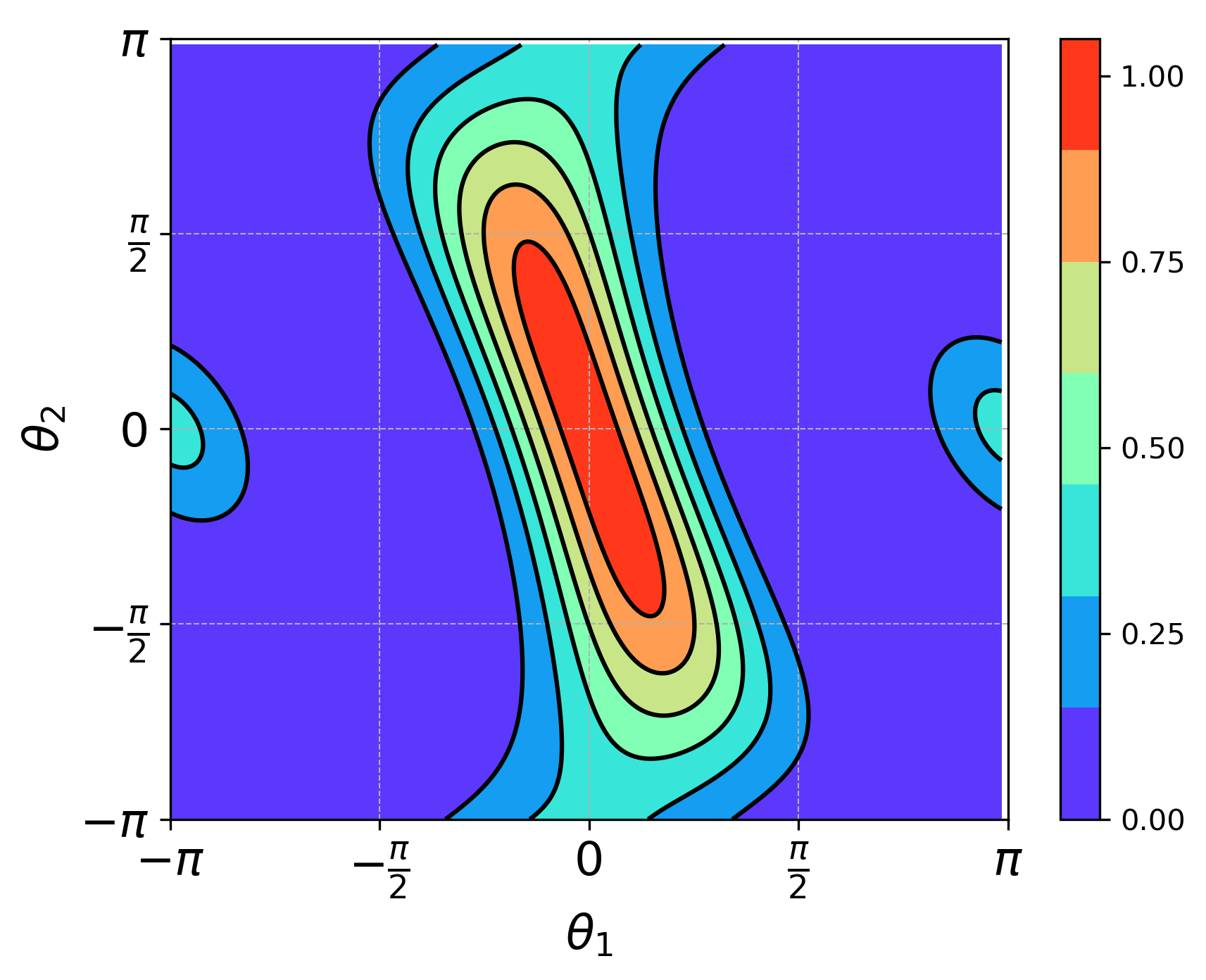}}
    \caption{Jacobi symbol applied to anisotropic equation with  $\xi = 10^{-6}$, $\theta=0.1\pi$}
    \label{fig:ani_LFA}
\end{figure}

Selecting $\cu{B}$ as the damped Jacobi method, Figure \ref{fig:ani_LFA} shows the Jacobi symbol when solving the anisotropic discrete system corresponding to $\xi = 10^{-6}$, $\theta = 0.1\pi$. 
It can be seen that when the damped Jacobi method uses $\omega = 1/2$ and $M = 1$, error components with frequencies along the anisotropic direction are difficult to eliminate.
We denote the frequency interval along the anisotropic direction as $\Theta^{\mathcal{H}}$ and $\Theta^{\cu{B}} = [-\pi, \pi)^2 \setminus \Theta^{\mathcal{H}}$. 
To enhance the compressibility of error components with frequencies in $\Theta^{\cu{B}}$, we increase $M$ to 5. Figure \ref{fig:ani2} shows the distribution of the corresponding symbol, which has improved a lot.
The reason for using $\omega = 1/2$ is that $\Theta^{\cu{B}}$ remains connected under this weight. If $\omega$ is further increased, such as to $\omega = 2/3$, $\Theta^{\cu{B}}$ will appear in other parts of the domain, which will increase the difficulty of learning  for $\mathcal{H}$.

Next, we train FNS to ensure that $\mathcal{H}$ satisfies Assumption \ref{ass:H}. Corresponding to PPDE \eqref{ppde}, $\cu{\mu} = (\xi, \theta)$ in this example. We sample 500 sets of parameters $\cu{\mu}_i$ in the parameter interval $[10^{-6}, 1] \times [-\pi, \pi]$ according to the uniform distribution. For each parameter $\cu{\mu}_i$, we sample 20 RHS $\cu{f}_i \sim \mathcal{N}(0, \cu{I})$, and use the loss function \eqref{eq:loss_func} for unsupervised training. The hyperparameters used during training are the same as those used in the random diffusion equations.

\begin{figure}[!htbp]
    \centering
    \includegraphics[width=\textwidth]{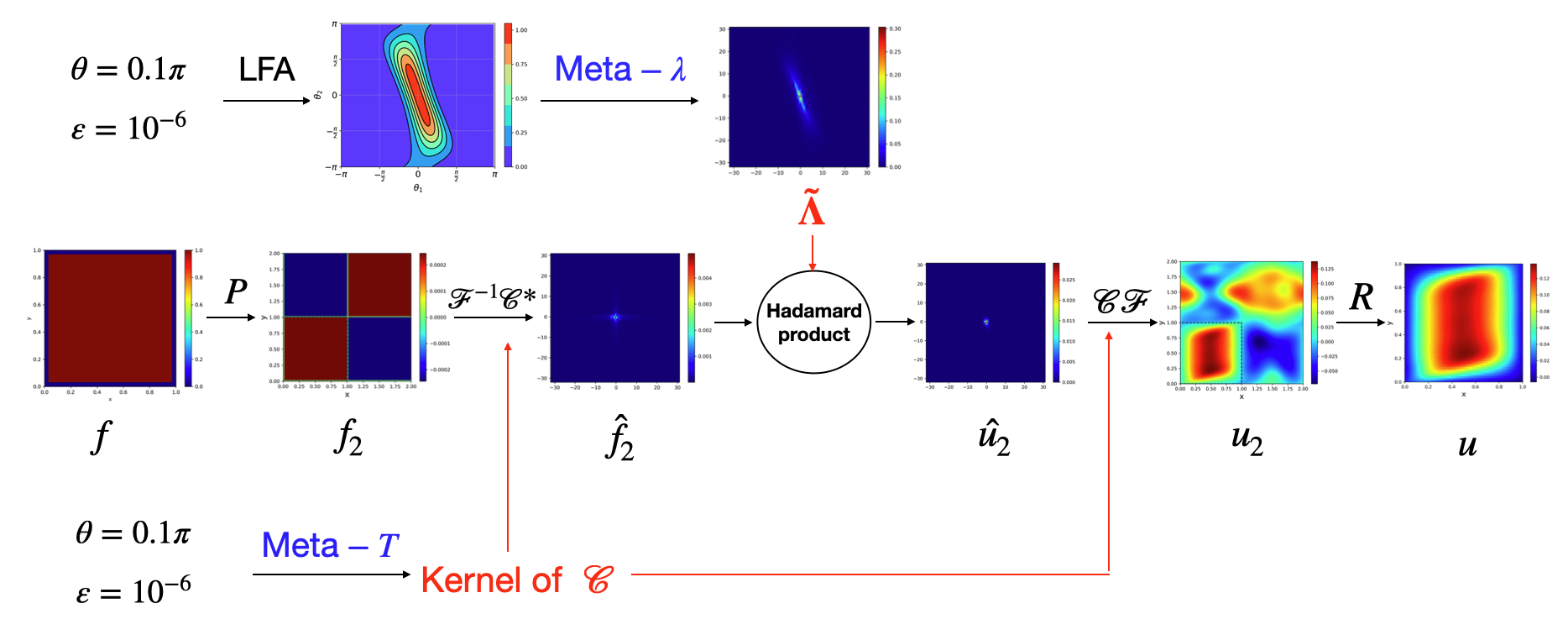}
    \caption{Calculation flow of $\mathcal{H}$ when solving the anisotropic diffusion equation.}
    \label{fig:Ani}
\end{figure}

Once trained, we test the performance of FNS on newly selected parameters and different scales. Figure \ref{fig:Ani} shows the calculation flow of $\mathcal{H}$ when receiving specific parameters $\cu{\mu_i}$ and $\cu{f_i}$. The difference compared to the random diffusion equation is that, instead of directly inputting the parameter $\cu{\mu}$ into the Meta$-\lambda$, we use the Jacobi symbol obtained by LFA as input. This is because the symbol contains more intuitive information than $\cu{\mu}$ and is easier to learn. From Figure \ref{fig:Ani}, it can be seen that the $\cu{\tilde{\Lambda}}$ given by Meta$-\lambda$ is large in $\Theta^{\mathcal{H}}$ but small in $\Theta^{\cu{B}}$, which is consistent with our expectations.

We next evaluate the performance of FNS in solving anisotropic diffusion equations with varying parameters $\cu{\mu}$ across different mesh resolutions $h$. FNS was trained separately at $63^2$ and $511^2$ grids. Due to the high training cost at $511^2$, the training dataset for this case was limited to 2,000 samples (200 $\cu{\mu}_i$, each paired with 10 different RHS $\cu{f}_i$).
Table~\ref{tab:ani_iters} reports the iteration counts of FNS trained at different resolutions. It can be observed that:
\begin{enumerate}
\item FNS trained at $n=63$ exhibits mesh- and parameter-independent convergence for test grids with $n \in [31, 127]$, but its performance deteriorates  on larger grids ($n \geq 255$).
\item FNS trained at $n=511$ achieves comparable convergence performance at its training scale to that of the $n=63$ model at its own scale.
\item The $n=511$-trained FNS performs worse than the $n=63$ model on smaller grids, which may be attributed to limited training data  at large scales.
\end{enumerate}
\begin{table}[!htb]
\centering
\caption{FNS iteration counts for anisotropic diffusion equations across grid sizes. Each entry reports the mean $\pm$ standard deviation over 10 random coefficients with $\log_{10} \varepsilon \sim \mathcal{U}(-6,0)$, $\theta = \frac{5\pi}{12}$, and  $f = 1$.}
\label{tab:ani_iters}
\footnotesize{
\begin{tabular}{lcccccc}
\toprule
Grid size $n$  & 31 & 63 & 127 & 255 & 511  \\ \midrule
Model trained at $n=63$  & $22.7 \pm 2.45$ & $\mathbf{21.5 \pm 1.97}$ & $24.4 \pm 2.96$ & $78.5 \pm 10.86$ & $>200$ \\
Model trained at $n=511$ & $78.3 \pm 9.53$ & $99.5 \pm 8.02$ & $83.6 \pm 5.16$ & $78.0 \pm 10.69$ & $\mathbf{24.5 \pm 6.44}$ \\
\bottomrule
\end{tabular}}
\end{table}

\begin{figure}[!htbp]
\centering
\subfigure[Random diffusion]{\label{fig:hints_random}
\includegraphics[width=0.45\textwidth]{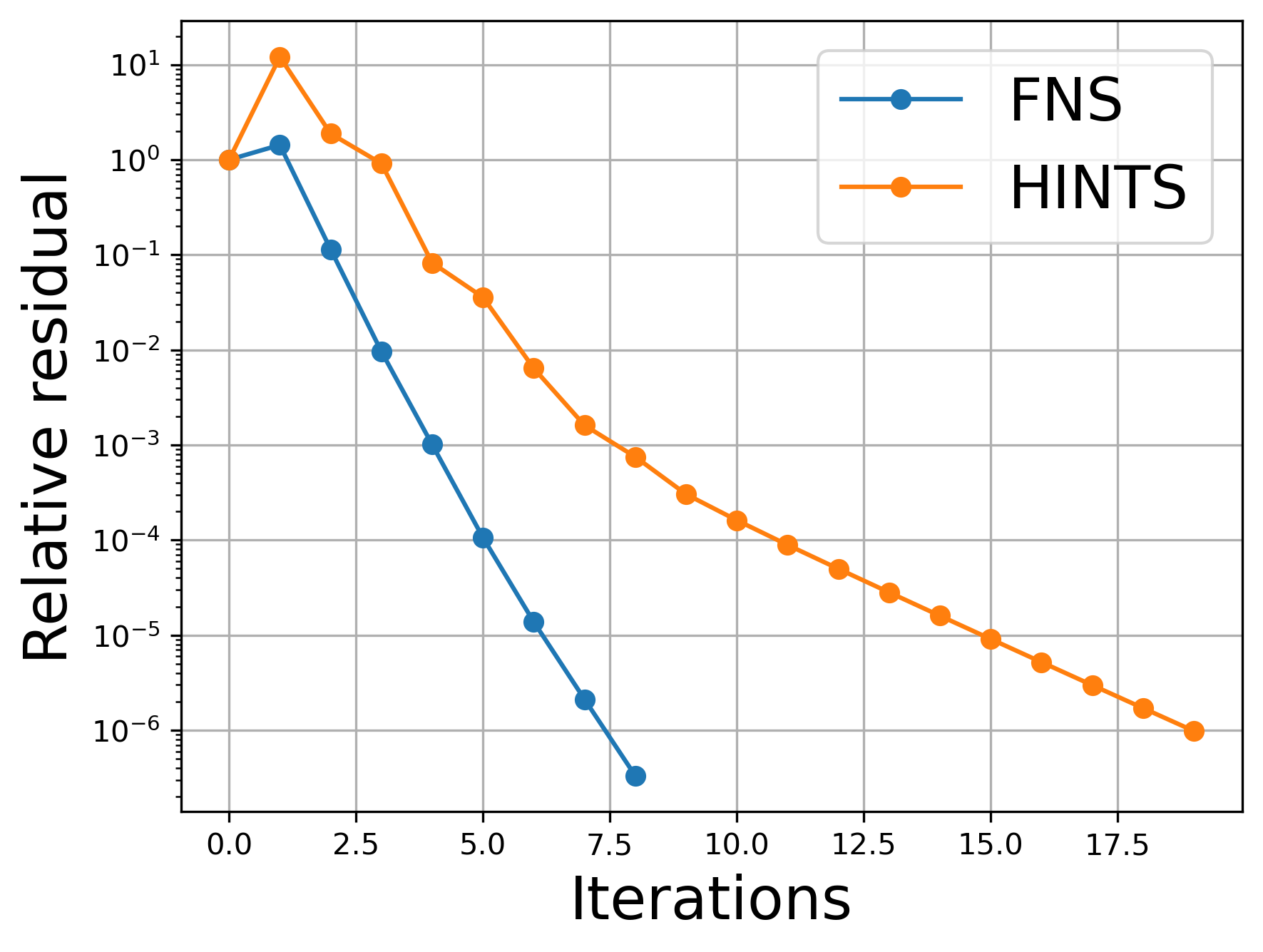}}
\subfigure[Anisotropic diffusion ($\varepsilon = 10^{-6}, \theta=\pi/4$)]{\label{fig:hints_ani}
\includegraphics[width=0.45\textwidth]{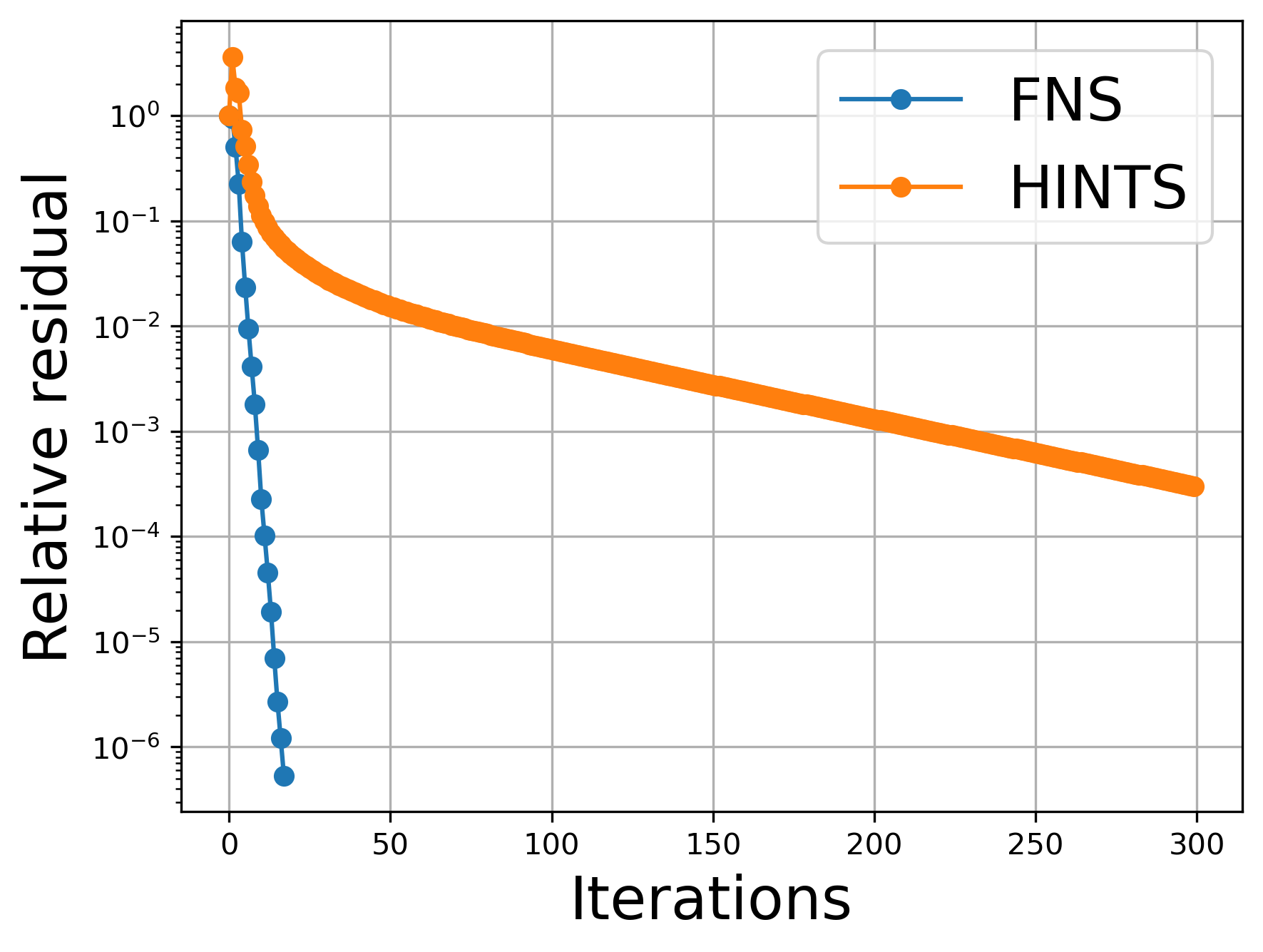}}
\caption{Convergence histories of FNS and HINTS for solving random and anisotropic diffusion equations, with $n=63$.}
\end{figure}

We then compare the performance of FNS with HINTS, which employs DeepONet as  $\mathcal{H}$. For the random diffusion equation \eqref{eq:random}, the damped Jacobi smoother $\cu{B}$  effectively removes high-frequency errors, while DeepONet efficiently learns the remaining low-frequency components due to the  spectral bias. As a result, HINTS achieves rapid convergence on small-scale problems, as shown in Figure~\ref{fig:hints_random}.
However, for anisotropic diffusion equations, where $\cu{B}$ struggles to eliminate high-frequency errors aligned with the anisotropy direction, DeepONet also fails to learn such errors effectively, again due to the spectral bias. In this case, HINTS converges slowly, as illustrated in Figure~\ref{fig:hints_ani}. 
In contrast, FNS effectively mitigates spectral bias by operating in the frequency domain and adopting an end-to-end training strategy. This design enables the network to better capture a broader range of error components, including those that are typically challenging for conventional neural solvers, thereby achieving more robust and efficient convergence.

\subsection{Convection-diffusion equations}\label{sec:convection}
Consider the following 2D convection-diffusion equation
\begin{equation}
    \begin{aligned} 
    -\varepsilon \Delta u + \vec{w} \cdot \nabla u &= f, \quad \cu{x} \in \Omega, \\
    u &= g, \quad \cu{x} \in \partial \Omega, \\
    \end{aligned}
\end{equation}
where $\varepsilon > 0$ is a viscosity parameter, and $\vec{w} = (w_x(x, y), w_y(x, y))$ is the flow velocity, assumed to be incompressible ($\mathrm{div}\, \vec{w} = 0$).
We are interested in the convection-dominated case, \ie, $\varepsilon \ll |\vec{w}|$. In this setting, the solution typically has steep gradients in some parts of the domain. If we apply the usual Galerkin FEM, a sharply oscillating solution will be obtained.
To avoid this problem, we use the streamline diffusion FEM, which finds $ u_h \in U_h $ such that
\begin{equation}
    \begin{aligned}
    a_h(u_h, v_h) &= \varepsilon (\nabla u_h, \nabla v_h) + (\vec{w} \cdot \nabla u_h, v_h) 
    + \sum_{K \in \Omega_h} \delta_K ( \vec{w} \cdot \nabla u_h, \vec{w} \cdot \nabla v_h)_K \\&= (f,v_h) + \sum_{K \in \Omega_h} \delta_K ( f, \vec{w} \cdot \nabla v_h)_K
    \quad \forall v_h\, \in V_h,
    \end{aligned}
\end{equation}
where $ \delta_K $ is a user-chosen non-negative stabilization parameter, which plays a key role in the accuracy of the numerical solution.
As shown in \cite{elman2014finite}, a good choice of $ \delta_K $ is
\begin{equation}
    \delta_K = 
    \begin{cases} 
    \frac{h_K}{2|\vec{w}_K|}\left(1 - \frac{1}{\mathcal{P}_h^K}\right) & \text{if } \mathcal{P}_h^K > 1, \\
    0 & \text{if } \mathcal{P}_h^K \leq 1,
    \end{cases}
    \label{eq:delta}
\end{equation}
where $\mathcal{P}_h^K = \frac{|\vec{w}_K| h_K}{2 \varepsilon}$
is the element Péclet number.
In the following, we consider steady flow and use bilinear finite element on a uniform grid with spacing $ h = 1/(n+1) $. The corresponding stencil is given by
\begin{equation}
    \begin{aligned}
    & \frac{\varepsilon}{3}\left[\begin{array}{ccc}
    -1 & -1 & -1 \\
    -1 & 8 & -1 \\
    -1 & -1 & -1
    \end{array}\right]+\left[\begin{array}{ccc}
    -w_x+w_y & 4 w_y & w_x+w_y \\
    -4 w_x & 0 & 4 w_x \\
    -\left(w_x+w_y\right) & -4 w_y & w_x-w_y
    \end{array}\right] \\
    & +\delta\left[\begin{array}{ccc}
    -\frac{1}{6}\left(w_x^2+w_y^2\right)+\frac{1}{2}w_x w_y & \frac{1}{3} w_x^2-\frac{2}{3} w_y^2 & -\frac{1}{6}\left(w_x^2+w_y^2\right)-\frac{1}{2} w_x w_y \\
    -\frac{2}{3} w_x^2+\frac{1}{3} w_y^2 & \frac{4}{3}\left(w_x^2+w_y^2\right) & -\frac{2}{3} w_x^2+\frac{1}{3} w_y^2 \\
    -\frac{1}{6}\left(w_x^2+w_y^2\right)-\frac{1}{2} w_x w_y & \frac{1}{3} w_x^2-\frac{2}{3} w_y^2 & -\frac{1}{6}\left(w_x^2+w_y^2\right)+\frac{1}{2} w_x w_y
    \end{array}\right].
    \end{aligned}
    \label{eq:convstencil}
\end{equation}

\begin{figure}[!htb]
    \centering
    \subfigure[$\omega=1/2,M=1$]{
    \includegraphics[width=0.22\textwidth]{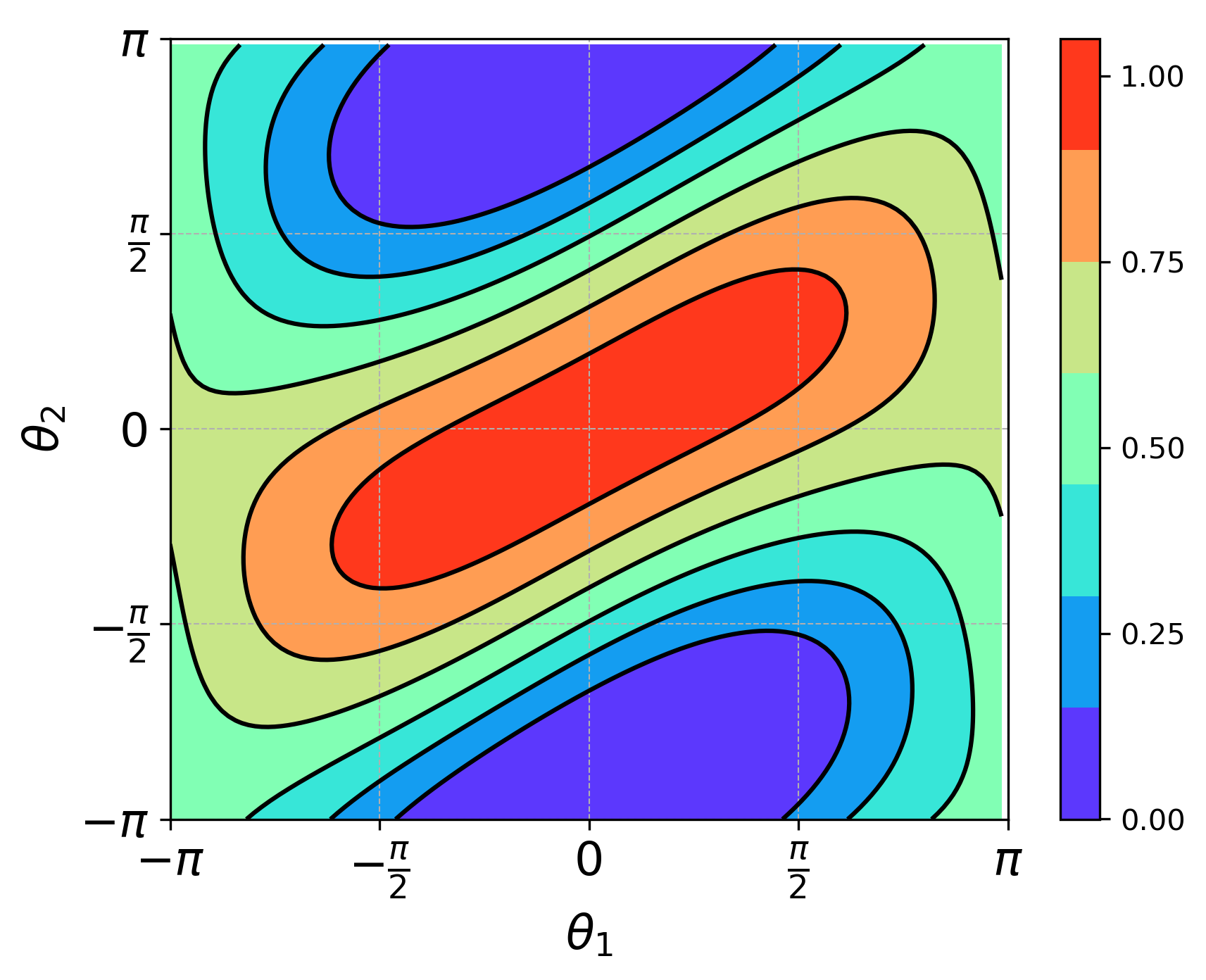}}
    \subfigure[$\omega=1/2,M=5$]{
    \includegraphics[width=0.22\textwidth]{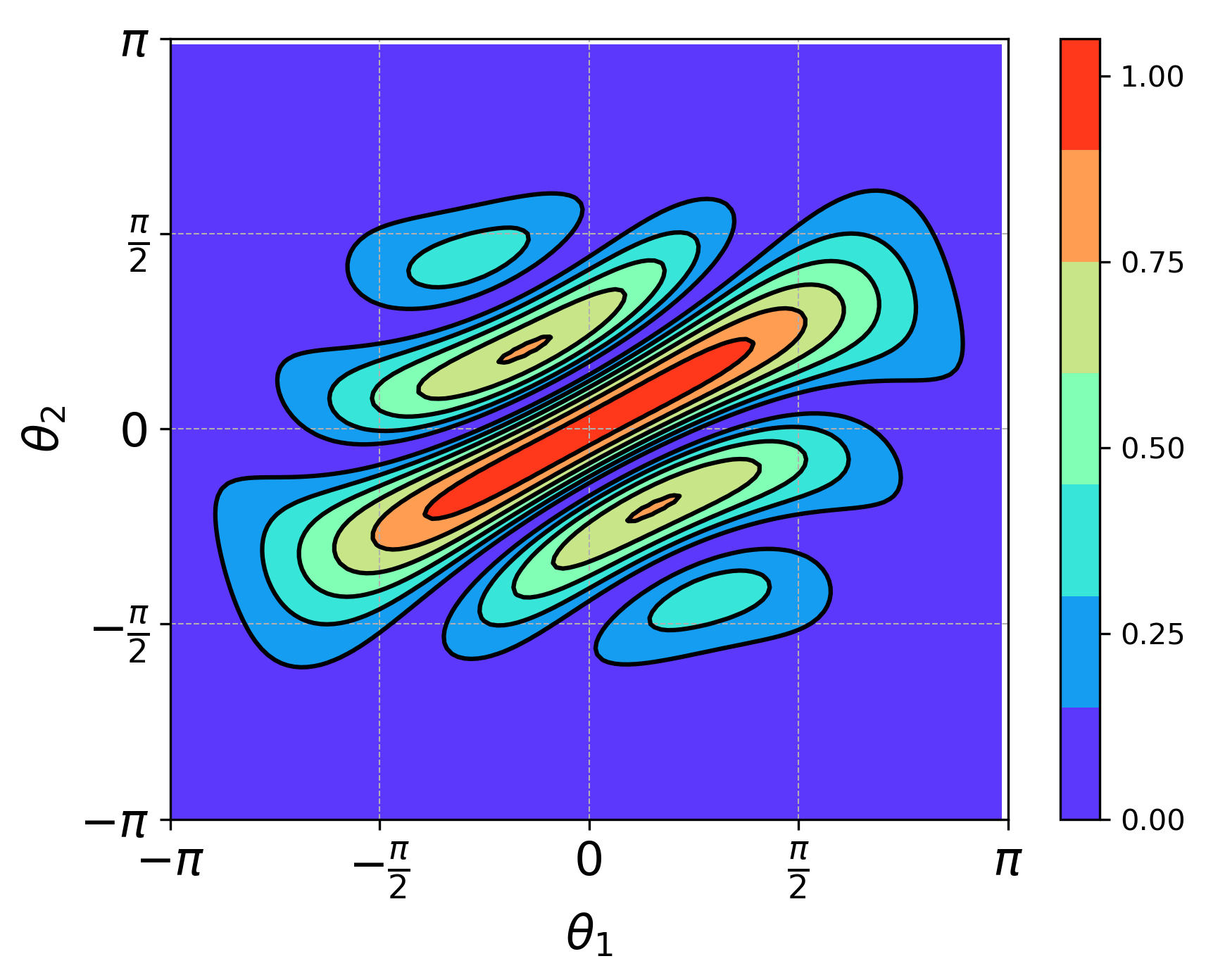}}
    \subfigure[$\omega=1/2,M=10$]{
    \includegraphics[width=0.22\textwidth]{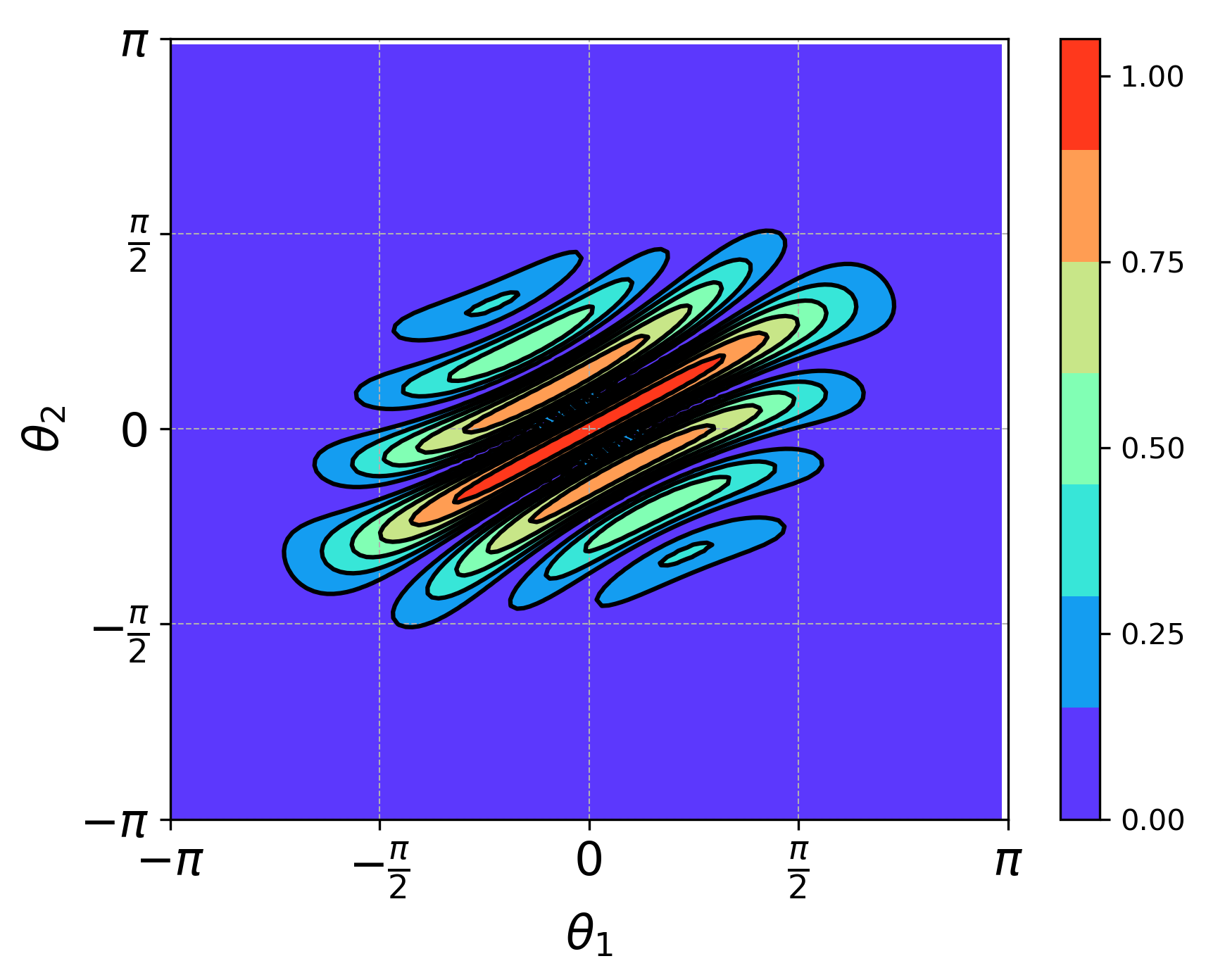}}
    \subfigure[$\omega=2/3,M=10$]{
    \includegraphics[width=0.22\textwidth]{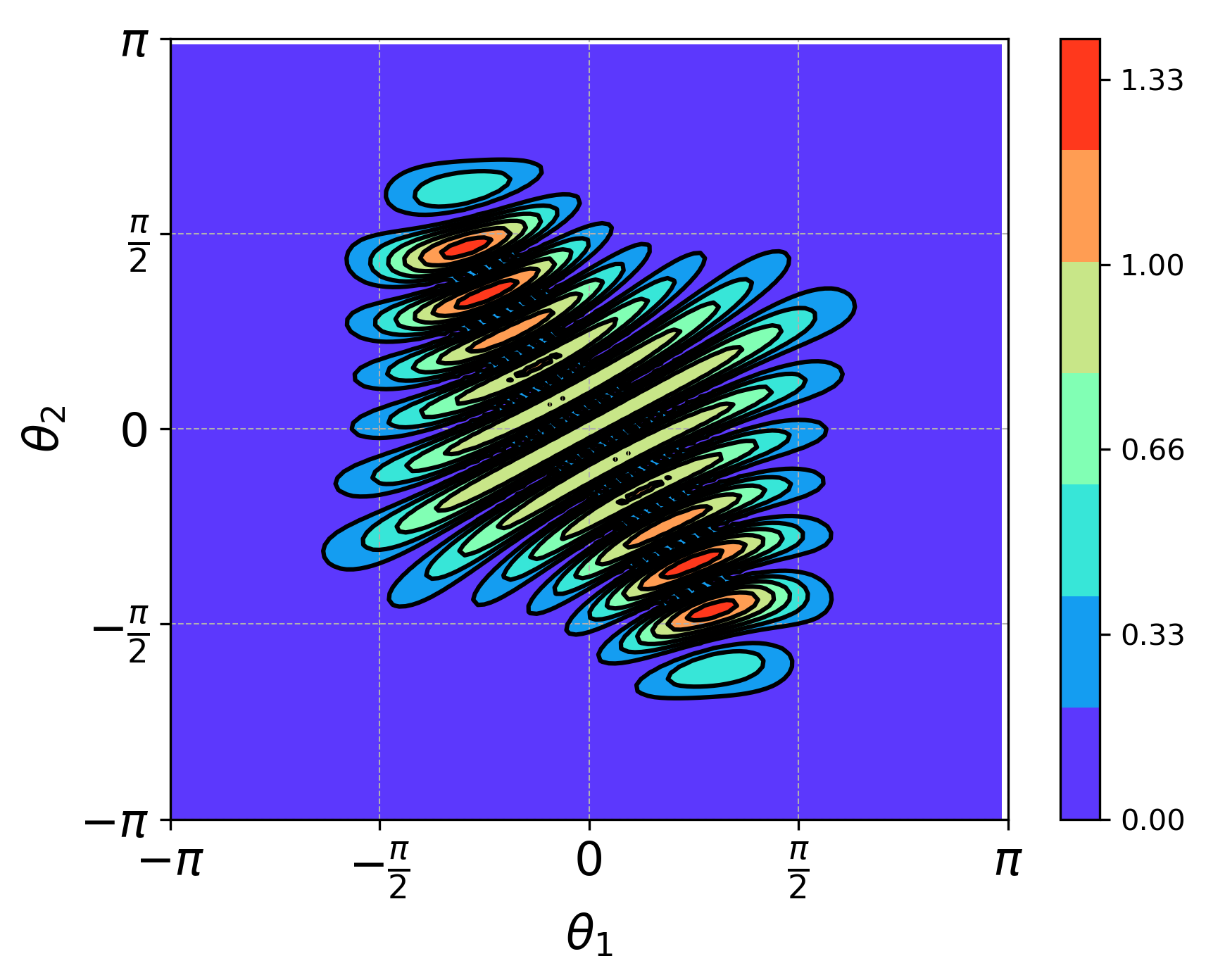}}
    \caption{Jacobi symbol applied to the convection-diffusion equation with $\varepsilon = 10^{-8}$, $w_x = -\sin(\pi/6)$, and $w_y = \cos(\pi/6)$.}
    \label{fig:lfa_convection}
\end{figure}

Selecting $\cu{B}$ as the damped Jacobi method, we estimate its compressibility for error components with different frequencies using LFA. Figure \ref{fig:lfa_convection} shows the Jacobi symbol for different weights and iteration counts when applied to the discrete system corresponding to $\varepsilon = 10^{-8}$, $w_x = -\sin(\pi/6)$, and $w_y = \cos(\pi/6)$.
The first three plots show that when $\omega = 1/2$, the damped Jacobi method has a better smoothing effect for error components with frequencies except along the streamline direction as $M$ increases, which means Assumption \ref{ass:B} can be met. However, when $\omega$ increases, such as to $\omega = 2/3$, the damped Jacobi method not only fails to achieve a better smoothing effect but also amplifies the error components with certain frequencies. Therefore, we take $\omega = 1/2$ and $M = 10$.

Next, we train FNS to make $\mathcal{H}$ satisfy \eqref{eq:assl}. Corresponding to the PPDE \eqref{ppde}, the parameters of this problem are $\boldsymbol{\mu} = (\varepsilon, \vec{w})$.
We sample 100 sets of parameters $\{\varepsilon, w_x, w_y\}$ as follows: $\log \frac{1}{\varepsilon} \sim U[0, 8]$ and $w_x, w_y \sim U[-1, 1]$. Taking $h = 1/64$ as the training size, we calculate the corresponding $\delta$ according to \eqref{eq:delta} and then obtain the corresponding stencil \eqref{eq:convstencil}. For each set of parameters $\boldsymbol{\mu}_i$, we randomly generate 100 RHS $\cu{f}_i$ according to the standard normal distribution, obtaining a total of 10,000 training data pairs.

\begin{figure}[!htbp]
    \centering
    \includegraphics[width=\textwidth]{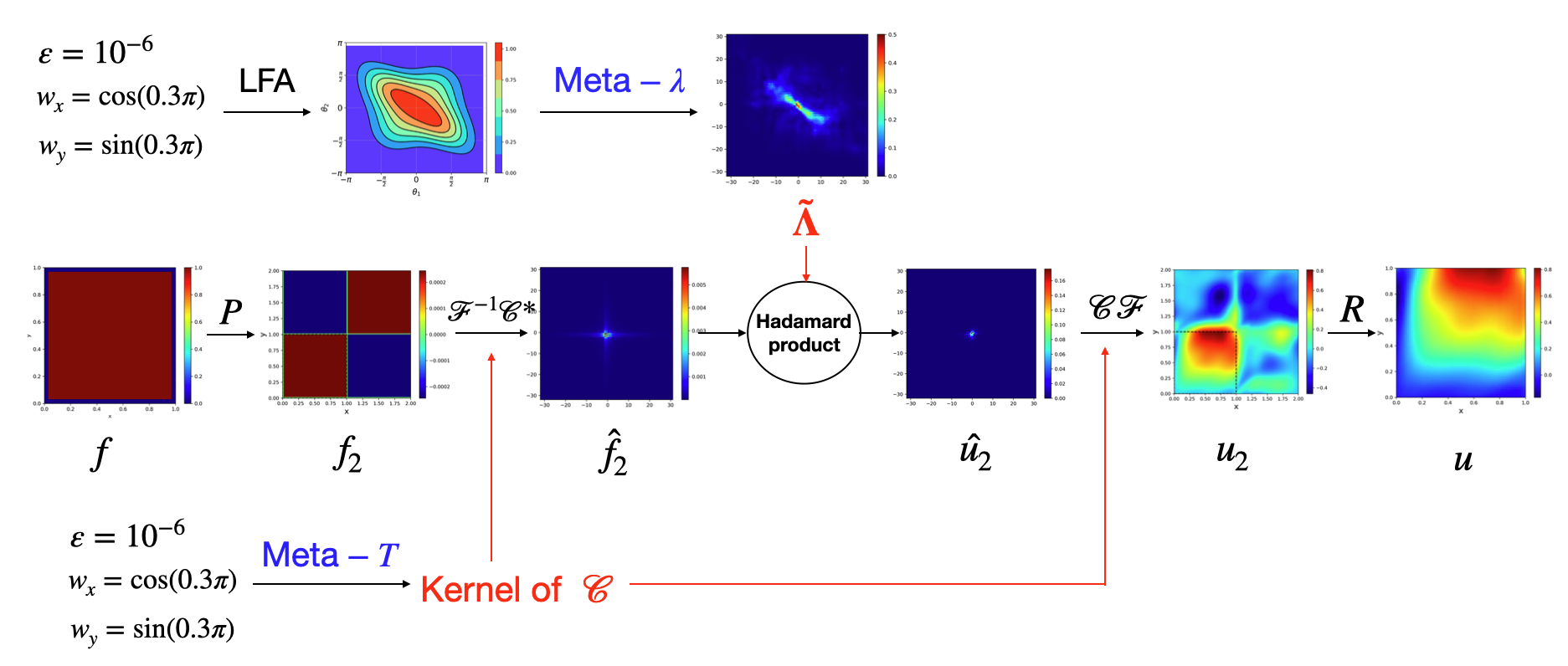}
    \caption{Calculation flow of $\mathcal{H}$ when solving the convection-diffusion equation.}
    \label{fig:convec_flow}
\end{figure}
After training, we test the performance of FNS. Figure \ref{fig:convec_flow} shows the calculation flow of $\mathcal{H}$ when it receives a pair of $\{\cu{\mu}_i, \cu{f}_i\}$. Since the smoothing effect of $\cu{B}$ on the convection-diffusion equation is similar to that of the anisotropic diffusion equation, we also first perform LFA on $\cu{B}$ to obtain its symbol, and then input the symbol into the Meta$-\lambda$. From Figure \ref{fig:convec_flow}, we can see that the $\cu{\tilde{\Lambda}}$ given by Meta$-\lambda$ is large in $\Theta^{\mathcal{H}}$ but small in $\Theta^{\cu{B}}$, which is consistent with our expectations.

Table \ref{tab:convec} shows the required iteration counts of FNS, which is trained on the scale $n=63$ but test on other scales. It can be seen that the convergence speed of FNS is weakly dependent on the grid size $h$.
\begin{table}[!htb]
    \centering
    \caption{FNS iteration counts to satisfy $\|\cu{r}^{(k)}\|/\|\cu{f}\|<10^{-6}$, starting with zero initial value, for convection-diffusion equation with $\varepsilon=10^{-7}, w_x=\cos(0.8\pi), w_y=\sin(0.8\pi)$ on different scales.}
    \label{tab:convec}
    \begin{tabular}{cccccc}
    \toprule
    $n$  & 31 & 63 & 127 & 255 & 511\\ \midrule
     iters & 8 & 12 & 16 & 19 & 23 \\ \bottomrule
    \end{tabular}
\end{table}

\subsection{Jumping diffusion equations} 
\begin{figure}[!htbp]
    \centering
    \subfigure[$a(\cu{x})$]{\label{fig:jumpa}
    \includegraphics[width=0.268\textwidth]{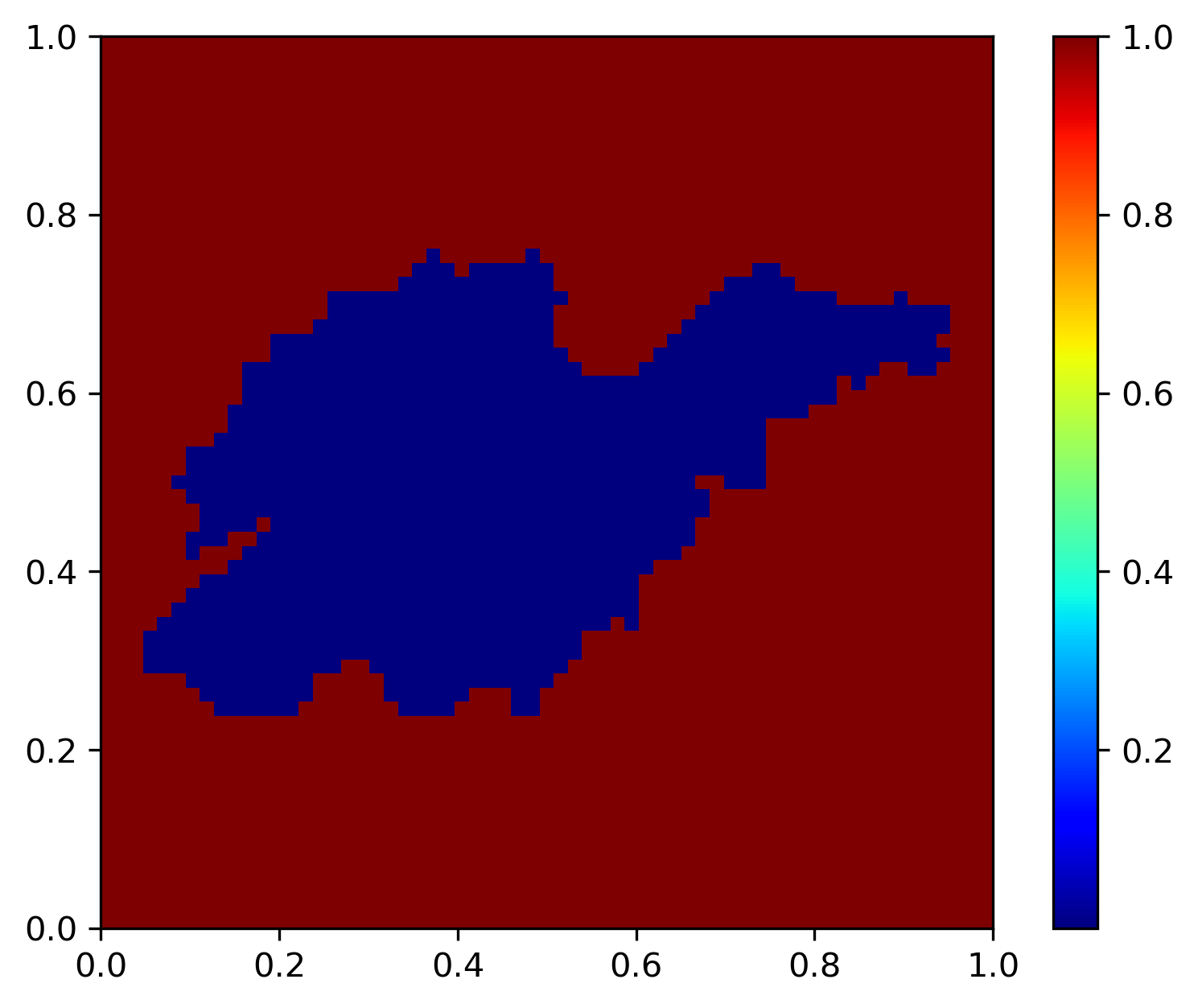}}
    \subfigure[$u(\cu{x})$]{\label{fig:jumpu}
    \includegraphics[width=0.268\textwidth]{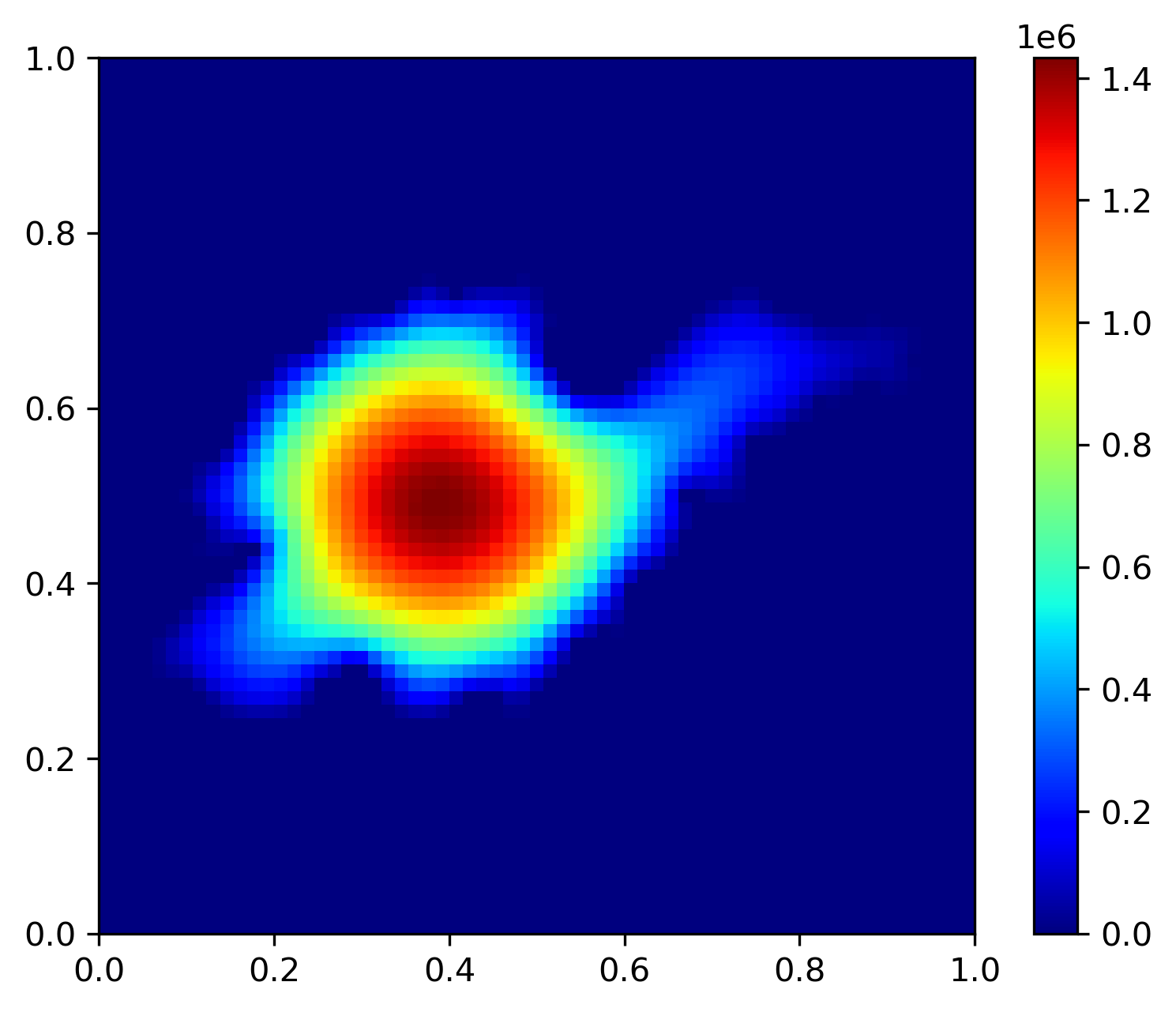}}
    \subfigure[Dual element and interface]{\label{fig:fvm_grid}
    \includegraphics[width=0.3\textwidth]{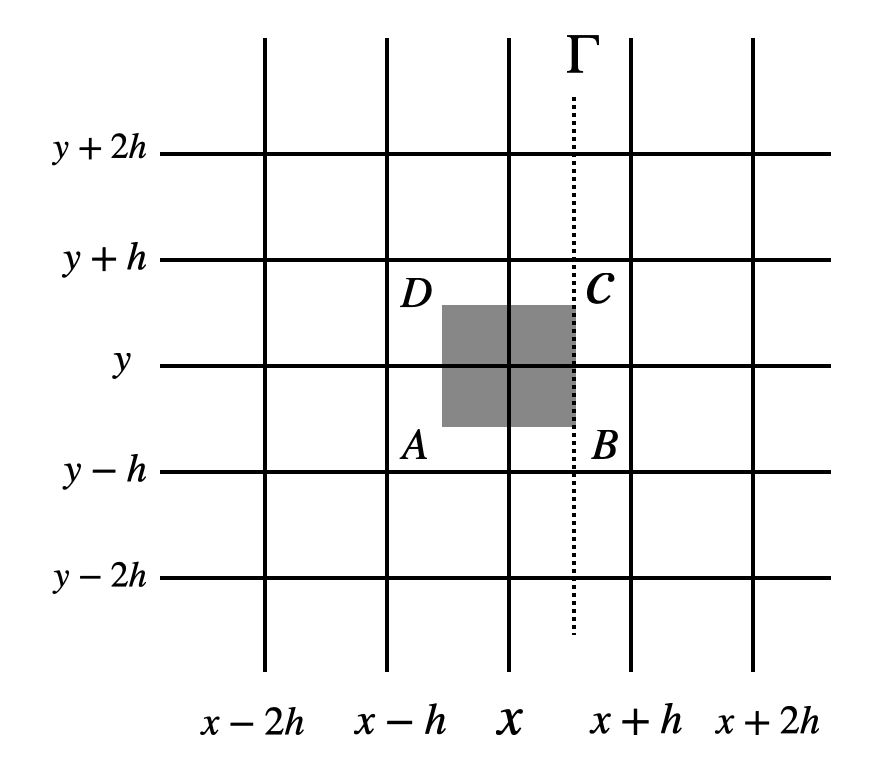}}
    \caption{An example of the jumping diffusion equation and discretization grid.}
\end{figure}

Consider the 2D diffusion equation with jumping coefficients
\begin{equation} 
    \left\{
    \begin{aligned}
    -\nabla \cdot(a(\cu{x}) \nabla u) & = f, \text{ in } \Omega, \\
    u & = 0, \text{ on } \partial \Omega, \\
    \llbracket u \rrbracket = 0, \llbracket a \nabla u \cdot \mathbf{n} \rrbracket & = 0, \text{ on } \Gamma,
    \end{aligned}
    \right.
    \label{eq:jumping} 
\end{equation} 
where
$$\Omega = \Omega_1 \cup \Omega_2, \quad \Omega_1 \cap \Omega_2 = \emptyset, \quad \Gamma = \partial \Omega_1 \cap \partial \Omega_2.$$
The function $a(\cu{x})$ is a high-contrast piecewise constant function that jumps across the interface
$$a(\cu{x}) = \left\{
\begin{array}{cc}
1 & \text{ in } \Omega_1, \\
10^{-m} & \text{ in } \Omega_2.
\end{array}
\right.$$
Figure \ref{fig:jumpa}\ref{fig:jumpu} presents an example of the coefficients and the corresponding solution, where the red region represents $\Omega_1$ and the blue region represents $\Omega_2$.

We employ the cell-centered finite volume discretization method \cite{wesseling1995introduction}. To achieve this, a uniform grid with spacing $ h = 1/(n+1) $ is used, ensuring that the interface (dashed line) is fitted with the dual element (gray region), as illustrated in Figure \ref{fig:fvm_grid}. The resulting discrete system is
$$ \begin{bmatrix}&s_n\\s_w&s_c&s_e\\&s_s&\end{bmatrix}_h u_h = h^2 f_h, $$
where the coefficients are 
$$
\begin{aligned}
&s_{w} = -\frac{2a(x-h, y)a(x, y)}{a(x-h, y) + a(x, y)}, 
&s_{e} = -\frac{2a(x+h, y)a(x, y)}{a(x+h, y) + a(x, y)}, \\
&s_{n} = -\frac{2a(x, y+h)a(x, y)}{a(x, y+h) + a(x, y)}, 
&s_{s} = -\frac{2a(x, y-h)a(x, y)}{a(x, y-h) + a(x, y)}, \\
&s_{c} = -(s_{n} + s_{s} + s_{e} + s_{w}).
\end{aligned}
$$
In comparison with the multi-scale property of anisotropic diffusion  equations, the multi-scale property of this equation is also influenced by the node positions. 

\begin{figure}[!htb]
    \centering
    \subfigure[$a(\cu{x})$]{\label{fig:testa}
    \includegraphics[width=0.35\textwidth]{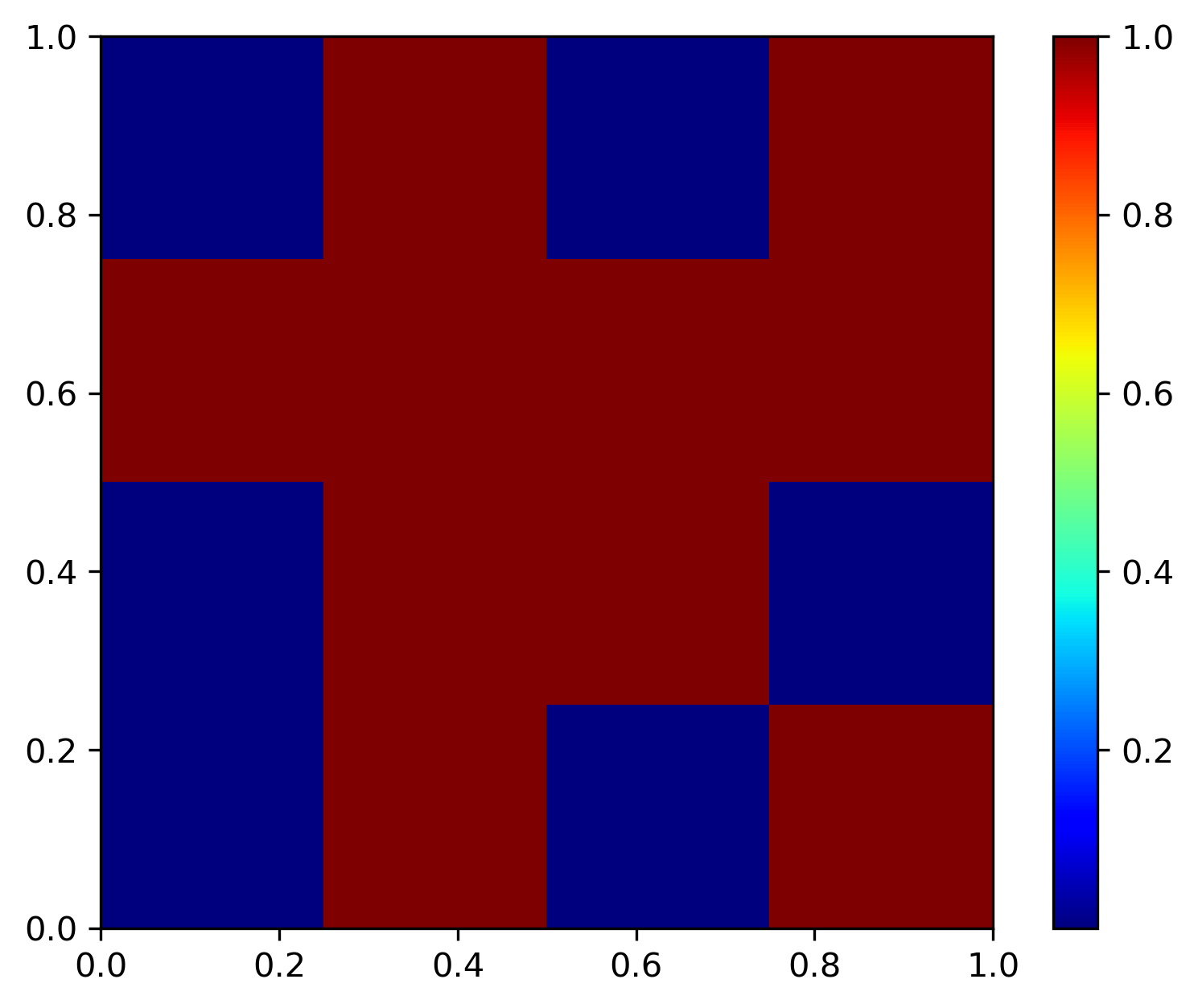}}\quad 
    \subfigure[Jacobi symbol]{
    \includegraphics[width=0.37\textwidth]{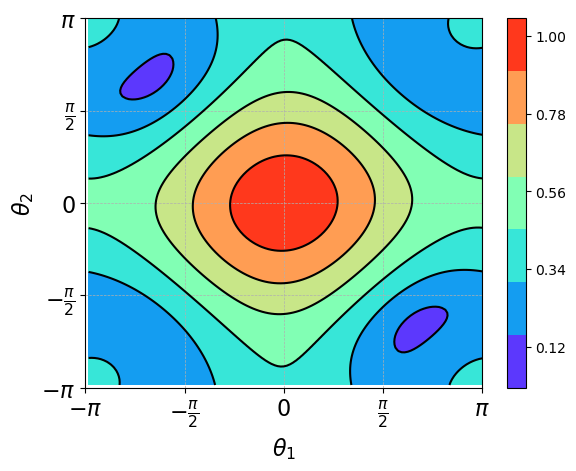}} 
    \caption{An example of jumping coefficient and corresponding Jacobi symbol.} 
    \label{fig:lfa_jump} 
\end{figure}

Consider $\cu{B}$ as the damped Jacobi method. We calculate its compressibility for error components with different frequencies using non-standard LFA \cite{kumar2019local}. Figure \ref{fig:lfa_jump} displays an example of $a(\cu{x})$ and the Jacobi symbol with $\omega = 2/3$ when solving corresponding discrete system. It can be seen that this $\cu{B}$ is effective for eliminating high-frequency error components, while the compressibility for low-frequency errors is poor. Let $\Theta^{\mathcal{H}} = [-\pi/2, \pi/2)^2$ and $\Theta^{\cu{B}} = [-\pi, \pi)^2 \backslash \Theta^{\mathcal{H}}$, then Assumption \ref{ass:B} is satisfied.

Next, we train FNS to ensure that $\mathcal{H}$ satisfies assumption \ref{ass:H}. Corresponding to PPDE \eqref{ppde}, $\boldsymbol{\mu} = a(\cu{x})$ in this example. We generate 10,000 different functions $a$ in the following way: split $\Omega$ into $4 \times 4$ checkerboard blocks, where the value of $a$ in each block is a constant, either $1$ or $10^{-m}$, with $m \in [4, 8]$. 
For each $a_i$, a random RHS $\cu{f_i} \sim \mathcal{N}(0, \cu{I})$ is sampled, and unsupervised training is then performed using the loss function \eqref{eq:loss_func}. The hyperparameters used during the training phase are the same as those used in the previous experiments.

\begin{figure}[!htbp]
    \centering
    \includegraphics[width=0.8\textwidth]{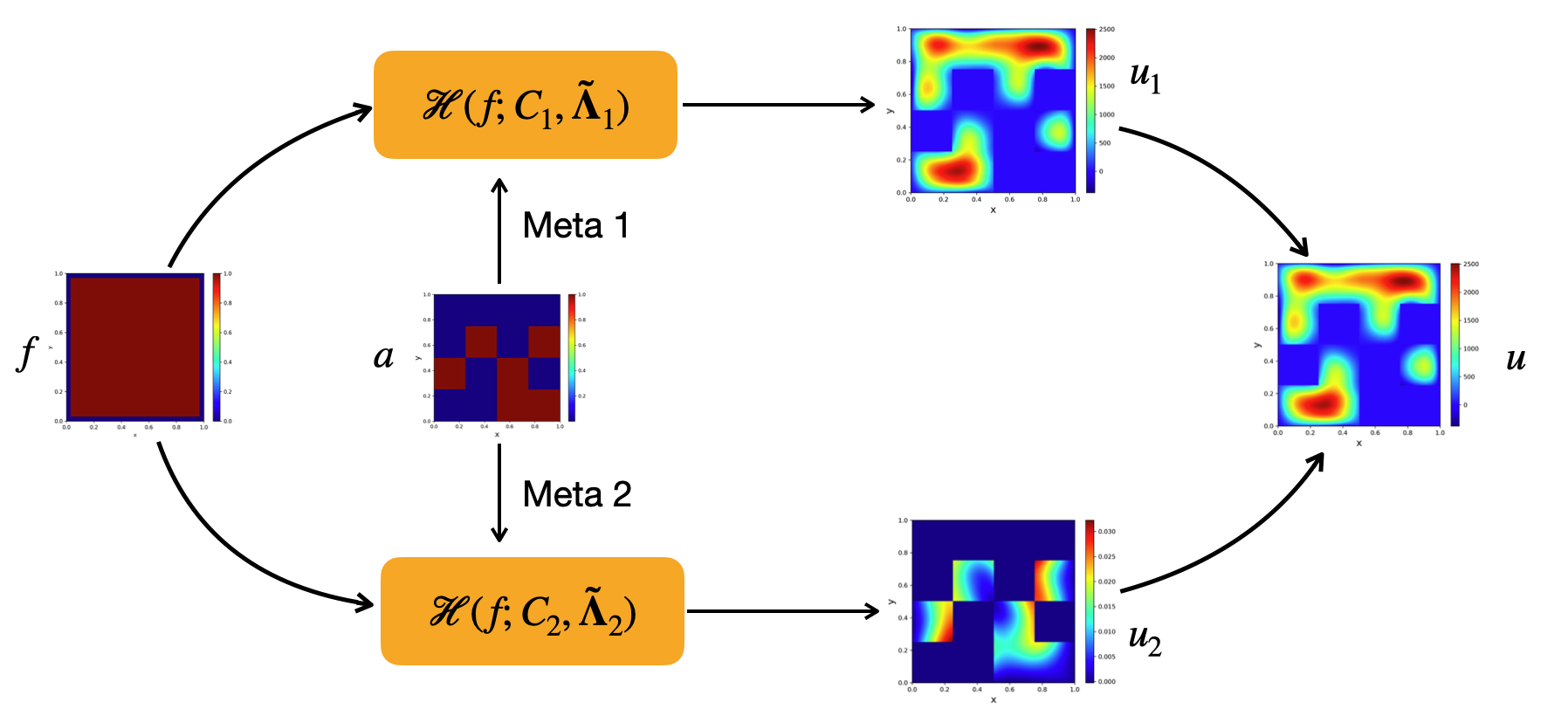}
    \caption{Calculation flow of $\mathcal{H}$ when solving the jumping diffusion equation.}
    \label{fig:Jump_flow}
\end{figure}

It is worth noting that due to the multi-scale property with respect to position, it is challenging to use one $\mathcal{H}$ to learn the correction values for all components simultaneously. Therefore, we classify all nodes into two parts on a geometric level: nodes in $\Omega_1$ form one category, and nodes in $\Omega_2$ form another category. Correspondingly, two separate $\mathcal{H}$ networks are used to learn the correction values for the first and second groups of nodes, respectively.
Figure \ref{fig:Jump_flow} shows the calculation flow of $\mathcal{H}$, where Meta 1 is used to learn corrections for nodes in $\Omega_1$ (blue region), and Meta 2 is used to learn corrections for nodes in $\Omega_2$ (red region). It can be seen that this geometric split naturally corresponds to the algebraic multi-resolution decomposition. This allows the network to learn the correction for nodes with the same order of magnitude only, making the learning process easier.

Next, we test the performance of FNS on varying scales. Table \ref{tab:jump} shows the required iteration counts to solve the discrete system when the relative residual is less than $10^{-6}$ for the coefficient $a$ shown in Figure \ref{fig:testa} and $m=8, f=1$. It can be seen that if FNS is only trained on the scale of $n=63$, it has poor generalization to other scales. However, if FNS is trained on the testing scale, it can still learn suitable parameters on different scales, which significantly reducing the number of iterations. In the future, we aim to design better network architectures and use more powerful optimization algorithms to facilitate training.
\begin{table}[!htb]
    \centering
    \caption{FNS iteration counts to satisfy $\|\cu{r}^{(k)}\|/\|\cu{f}\|<10^{-6}$, starting with zero initial value, for  jumping diffusion equations on different scales.}
    \label{tab:jump}
    \begin{tabular}{cccccc}
    \toprule
    $n$ & 15 & 31 & 63 & 127 & 255 \\ \midrule
    Only trained on $n=63$ & 60 & 44 & 23 & 106 & 435 \\ \midrule
    Trained on the testing scale & 12 & 23 & 23 & 41 & 70 \\ \bottomrule
    \end{tabular}
\end{table}

\begin{figure}[!htb]
    \centering
    \subfigure[$2\times 2$]{
    \includegraphics[width=0.18\textwidth]{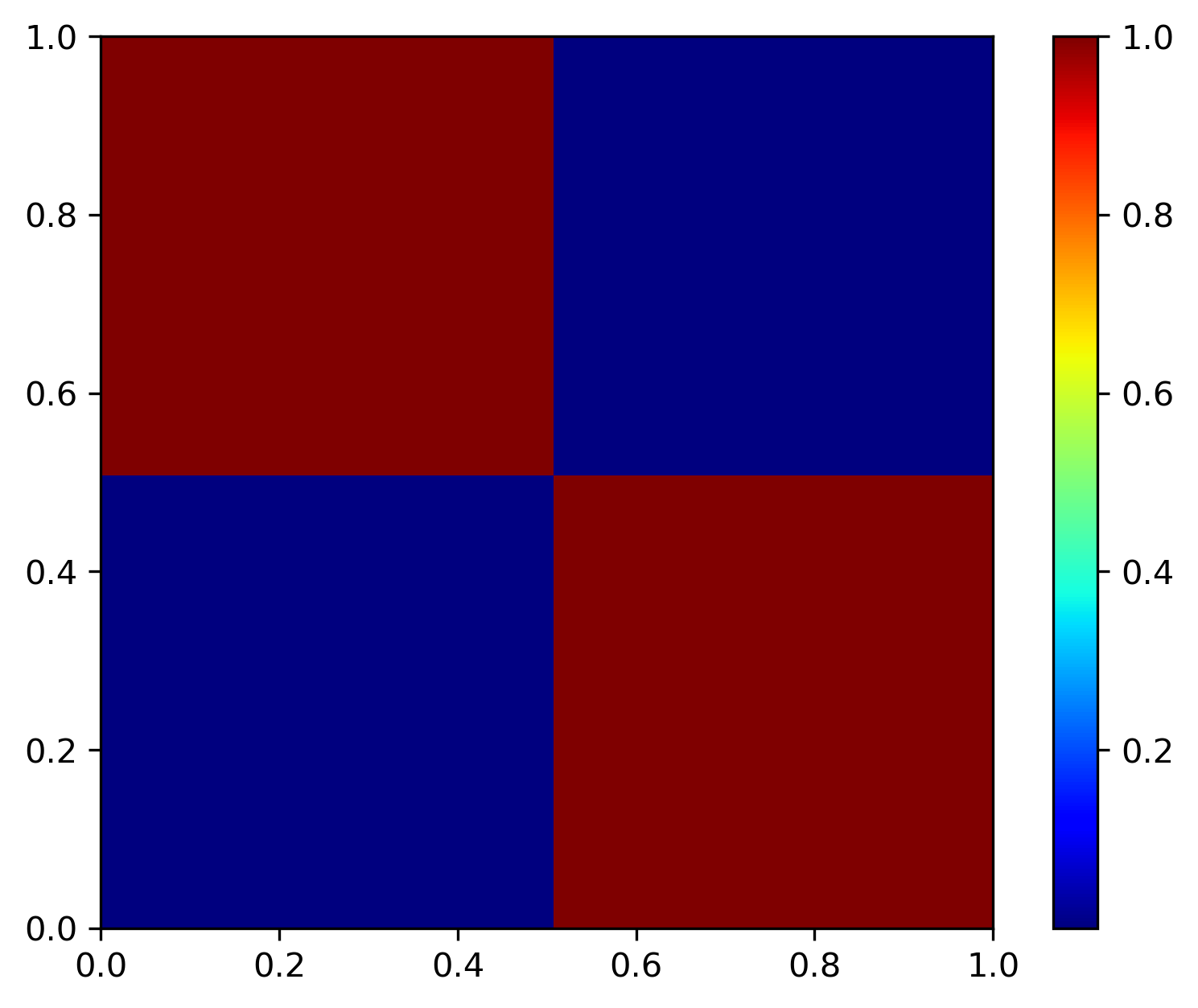}}
    \subfigure[$4\times 4$]{
    \includegraphics[width=0.18\textwidth]{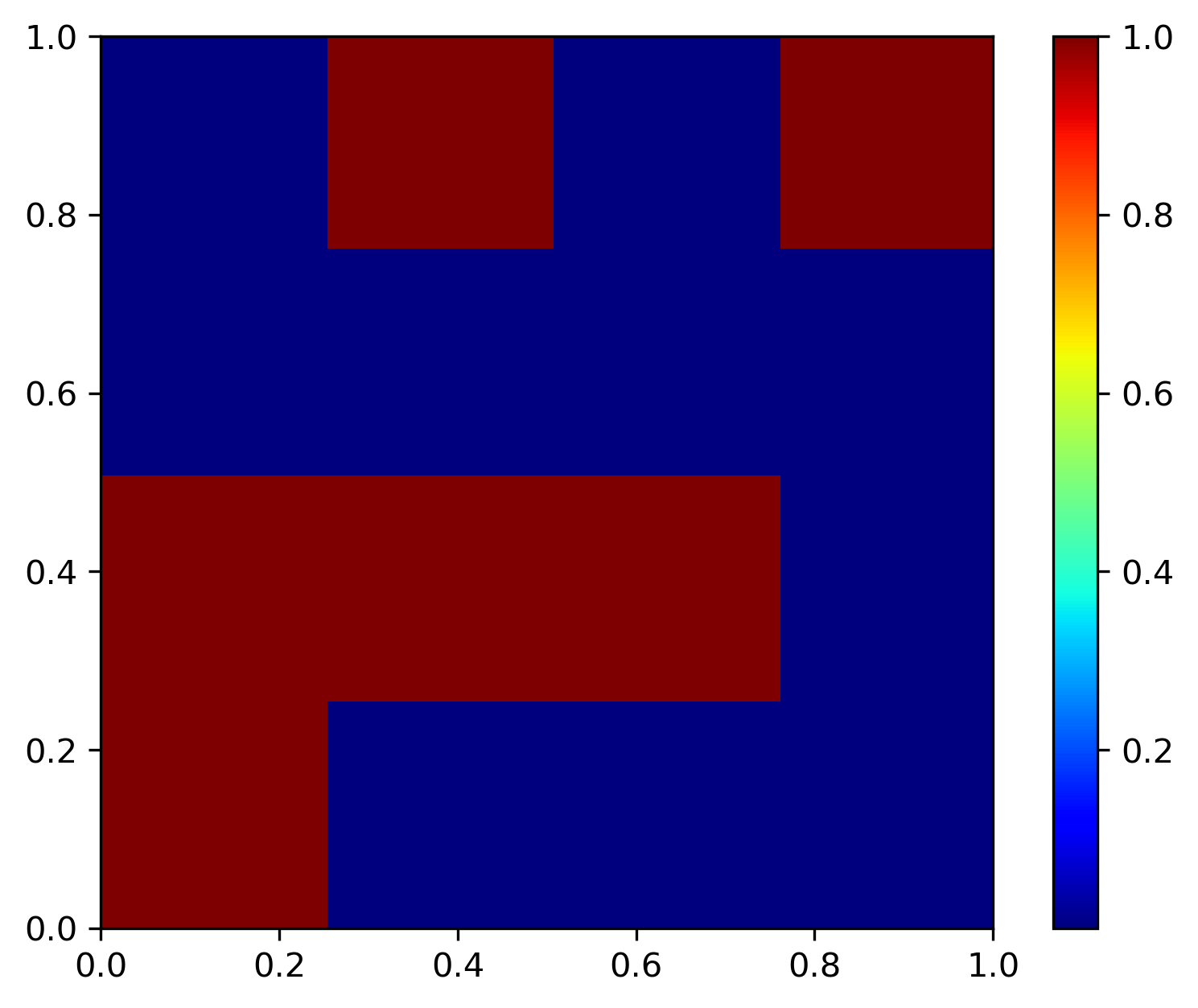}}
    \subfigure[$6\times 6$]{
    \includegraphics[width=0.18\textwidth]{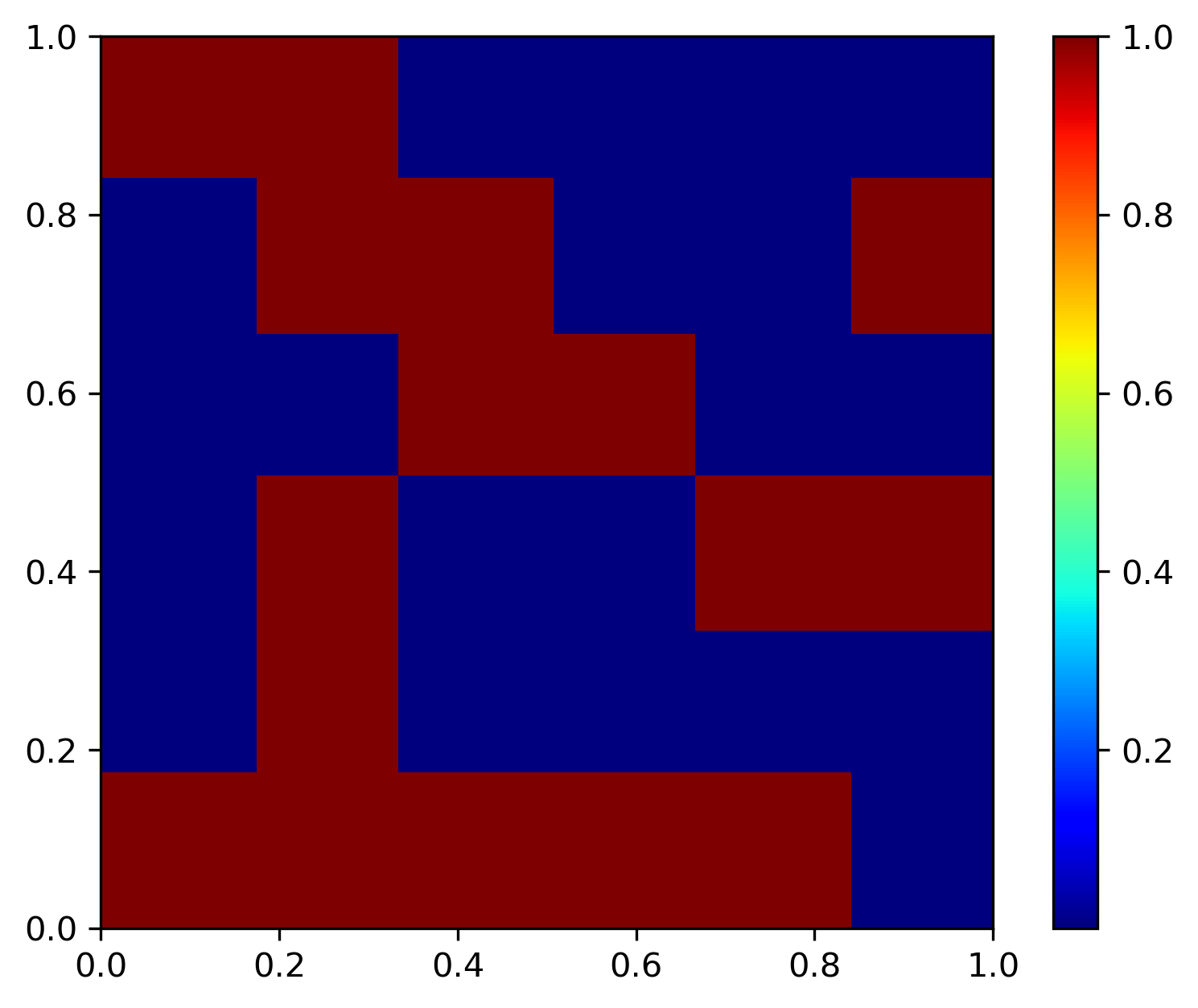}}   
    \subfigure[$40\times 40$]{
    \includegraphics[width=0.18\textwidth]{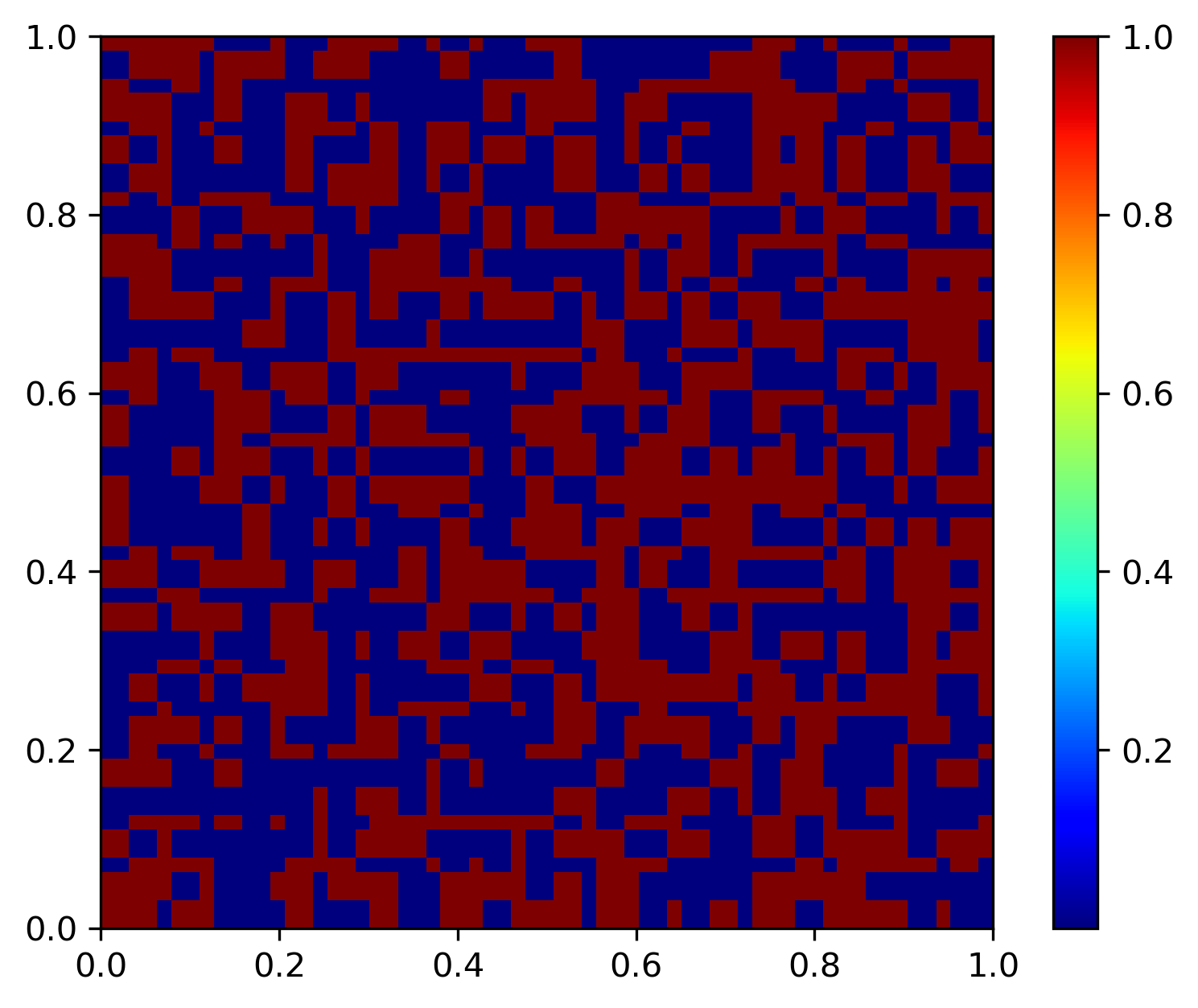}}   
    \subfigure[$n\times n$]{
    \includegraphics[width=0.18\textwidth]{coefN.png}}   \\
    \subfigure[$2\times 2$]{
    \includegraphics[width=0.18\textwidth]{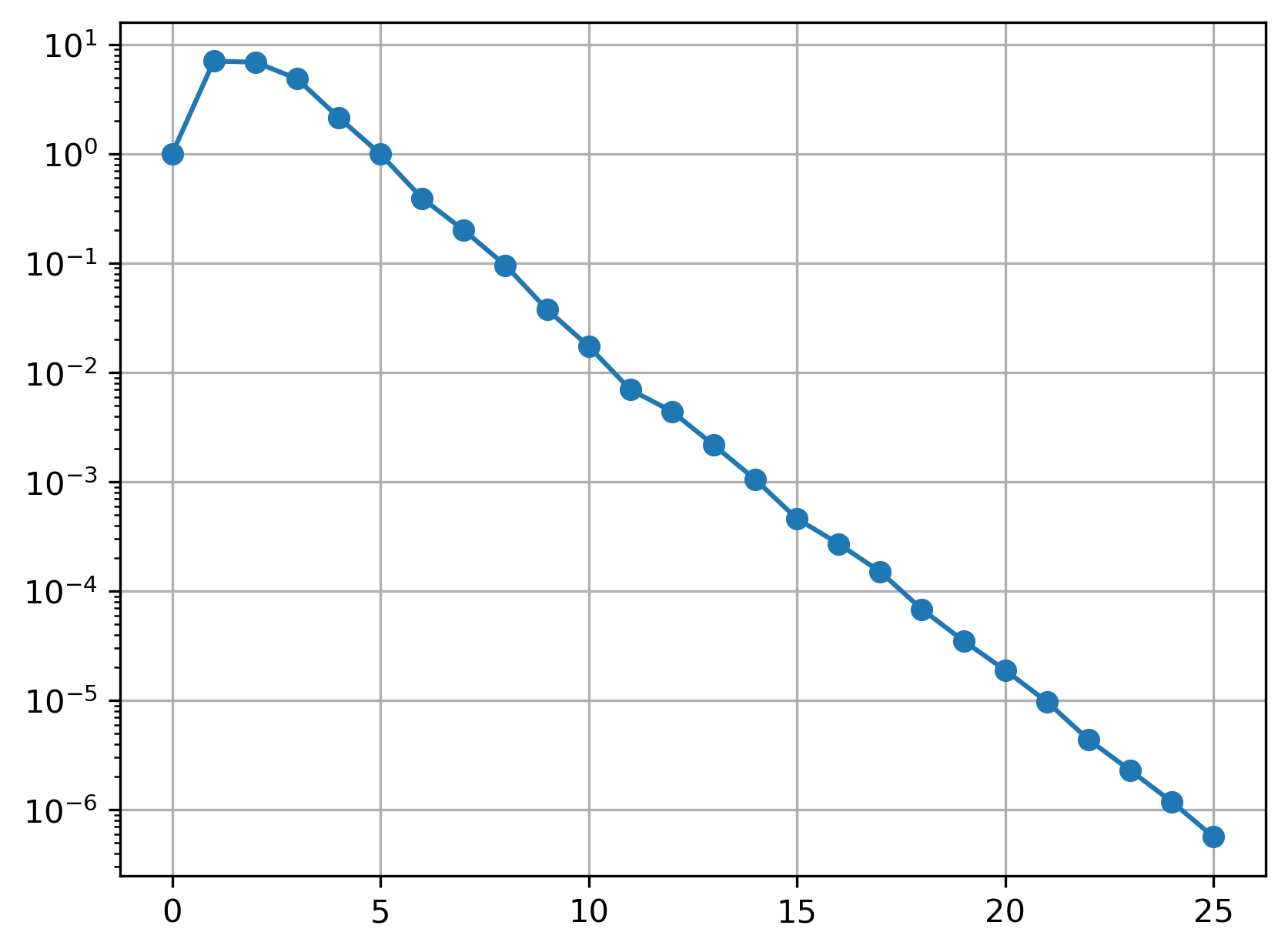}}
    \subfigure[$4\times 4$]{
    \includegraphics[width=0.18\textwidth]{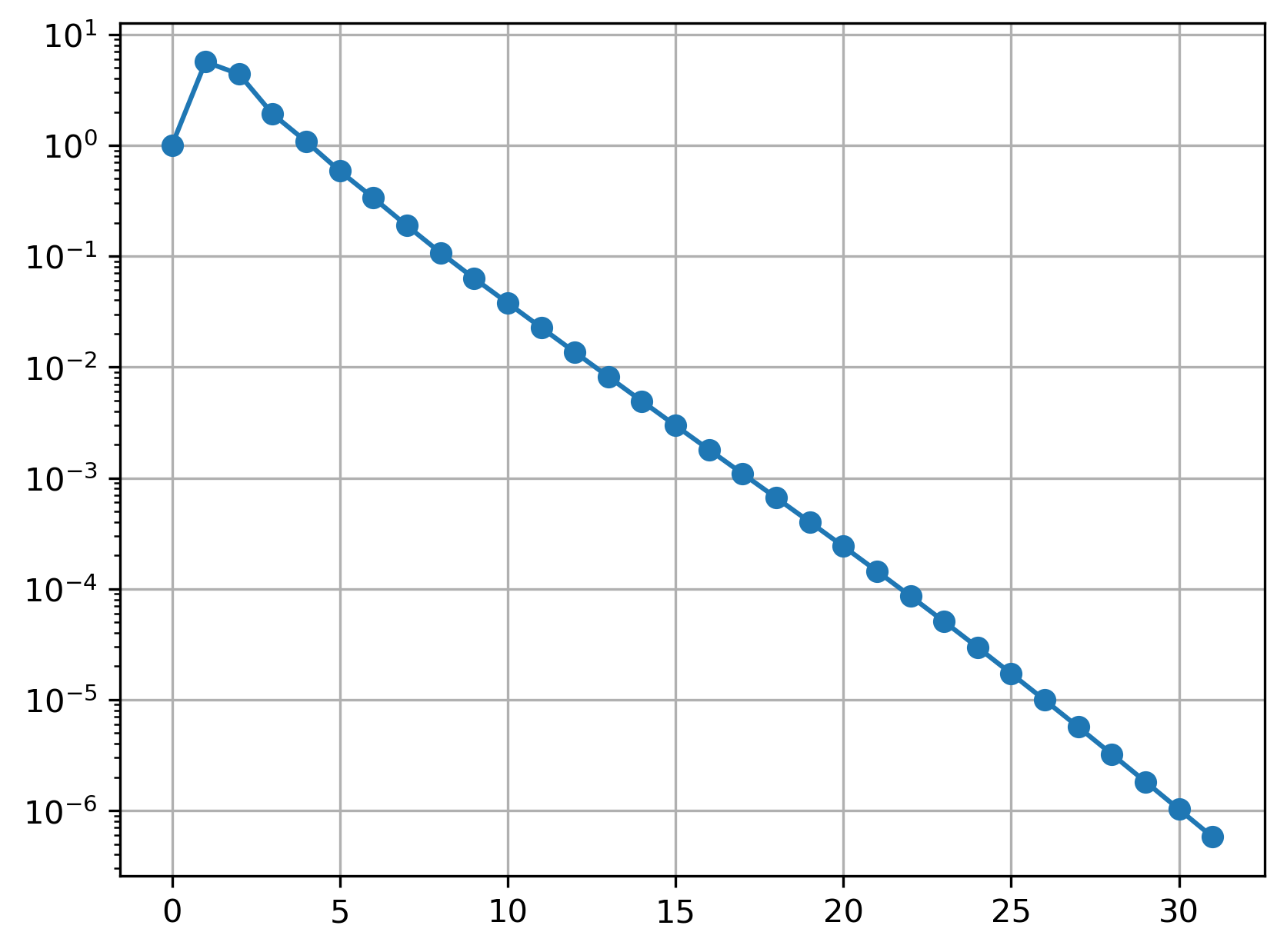}}
    \subfigure[$6\times 6$]{
    \includegraphics[width=0.18\textwidth]{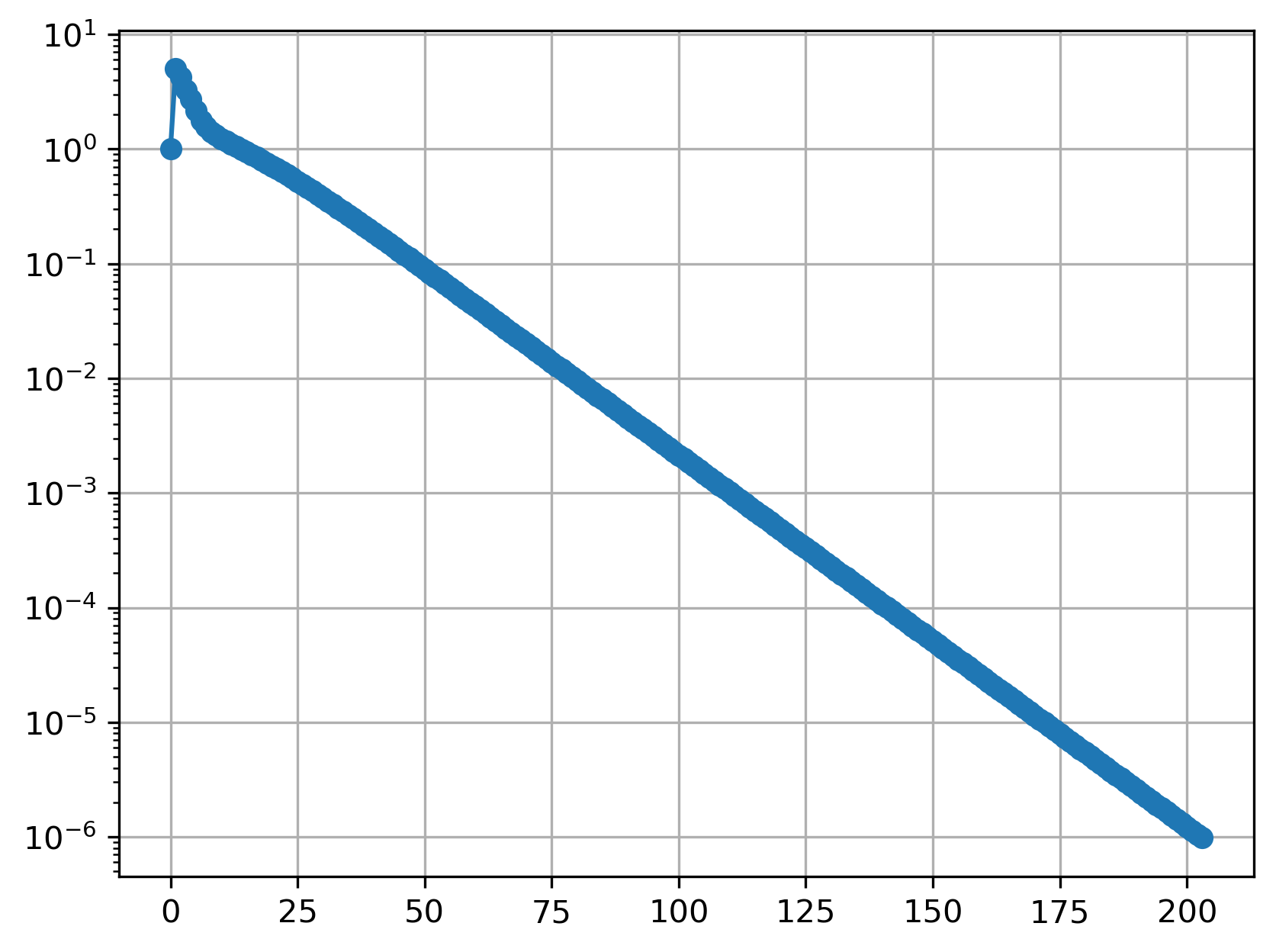}}   
    \subfigure[$40\times 40$]{
    \includegraphics[width=0.18\textwidth]{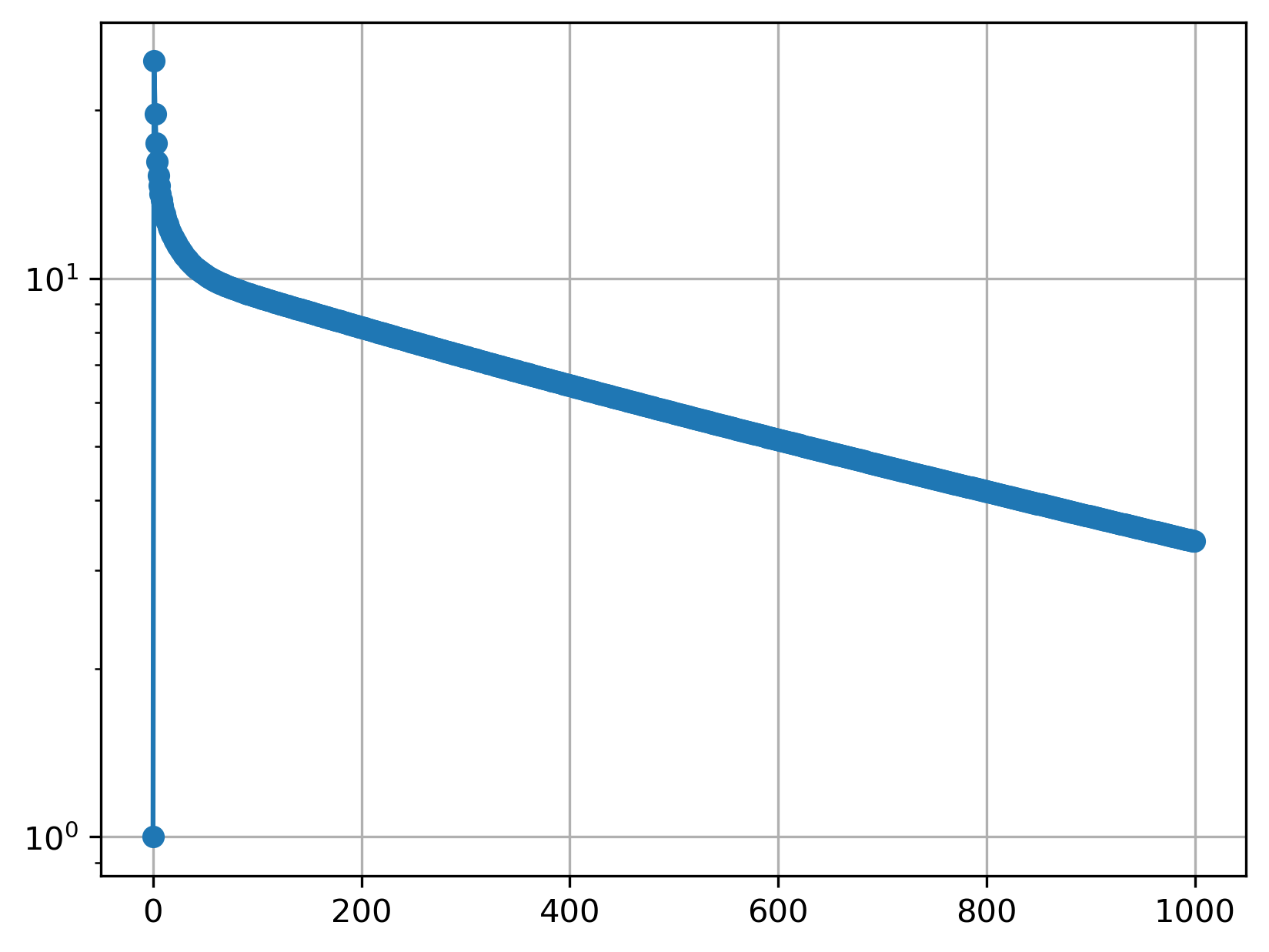}}   
    \subfigure[$n\times n$]{
    \includegraphics[width=0.18\textwidth]{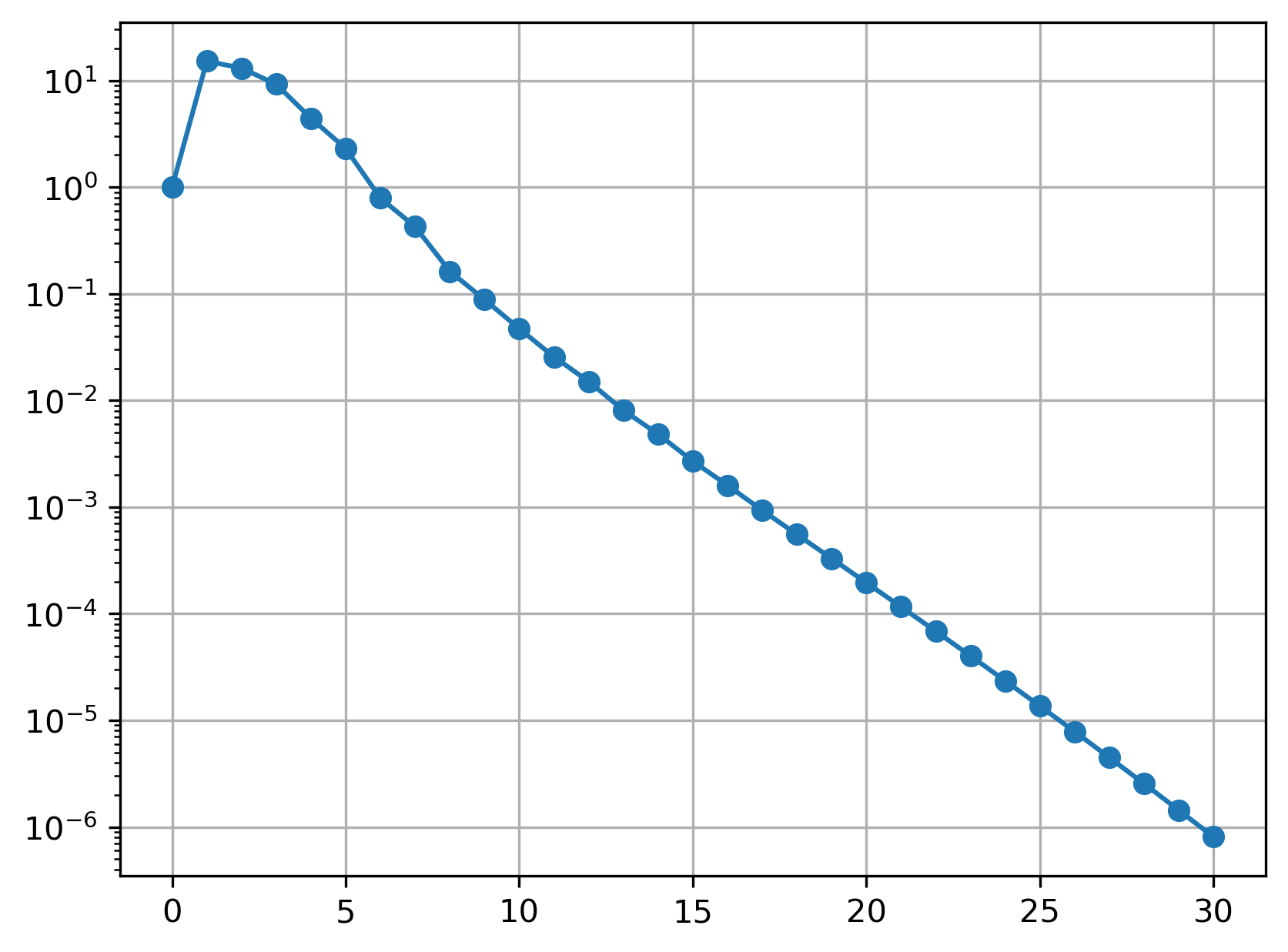}}  
    \caption{Top: Diffusion coefficients obtained randomly by using different checkerboard block sizes. Bottom: Convergence history of FNS for solving the corresponding discrete system.}
    \label{fig:diff_jump}
\end{figure}

It should also be pointed out that the distribution of the coefficient $a$ determines the multi-scale property of the discrete system. The more random the distribution, the stronger the multi-scale property, making the system more challenging to solve.
Figure \ref{fig:diff_jump} shows the behavior of FNS applied to discrete systems with different coefficients and $m=8$, $f=1$. The first four columns indicate that as the number of checkerboard blocks increases, the interface becomes more complex, resulting in reduced connectivity within $\Omega_1$ or $\Omega_2$. Consequently, the efficiency of FNS decreases. However, as shown in the last column, if $\Omega_1$ and $\Omega_2$ are connected, even if the interface is highly discontinuous, the efficiency of FNS remains comparable to that of the $4 \times 4$ block configuration.

\subsection{Helmholtz equation}
Consider the 2D Helmholtz equation with impedance boundary condition
\begin{equation} 
    \left\{
    \begin{aligned}
        -\Delta u-k^{2}u &=f, \text{ in } \Omega=(0,1)^2, \\
        \nabla u \cdot \mathbf{n} -iku&=0, \text{ on } \partial \Omega, \\
    \end{aligned}
    \right.
    \label{eq:helm}
\end{equation}
where $k(\mathbf{x})$ is the wave number, and $f(\mathbf{x}) = \delta(\mathbf{x} - \mathbf{x}_0)$ represents a point source in $\Omega$.
Using the second-order central finite difference method on a uniform grid, the discretized Helmholtz operator is given by
\begin{equation}
    \frac{1}{h^2}\left[\begin{array}{ccc} 
     0 & -1 & 0 \\
     -1 & 4-k_{i,j}^2 h^2 & -1 \\
     0 & -1 & 0  
     \end{array}\right],
\end{equation}
where $h = 1/(n+1)$ in both the $x$- and $y$-directions.
According to the Shannon sampling principle, at least 10 grid points per wavelength are required, resulting in a \textit{large-scale} linear system.

\begin{figure}[!htbp]
    \centering
    \subfigure[$\omega=2/3$]{\label{fig:helm1}
    \includegraphics[width=0.3\textwidth]{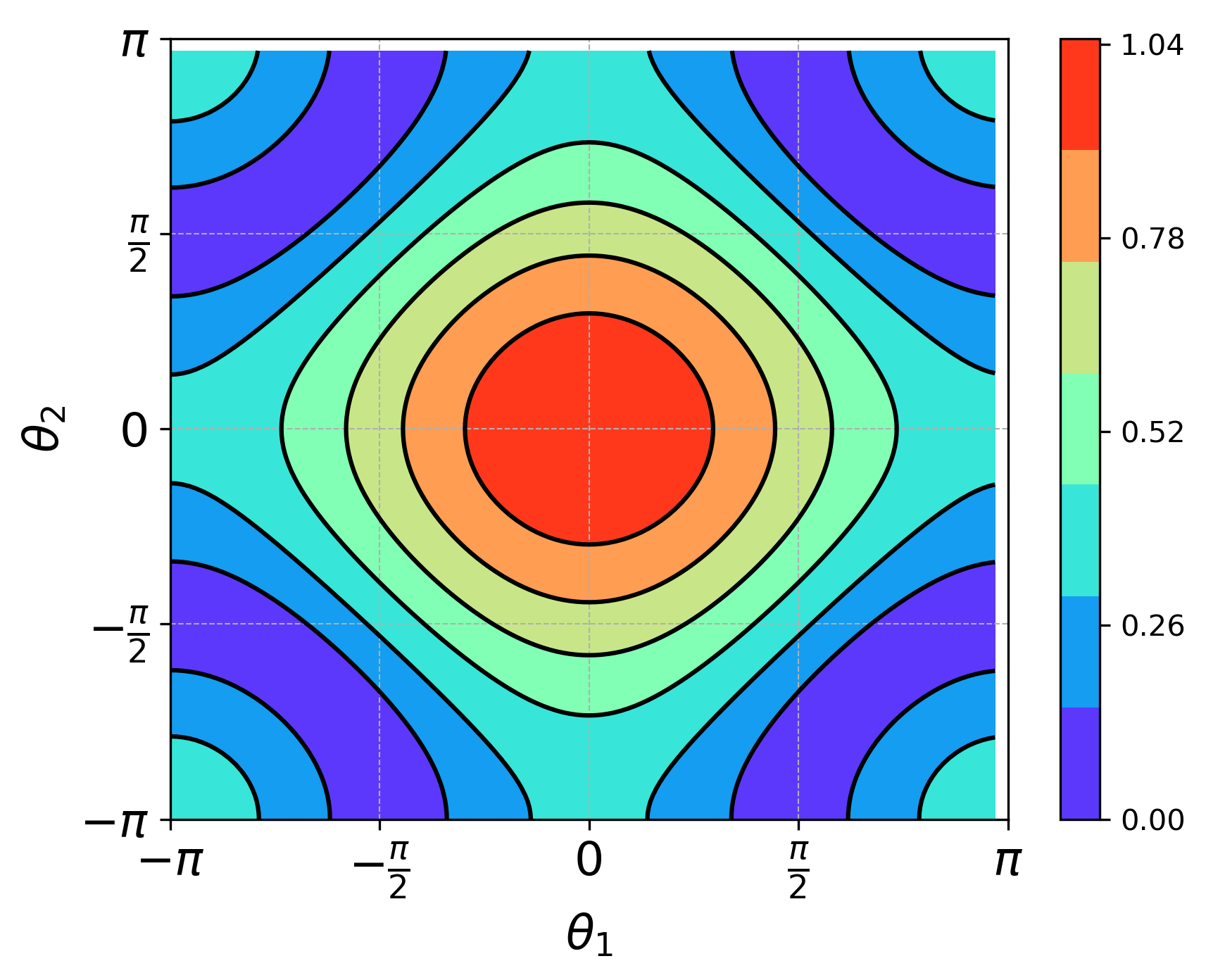}}\quad
    \subfigure[$\omega=3/4$]{\label{fig:helm2}
        \includegraphics[width=0.3\textwidth]{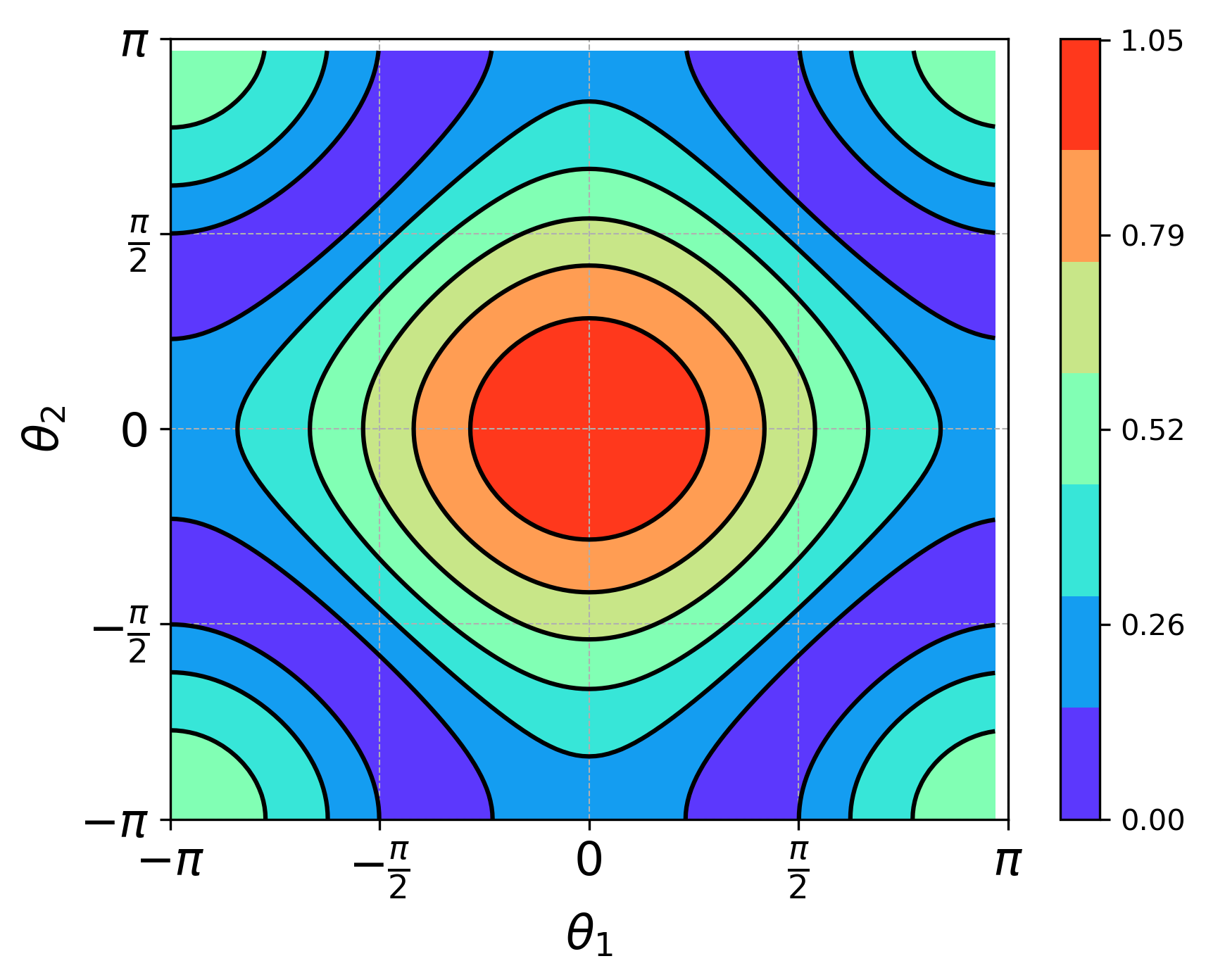}}
    \subfigure[$\omega=1$]{\label{fig:helm3}
        \includegraphics[width=0.3\textwidth]{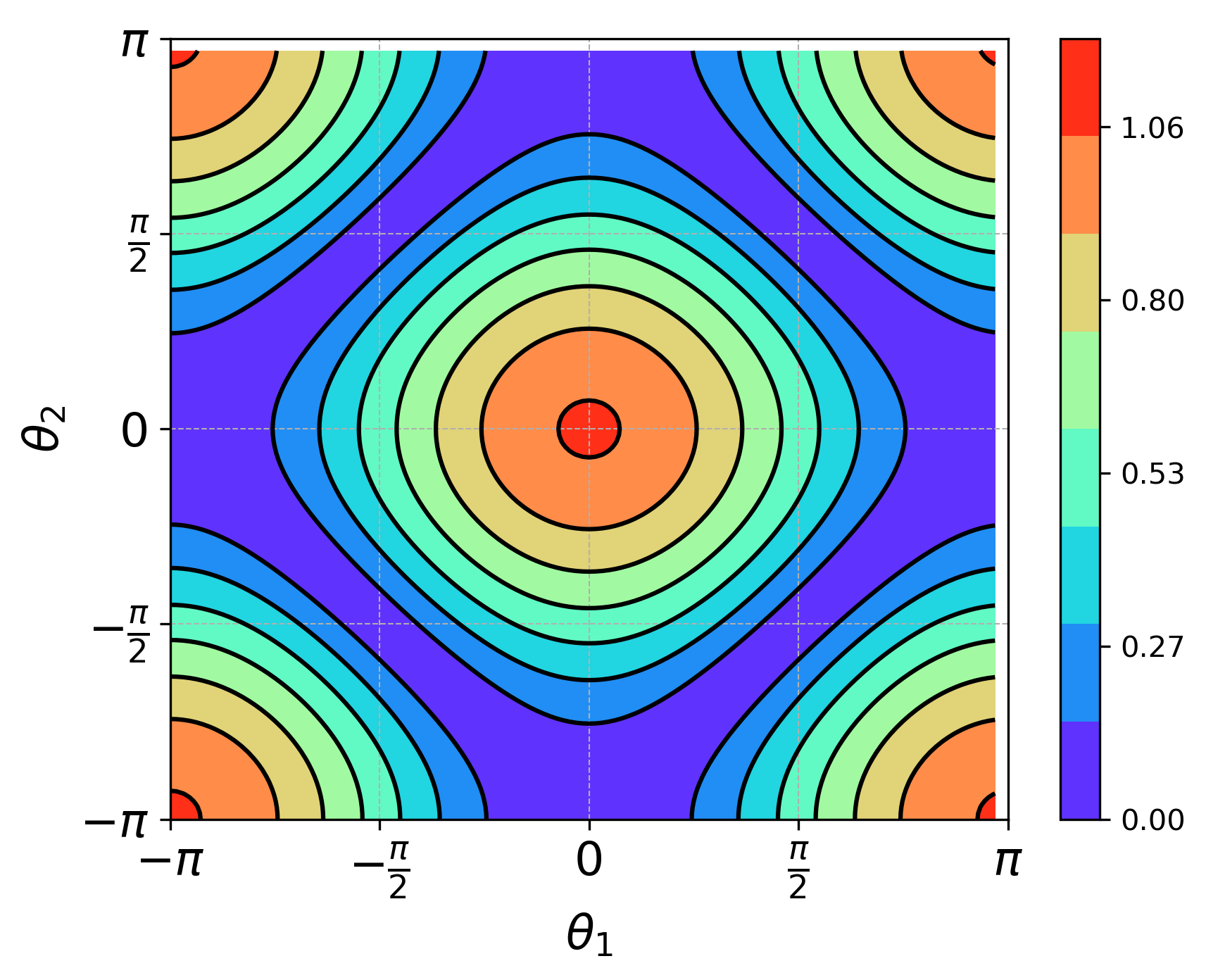}}
    \caption{Jacobi symbol applied to Helmholtz equation with  $k=20\pi$}
    \label{fig:helm_LFA}
\end{figure}

We first perform LFA on the damped Jacobi smoother using an example of a constant wave number. Figure \ref{fig:helm_LFA} shows the Jacobi symbol with different $\omega$ when $k(\mathbf{x}) = 20\pi$.
It can be observed that regardless of the choice of $\omega$, the damped Jacobi method always amplifies the lowest-frequency error, corresponding to the case of $\varepsilon_{\cu{B}} > 0$ in Assumption \ref{ass:B}. However, the damped Jacobi method still exhibits good smoothing effect for high-frequency errors. Let $\Theta^{\mathcal{H}} = [-\pi/2, \pi/2)^2$ and $\Theta^{\cu{B}} = [-\pi, \pi)^2 \setminus \Theta^{\mathcal{H}}$. When $\omega = 2/3$, the smoothing effect of the damped Jacobi method is optimal.
Therefore, we choose $\cu{B}$ as the damped Jacobi method with $\omega = 2/3, M=1$, ensuring that Assumption \ref{ass:B} is satisfied.

\begin{figure}[!htbp]
\centering
\subfigure[$k(\cu{x})$]{\label{fig:helm_k}
\includegraphics[width=0.3\textwidth]{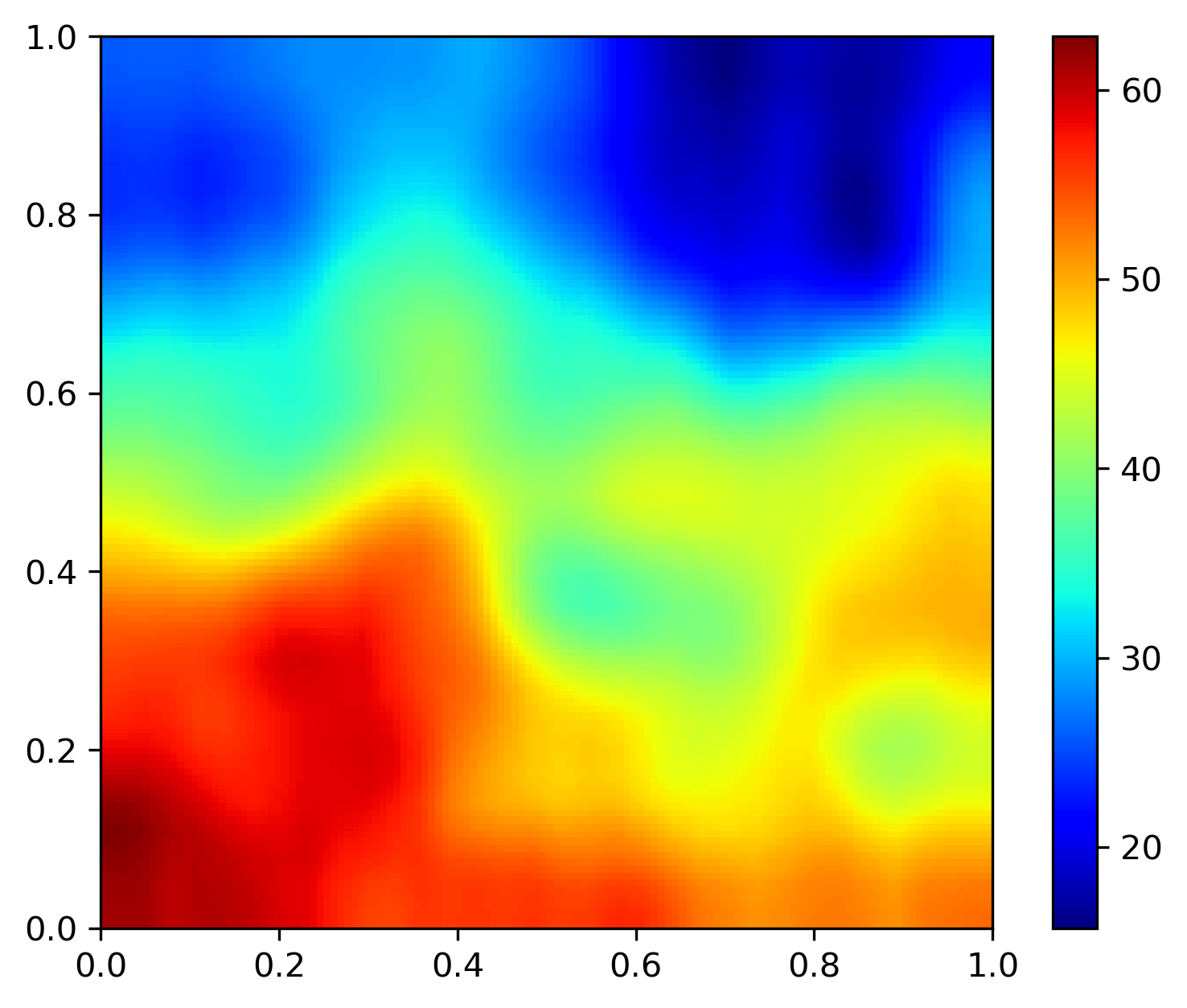}}
\subfigure[$u(\cu{x})$]{
\includegraphics[width=0.3\textwidth]{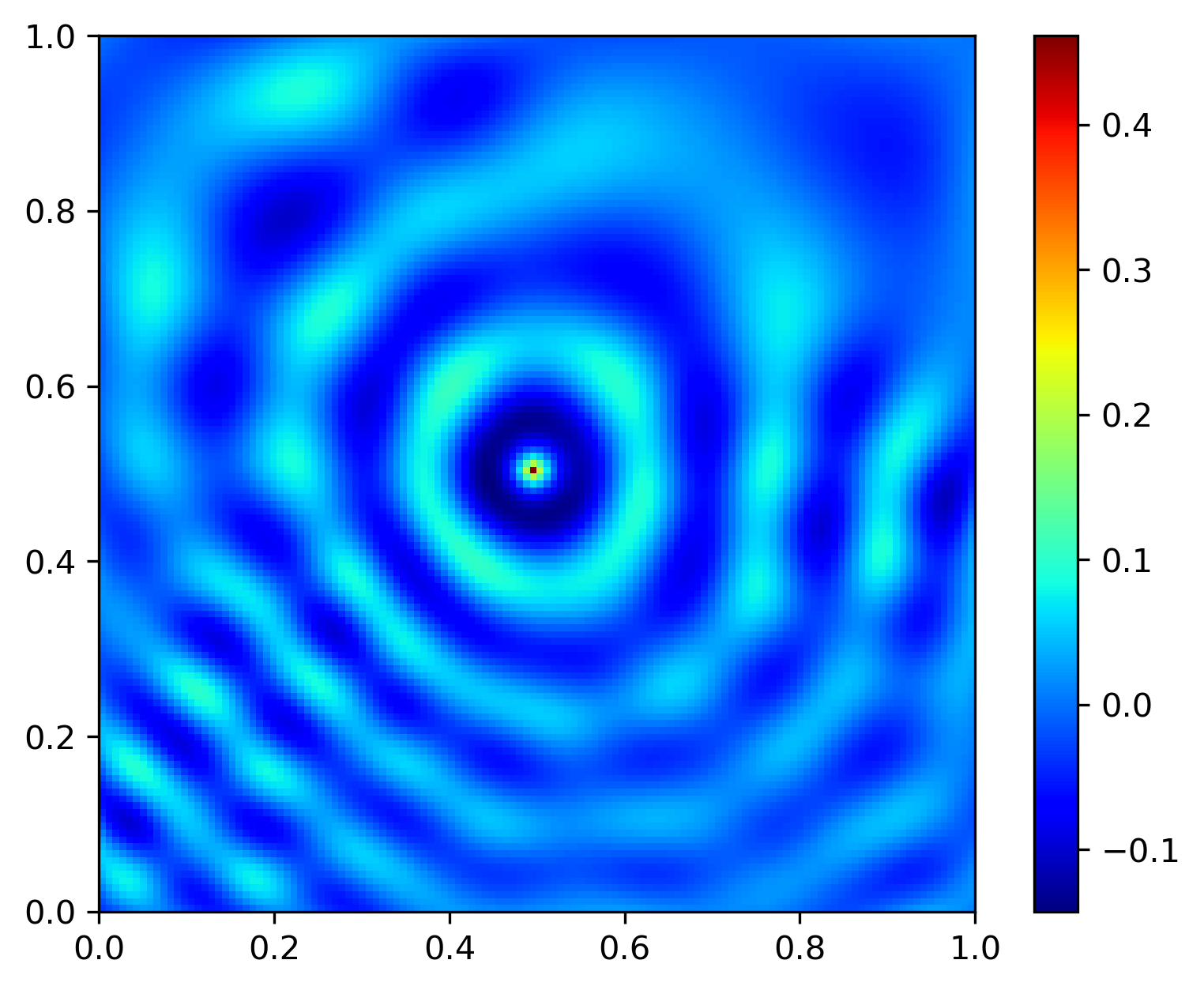}}
\subfigure[Convergence history]{\label{fig:helm_res}
\includegraphics[width=0.31\textwidth]{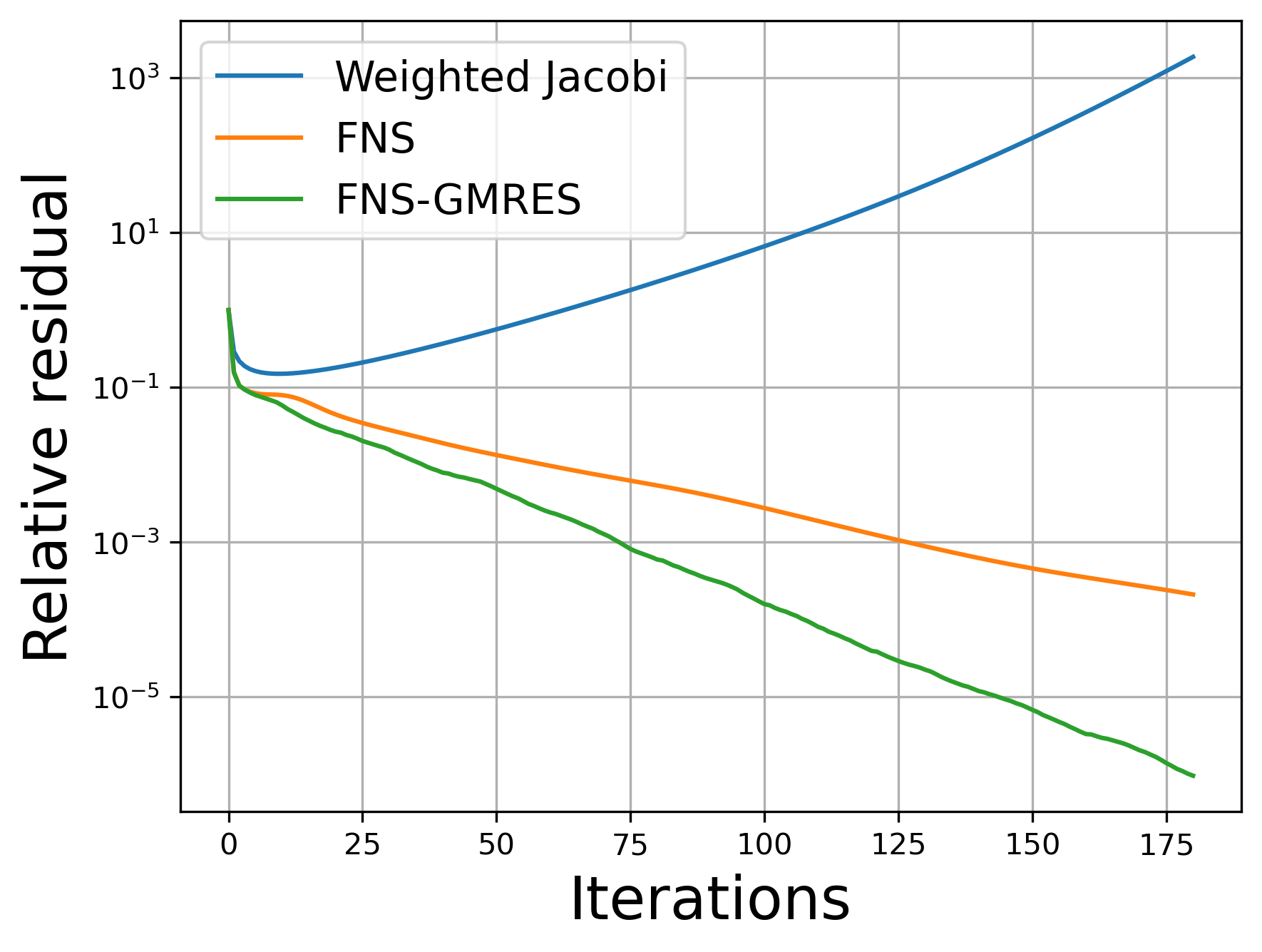}}
\caption{Helmholtz test example. (a) Wave number. (b) Solution. (c) Convergence history of FNS and damped Jacobi method for solving the corresponding discrete system.}
\end{figure}

Next, we train FNS to ensure that $\mathcal{H}$ satisfies assumption \ref{ass:H}. Corresponding to PPDE \eqref{ppde}, $\boldsymbol{\mu} = k(\cu{x})$ in this example. We generate 10,000 different functions $k$ converted from the CIFAR-10 dataset following \cite{azulay2022multigrid}.
Figure \ref{fig:helm_k} shows a test example $k(\cu{x})$.
The RHS are point source located at the center of the domain and unsupervised training is then performed using the loss function \eqref{eq:loss_func}. The hyperparameters used during the training phase are the same as those used in the previous experiments.

\begin{figure}[!htbp]
\centering
\includegraphics[width=\textwidth]{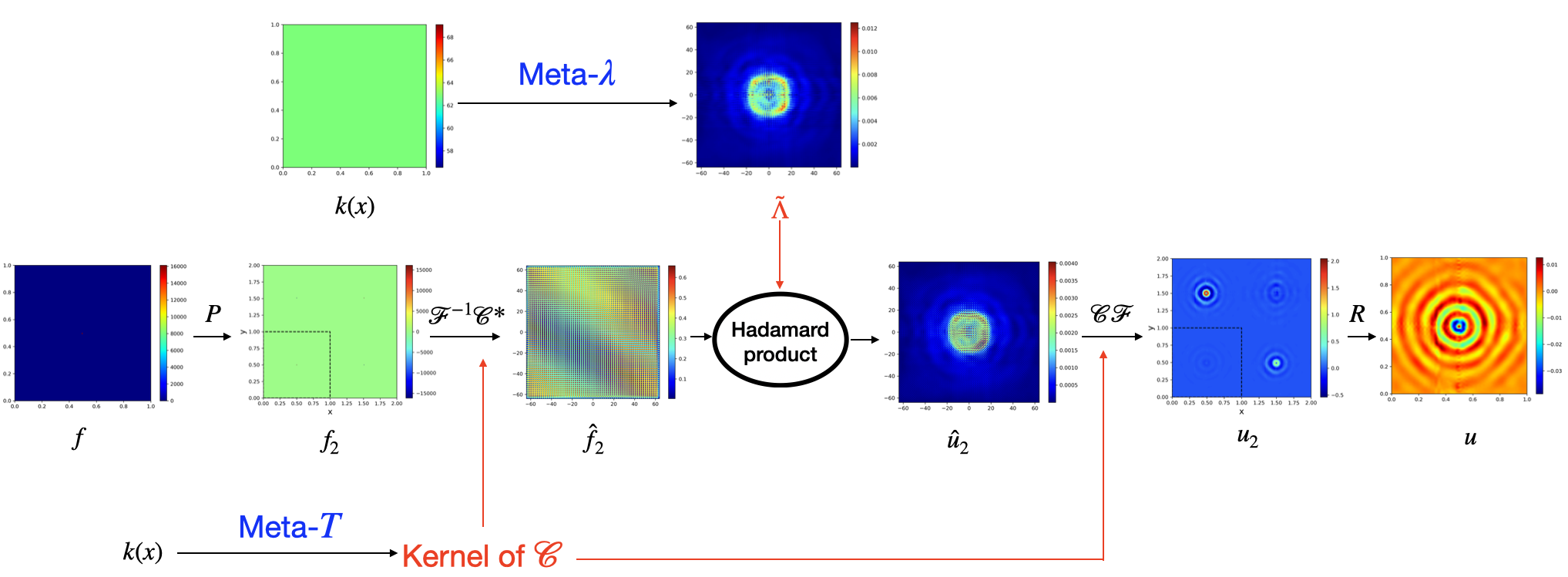}
\caption{Calculation flow of $\mathcal{H}$ when solving the Helmholtz equation.}
\label{fig:helm_flow}
\end{figure}

After training, we performed test experiments. Figure \ref{fig:helm_flow} illustrates the calculation flow of $\mathcal{H}$ when it receives a pair of ${\cu{\mu}_i, \cu{f}_i}$. Here, $k = 20\pi$, and $f$ is a point source located at the center.
From Figure \ref{fig:helm_flow}, we observe that the $\cu{\tilde{\Lambda}}$ provided by Meta$-\lambda$ is large in $\Theta^{\mathcal{H}}$ but small in $\Theta^{\cu{B}}$, which aligns with our expectations.

Next, we conducted convergence tests for a variable wave number (Figure \ref{fig:helm_k}). Figure \ref{fig:helm_res} presents the convergence results. It is evident that even in scenarios where the smoother $\cu{B}$ fails to converge, $\mathcal{H}$ effectively learns to correct the errors.
Additionally, FNS can function as a preconditioner for GMRES, leading to further acceleration.

\section{Conclusion and discussion}\label{sec:final}

\subsection{Conclusions}
In this paper, we first conduct a convergence analysis from a spectral viewpoint for general DL-HIM. Under reasonable assumptions on smoothing operator $\cu{B}$ and neural operator $\mathcal{H}$, the convergence rate of DL-HIM is shown to be independent of PDE parameters and grid size.
Guided by this framework, we design an FFT-based neural operator $\mathcal{H}$ to overcome spectral bias and meet the corresponding assumptions, resulted DL-HIM known as FNS. The convergence speed of FNS is guaranteed by the universal approximation and discretization  invariance of the meta subnets, which requires sufficient training and is easier said than done.
We verify the theoretical results and test the computational efficiency of FNS through numerical experiments on a variety of PDE discrete systems. For single-scale problems, $\mathcal{H}$ can easily learn error components with frequencies complementary to $\cu{B}$. For multi-scale problems, $\mathcal{H}$ needs additional help to learn a fast FNS.

\subsection{Scope, limitations, and future work}\label{sec:Limitations}

\textbf{On the scale-independence of FNS convergence rates.}
Numerical experiments on random and anisotropic diffusion equations indicate that FNS achieves convergence rates that are largely independent of discretization scale when tested on systems that are either smaller than or moderately larger than (e.g., up to twice the size of) the training scale. However, when the test scale significantly exceeds the training scale, the number of iterations required increases substantially. These findings suggest that FNS exhibits scale-independent convergence primarily within a local neighborhood of the training scale, but has not yet achieved universal scale invariance. Incorporating training data from larger-scale problems can improve generalization and enable convergence performance at extended scales comparable to that at smaller ones.
On the other hand, this limitation may also be influenced by the choice of Meta-network architecture. We experimented with both UNet and FNO as Meta-$\lambda$. While FNO generally outperformed UNet and is therefore adopted in this work, UNet achieved better results in the case of random diffusion equations. This suggests that the Meta-network architecture plays a non-negligible role in determining convergence behavior. Designing more effective Meta-network architectures to improve scalability and generalization of FNS remains a key direction for future research.

\textbf{On the consumption of training resources.}
The experiments in this paper are mainly conducted at the grid scale of $n=63$, which is a representative training scale based on a balance between computational efficiency and model performance.
Table~\ref{tab:train_resources} summarizes the training resource consumption observed on an NVIDIA A100 GPU for various problem sizes, providing justification for this choice:
\begin{enumerate}
\item \textbf{Hardware utilization}: At $n = 63$,
\begin{itemize}
    \item Computational resources are effectively utilized. Increasing $n$ from 31 to 63 leads to a VRAM usage increase of less than 3x and a per-epoch training time increase of only 2x, indicating underutilization at smaller scales.
    \item Further increasing $n$ beyond 63 results in a more than 3x increase in training time, suggesting that $n = 63$ achieves near-optimal parallelization efficiency on our hardware.
\end{itemize}
\item \textbf{Acceptable training cost}: Training at $n = 63$ remains computationally feasible while producing performant models.
\end{enumerate}
\begin{table}[!htb]
    \centering
    \caption{Training resource consumption on an NVIDIA A100 GPU for different problem sizes of random diffusion equations. Each configuration uses $K=3$ in the loss function.}
    \label{tab:train_resources}
    \begin{tabular}{cccccc}
    \toprule
    Grid size $n$ & Data size & Batch size & VRAM usage (MB) & Time per epoch (s) \\
    \midrule
    31  & $10^5$ & 40 & 5837  & 21  \\
    63  & $10^5$ & 40 & 14179 & 40  \\
    127 & $10^5$ & 40 & 44083 & 125 \\
    \bottomrule
    \end{tabular}
\end{table}
Since the current network architecture indeed demonstrates strong generalization capability primarily near the training scale. To address large-scale computational demands, scaling up the training regime—i.e., training on larger problem sizes—becomes imperative.
This is a key direction in our ongoing work, aimed at improving the robustness and scalability of FNS across a broader range of discretization levels.

\textbf{On the applicability of FNS to nonsymmetric PDEs.} The design of the operator $\mathcal{H}$ is motivated by the eigendecomposition of $\cu{A}$, which fundamentally assumes that $\cu{A}$ is symmetric or normal. Nevertheless, our numerical experiments in Section \ref{sec:convection} on convection-diffusion equations demonstrate that FNS retains weak dependence on the discretization scale and achieves near-linear runtime, even when applied to nonsymmetric SUPG discretizations. Extending this capability to more general nonsymmetric discrete systems remains an important direction for future research.

\textbf{On the applicability of FNS to general discrete systems.}
The discrete systems considered in this work are all derived from structured grids on regular domains. The applicability of FNS to more general cases—such as discrete systems arising from irregular domains, unstructured grids, three-dimensional PDEs, or PDE systems—remains an open question. A key challenge is the current reliance of FNS on the FFT, which assumes a regular and consistent grid structure. Extending the framework to accommodate arbitrary discretizations will be an important direction for future work, potentially through the integration of graph neural networks or geometry-aware architectures.

\section*{Acknowledgments}
We would like to thank the referees for their insightful comments that greatly improved the paper. 
This work is funded by the NSFC grants (12371373, 12171412). SS is supported by Science Challenge Project (TZ2024009).
KJ is partly supported by the National Key R\&D Program of China (2023YFA1008802), the Science and Technology Innovation Program of Hunan Province (2024RC1052).
This work was carried out in part using computing resources at the High Performance Computing Platform of Xiangtan University.

\bibliographystyle{elsarticle-num} 
\bibliography{reference}

\begin{thebibliography}{10}
\expandafter\ifx\csname url\endcsname\relax
  \def\url#1{\texttt{#1}}\fi
\expandafter\ifx\csname urlprefix\endcsname\relax\def\urlprefix{URL }\fi
\expandafter\ifx\csname href\endcsname\relax
  \def\href#1#2{#2} \def\path#1{#1}\fi

\bibitem{saad2003iterative}
Y.~Saad, Iterative methods for sparse linear systems, SIAM, 2003.

\bibitem{golub1961chebyshev}
G.~H. Golub, R.~S. Varga, Chebyshev semi-iterative methods, successive overrelaxation iterative methods, and second order richardson iterative methods, Numerische Mathematik 3~(1) (1961) 157--168.

\bibitem{hestenes1952methods}
M.~R. Hestenes, E.~Stiefel, et~al., Methods of conjugate gradients for solving linear systems, Vol.~49, NBS Washington, DC, 1952.

\bibitem{nesterov1983method}
Y.~Nesterov, A method of solving a convex programming problem with convergence rate o (1/k** 2), Doklady Akademii Nauk SSSR 269~(3) (1983) 543.

\bibitem{luo2022differential}
H.~Luo, L.~Chen, From differential equation solvers to accelerated first-order methods for convex optimization, Mathematical Programming 195~(1) (2022) 735--781.

\bibitem{beck2013convergence}
A.~Beck, L.~Tetruashvili, On the convergence of block coordinate descent type methods, SIAM journal on Optimization 23~(4) (2013) 2037--2060.

\bibitem{birkhoff1962alternating}
G.~Birkhoff, R.~S. Varga, D.~Young, Alternating direction implicit methods, in: Advances in computers, Vol.~3, Elsevier, 1962, pp. 189--273.

\bibitem{walker2011anderson}
H.~F. Walker, P.~Ni, Anderson acceleration for fixed-point iterations, SIAM Journal on Numerical Analysis 49~(4) (2011) 1715--1735.

\bibitem{saad1986gmres}
Y.~Saad, M.~H. Schultz, Gmres: A generalized minimal residual algorithm for solving nonsymmetric linear systems, SIAM Journal on scientific and statistical computing 7~(3) (1986) 856--869.

\bibitem{chow2015fine}
E.~Chow, A.~Patel, Fine-grained parallel incomplete lu factorization, SIAM journal on Scientific Computing 37~(2) (2015) C169--C193.

\bibitem{trottenberg2000multigrid}
U.~Trottenberg, C.~W. Oosterlee, A.~Schuller, Multigrid, Elsevier, 2000.

\bibitem{toselli2004domain}
A.~Toselli, O.~Widlund, Domain decomposition methods-algorithms and theory, Vol.~34, Springer Science \& Business Media, 2004.

\bibitem{notay2000flexible}
Y.~Notay, Flexible conjugate gradients, SIAM Journal on Scientific Computing 22~(4) (2000) 1444--1460.

\bibitem{saad1993flexible}
Y.~Saad, A flexible inner-outer preconditioned gmres algorithm, SIAM Journal on Scientific Computing 14~(2) (1993) 461--469.

\bibitem{cuomo2022scientific}
S.~Cuomo, V.~S. Di~Cola, F.~Giampaolo, G.~Rozza, M.~Raissi, F.~Piccialli, Scientific machine learning through physics--informed neural networks: Where we are and what’s next, Journal of Scientific Computing 92~(3) (2022) 88.

\bibitem{kovachki2023neural}
N.~Kovachki, Z.~Li, B.~Liu, K.~Azizzadenesheli, K.~Bhattacharya, A.~Stuart, A.~Anandkumar, Neural operator: Learning maps between function spaces with applications to pdes, Journal of Machine Learning Research 24~(89) (2023) 1--97.

\bibitem{cui2024neural}
C.~Cui, K.~Jiang, S.~Shu, A neural multigrid solver for helmholtz equations with high wavenumber and heterogeneous media, SIAM Journal on Scientific Computing 47~(3) (2025) C655--C679.

\bibitem{kaneda2022deep}
A.~Kaneda, O.~Akar, J.~Chen, V.~A.~T. Kala, D.~Hyde, J.~Teran, A deep conjugate direction method for iteratively solving linear systems, Proceedings of the 40th International Conference on Machine Learning 202 (2023) 15720--15736.

\bibitem{trifonov2024learning}
V.~Trifonov, A.~Rudikov, O.~Iliev, I.~Oseledets, E.~Muravleva, Learning from linear algebra: A graph neural network approach to preconditioner design for conjugate gradient solvers, arXiv preprint arXiv:2405.15557 (2024).

\bibitem{katrutsa2020black}
A.~Katrutsa, T.~Daulbaev, I.~Oseledets, Black-box learning of multigrid parameters, Journal of Computational and Applied Mathematics 368 (2020) 112524.

\bibitem{huang2022learning}
R.~Huang, R.~Li, Y.~Xi, Learning optimal multigrid smoothers via neural networks, SIAM Journal on Scientific Computing~(0) (2022) S199--S225.

\bibitem{chen2022meta}
Y.~Chen, B.~Dong, J.~Xu, Meta-mgnet: Meta multigrid networks for solving parameterized partial differential equations, Journal of Computational Physics 455 (2022) 110996.

\bibitem{greenfeld2019learning}
D.~Greenfeld, M.~Galun, R.~Basri, I.~Yavneh, R.~Kimmel, Learning to optimize multigrid pde solvers, in: International Conference on Machine Learning, PMLR, 2019, pp. 2415--2423.

\bibitem{luz2020learning}
I.~Luz, M.~Galun, H.~Maron, R.~Basri, I.~Yavneh, Learning algebraic multigrid using graph neural networks, in: International Conference on Machine Learning, PMLR, 2020, pp. 6489--6499.

\bibitem{kopanivcakova2024deeponet}
A.~Kopani{\v{c}}{\'a}kov{\'a}, G.~E. Karniadakis, Deeponet based preconditioning strategies for solving parametric linear systems of equations, SIAM Journal on Scientific Computing 47~(1) (2025) C151--C181.

\bibitem{taghibakhshi2021optimization}
A.~Taghibakhshi, S.~MacLachlan, L.~Olson, M.~West, Optimization-based algebraic multigrid coarsening using reinforcement learning, Advances in Neural Information Processing Systems 34 (2021) 12129--12140.

\bibitem{caldana2023deep}
M.~Caldana, P.~F. Antonietti, et~al., A deep learning algorithm to accelerate algebraic multigrid methods in finite element solvers of 3d elliptic pdes, Computers \& Mathematics with Applications 167 (2024) 217--231.

\bibitem{zou2023autoamg}
H.~Zou, X.~Xu, C.-S. Zhang, Z.~Mo, Autoamg ($\theta$): An auto-tuned amg method based on deep learning for strong threshold, Communications in Computational Physics 36~(1) (2024) 200--220.

\bibitem{klawonn2024learning}
A.~Klawonn, M.~Lanser, J.~Weber, Learning adaptive coarse basis functions of feti-dp, Journal of Computational Physics 496 (2024) 112587.

\bibitem{taghibakhshi2022learning}
A.~Taghibakhshi, N.~Nytko, T.~Zaman, S.~MacLachlan, L.~Olson, M.~West, Learning interface conditions in domain decomposition solvers, Advances in Neural Information Processing Systems 35 (2022) 7222--7235.

\bibitem{knoke2023domain}
T.~Knoke, S.~Kinnewig, S.~Beuchler, A.~Demircan, U.~Morgner, T.~Wick, Domain decomposition with neural network interface approximations for time-harmonic maxwell's equations with different wave numbers, arXiv preprint arXiv:2303.02590 (2023).

\bibitem{taghibakhshi2023mg}
A.~Taghibakhshi, N.~Nytko, T.~U. Zaman, S.~MacLachlan, L.~Olson, M.~West, Mg-gnn: Multigrid graph neural networks for learning multilevel domain decomposition methods, Proceedings of the 40th International Conference on Machine Learning 202 (2023) 33381--33395.

\bibitem{rahaman2019spectral}
N.~Rahaman, A.~Baratin, D.~Arpit, F.~Draxler, M.~Lin, F.~Hamprecht, Y.~Bengio, A.~Courville, On the spectral bias of neural networks, in: International Conference on Machine Learning, PMLR, 2019, pp. 5301--5310.

\bibitem{hong2022activation}
Q.~Hong, Q.~Tan, J.~W. Siegel, J.~Xu, On the activation function dependence of the spectral bias of neural networks, arXiv preprint arXiv:2208.04924 (2022).

\bibitem{xu2019frequency}
Z.-Q.~J. Xu, Y.~Zhang, T.~Luo, Y.~Xiao, Z.~Ma, Frequency principle: Fourier analysis sheds light on deep neural networks, Communications in Computational Physics 28~(5) (2020) 1746--1767.

\bibitem{zhang2022hybrid}
E.~Zhang, A.~Kahana, A.~Kopani{\v{c}}{\'a}kov{\'a}, E.~Turkel, R.~Ranade, J.~Pathak, G.~E. Karniadakis, Blending neural operators and relaxation methods in pde numerical solvers, Nature Machine Intelligence (2024) 1--11.

\bibitem{lu2019deeponet}
L.~Lu, P.~Jin, G.~Pang, Z.~Zhang, G.~E. Karniadakis, Learning nonlinear operators via deeponet based on the universal approximation theorem of operators, Nature machine intelligence 3~(3) (2021) 218--229.

\bibitem{kahana2023geometry}
A.~Kahana, E.~Zhang, S.~Goswami, G.~Karniadakis, R.~Ranade, J.~Pathak, On the geometry transferability of the hybrid iterative numerical solver for differential equations, Computational Mechanics 72~(3) (2023) 471--484.

\bibitem{zou2024large}
Z.~Zou, A.~Kahana, E.~Zhang, E.~Turkel, R.~Ranade, J.~Pathak, G.~E. Karniadakis, Large scale scattering using fast solvers based on neural operators, arXiv preprint arXiv:2405.12380 (2024).

\bibitem{cui2022fourier}
C.~Cui, K.~Jiang, Y.~Liu, S.~Shu, Fourier neural solver for large sparse linear algebraic systems, Mathematics 10~(21) (2022).

\bibitem{brandt1977multi}
A.~Brandt, Multi-level adaptive solutions to boundary-value problems, Mathematics of computation 31~(138) (1977) 333--390.

\bibitem{fortunato2020fast}
D.~Fortunato, A.~Townsend, Fast poisson solvers for spectral methods, IMA Journal of Numerical Analysis 40~(3) (2020) 1994--2018.

\bibitem{strang2007computational}
G.~Strang, Computational science and engineering, SIAM, 2007.

\bibitem{xie2023mgcnn}
Y.~Xie, M.~Lv, C.~Zhang, Mgcnn: a learnable multigrid solver for linear pdes on structured grids, arXiv preprint arXiv:2312.11093 (2023).

\bibitem{rudikov2024neural}
A.~Rudikov, V.~Fanaskov, E.~Muravleva, Y.~M. Laevsky, I.~Oseledets, Neural operators meet conjugate gradients: The fcg-no method for efficient pde solving, arXiv preprint arXiv:2402.05598 (2024).

\bibitem{fanaskov2023spectral}
V.~Fanaskov, I.~V. Oseledets, Spectral neural operators, in: Doklady Mathematics, Vol. 108, Pleiades Publishing Moscow, 2023, pp. S226--S232.

\bibitem{hu2024hybrid}
J.~Hu, P.~Jin, A hybrid iterative method based on mionet for pdes: Theory and numerical examples, arXiv preprint arXiv:2402.07156 (2024).

\bibitem{jin2022mionet}
P.~Jin, S.~Meng, L.~Lu, Mionet: Learning multiple-input operators via tensor product, SIAM Journal on Scientific Computing 44~(6) (2022) A3490--A3514.

\bibitem{chen2024graph}
J.~Chen, Graph neural preconditioners for iterative solutions of sparse linear systems, arXiv preprint arXiv:2406.00809 (2024).

\bibitem{hsieh2019learning}
J.-T. Hsieh, S.~Zhao, S.~Eismann, L.~Mirabella, S.~Ermon, Learning neural pde solvers with convergence guarantees, International Conference on Learning Representations (2019).

\bibitem{bolten2018fourier}
M.~Bolten, H.~Rittich, Fourier analysis of periodic stencils in multigrid methods, SIAM journal on scientific computing 40~(3) (2018) A1642--A1668.

\bibitem{kumar2019local}
P.~Kumar, C.~Rodrigo, F.~J. Gaspar, C.~W. Oosterlee, On local fourier analysis of multigrid methods for pdes with jumping and random coefficients, SIAM Journal on Scientific Computing 41~(3) (2019) A1385--A1413.

\bibitem{brown2019local}
J.~Brown, Y.~He, S.~MacLachlan, Local fourier analysis of balancing domain decomposition by constraints algorithms, SIAM Journal on Scientific Computing 41~(5) (2019) S346--S369.

\bibitem{li2020fourier}
Z.~Li, N.~B. Kovachki, K.~Azizzadenesheli, B.~liu, K.~Bhattacharya, A.~Stuart, A.~Anandkumar, Fourier neural operator for parametric partial differential equations, International Conference on Learning Representations (2021).

\bibitem{hornik1989multilayer}
K.~Hornik, M.~Stinchcombe, H.~White, Multilayer feedforward networks are universal approximators, Neural networks 2~(5) (1989) 359--366.

\bibitem{CSIAM-AM-4-13}
C.~Cui, K.~Jiang, S.~Shu, Solving time-dependent parametric pdes by multiclass classification-based reduced order model, CSIAM Transactions on Applied Mathematics 4~(1) (2023) 13--40.

\bibitem{Paszke2019PyTorchAI}
A.~Paszke, S.~Gross, F.~Massa, A.~Lerer, J.~Bradbury, G.~Chanan, T.~Killeen, Z.~Lin, N.~Gimelshein, L.~Antiga, A.~Desmaison, A.~K{\"o}pf, E.~Yang, Z.~DeVito, M.~Raison, A.~Tejani, S.~Chilamkurthy, B.~Steiner, L.~Fang, J.~Bai, S.~Chintala, Pytorch: An imperative style, high-performance deep learning library, in: Neural Information Processing Systems, 2019.

\bibitem{elman2014finite}
H.~Elman, D.~Silvester, A.~Wathen, Finite elements and fast iterative solvers: with applications in incompressible fluid dynamics (2014).

\bibitem{wesseling1995introduction}
P.~Wesseling, Introduction to multigrid methods, Tech. rep. (1995).

\bibitem{azulay2022multigrid}
Y.~Azulay, E.~Treister, Multigrid-augmented deep learning preconditioners for the helmholtz equation, SIAM Journal on Scientific Computing (2022) S127--S151.

\end{thebibliography}

\end{document}